\documentclass[12pt]{amsart}
\usepackage{mathrsfs}
\usepackage{cases}
\usepackage{epic}
\usepackage{amsfonts}
\usepackage{tikz}
\usepackage{tikz-cd}
\usepackage{graphicx}
\usepackage{amsmath}
\usepackage{amssymb, upgreek}
\usepackage{bm}
\usepackage{latexsym,todonotes}
\usepackage{pdflscape}
\usepackage[all]{xypic}
\usepackage[all]{xy}
\usepackage{tikz-cd}
\usepackage{color}
\usepackage{colordvi}
\usepackage{multicol}
\usepackage{mathtools}
\usepackage[normalem]{ulem}
\usepackage[linktocpage=true]{hyperref}
\hypersetup{colorlinks,linkcolor=blue,urlcolor=cyan,citecolor=blue}
\usepackage{graphicx}
\NewDocumentCommand{\sslash}{s}{%
	\IfBooleanTF{#1}
	{\big/\mkern-7mu\big/}
	{/\mkern-6mu/}%
}
\makeatletter
\newsavebox{\@brx}
\newcommand{\llangle}[1][]{\savebox{\@brx}{\(\m@th{#1\langle}\)}%
	\mathopen{\copy\@brx\kern-0.5\wd\@brx\usebox{\@brx}}}
\newcommand{\rrangle}[1][]{\savebox{\@brx}{\(\m@th{#1\rangle}\)}%
	\mathclose{\copy\@brx\kern-0.5\wd\@brx\usebox{\@brx}}}
\makeatother

\tikzcdset{scale cd/.style={every label/.append style={scale=#1},
		cells={nodes={scale=#1}}}}

\DeclareMathAlphabet{\mathbbb}{U}{bbold}{m}{n}
\DeclareMathOperator{\dimv}{\underline{\dim}}

\makeatletter

\newcommand{\Rmnum}[1]{\textup{\expandafter\@slowromancap\romannumeral #1@}}
\makeatother

\allowdisplaybreaks

\def \tUU{(\tU\otimes\tU)}
\def \tUUi {(\tU\otimes\tU)^\imath}

\begin{document}
	\input xy
	\xyoption{all}
	\newcommand{\iLa}{\Lambda^{\imath}}
	\newcommand{\iadd}{\operatorname{iadd}\nolimits}
	\renewcommand{\mod}{\operatorname{mod}\nolimits}
	\newcommand{\fproj}{\operatorname{f.proj}\nolimits}
	\newcommand{\Fac}{\operatorname{Fac}\nolimits}
	\newcommand{\ci}{{\I}_{\btau}}
	\newcommand{\proj}{\operatorname{proj}\nolimits}
	\newcommand{\inj}{\operatorname{inj}\nolimits}
	\newcommand{\rad}{\operatorname{rad}\nolimits}
	\newcommand{\Span}{\operatorname{Span}\nolimits}
	\newcommand{\soc}{\operatorname{soc}\nolimits}
	\newcommand{\ind}{\operatorname{inj.dim}\nolimits}
	\newcommand{\Ginj}{\operatorname{Ginj}\nolimits}
	\newcommand{\res}{\operatorname{res}\nolimits}
	\newcommand{\np}{\operatorname{np}\nolimits}
	\newcommand{\Mor}{\operatorname{Mor}\nolimits}
	\newcommand{\Mod}{\operatorname{Mod}\nolimits}
	\newcommand{\End}{\operatorname{End}\nolimits}
	\newcommand{\lf}{\operatorname{l.f.}\nolimits}
	\newcommand{\Iso}{\operatorname{Iso}\nolimits}
	\newcommand{\Aut}{\operatorname{Aut}\nolimits}
	\newcommand{\Rep}{\operatorname{Rep}\nolimits}
	
	\newcommand{\colim}{\operatorname{colim}\nolimits}
	\newcommand{\gldim}{\operatorname{gl.dim}\nolimits}
	\newcommand{\cone}{\operatorname{cone}\nolimits}
	\newcommand{\rep}{\operatorname{rep}\nolimits}
	\newcommand{\Ext}{\operatorname{Ext}\nolimits}
	\newcommand{\Tor}{\operatorname{Tor}\nolimits}
	\newcommand{\Hom}{\operatorname{Hom}\nolimits}
	\newcommand{\Top}{\operatorname{top}\nolimits}
	\newcommand{\Coker}{\operatorname{Coker}\nolimits}
	\newcommand{\thick}{\operatorname{thick}\nolimits}
	\newcommand{\rank}{\operatorname{rank}\nolimits}
	\newcommand{\Gproj}{\operatorname{Gproj}\nolimits}
	\newcommand{\Len}{\operatorname{Length}\nolimits}
	\newcommand{\RHom}{\operatorname{RHom}\nolimits}
	\renewcommand{\deg}{\operatorname{deg}\nolimits}
	\renewcommand{\Im}{\operatorname{Im}\nolimits}
	\newcommand{\Ker}{\operatorname{Ker}\nolimits}
	\newcommand{\Coh}{\operatorname{Coh}\nolimits}
	\newcommand{\Id}{\operatorname{Id}\nolimits}
	\newcommand{\Qcoh}{\operatorname{Qch}\nolimits}
	\newcommand{\CM}{\operatorname{CM}\nolimits}
	\newcommand{\sgn}{\operatorname{sgn}\nolimits}
	\newcommand{\utMH}{\operatorname{\cm\ch(\iLa)}\nolimits}
	\newcommand{\GL}{\operatorname{GL}}
	\newcommand{\Perv}{\operatorname{Perv}}
	
	\newcommand{\IC}{\operatorname{IC}}
	\def \hU{\widehat{\U}}
	\def \hUi{\widehat{\U}^\imath}
	\newcommand{\bb}{\psi_*}
	\newcommand{\bvs}{{\boldsymbol{\varsigma}}}
	\def \ba{\mathbf{a}}
	\newcommand{\vs}{\varsigma}
	\def \bfk {\mathbf{k}}

	\def \bd{\mathbf{d}}
	\newcommand{\e}{{\bf 1}}
	\newcommand{\EE}{E^*}
	\newcommand{\dbl}{\operatorname{dbl}\nolimits}
	\newcommand{\ga}{\gamma}
	\newcommand{\tM}{\cm\widetilde{\ch}}
	\newcommand{\la}{\lambda}
	
	\newcommand{\For}{\operatorname{{\bf F}or}\nolimits}
	\newcommand{\coker}{\operatorname{Coker}\nolimits}
	\newcommand{\rankv}{\operatorname{\underline{rank}}\nolimits}
	\newcommand{\diag}{{\operatorname{diag}\nolimits}}
	\newcommand{\swa}{{\operatorname{swap}\nolimits}}
	\newcommand{\supp}{{\operatorname{supp}}}
	
	\renewcommand{\Vec}{{\operatorname{Vec}\nolimits}}
	\newcommand{\pd}{\operatorname{proj.dim}\nolimits}
	\newcommand{\gr}{\operatorname{gr}\nolimits}
	\newcommand{\id}{\operatorname{id}\nolimits}
	\newcommand{\aut}{\operatorname{Aut}\nolimits}
	\newcommand{\Gr}{\operatorname{Gr}\nolimits}
	
	\newcommand{\pdim}{\operatorname{proj.dim}\nolimits}
	\newcommand{\idim}{\operatorname{inj.dim}\nolimits}
	\newcommand{\Gd}{\operatorname{G.dim}\nolimits}
	\newcommand{\Ind}{\operatorname{Ind}\nolimits}
	\newcommand{\add}{\operatorname{add}\nolimits}
	\newcommand{\pr}{\operatorname{pr}\nolimits}
	\newcommand{\oR}{\operatorname{R}\nolimits}
	\newcommand{\oL}{\operatorname{L}\nolimits}
	\def \brW{\mathrm{Br}(W_\btau)}
	\newcommand{\Perf}{{\mathfrak Perf}}
	\newcommand{\cc}{{\mathcal C}}
	\newcommand{\gc}{{\mathcal GC}}
	\newcommand{\ce}{{\mathcal E}}
	\newcommand{\calI}{{\mathcal I}}
	\newcommand{\cs}{{\mathcal S}}
	\newcommand{\cf}{{\mathcal F}}
	\newcommand{\cx}{{\mathcal X}}
	\newcommand{\cy}{{\mathcal Y}}
	\newcommand{\ct}{{\mathcal T}}
	\newcommand{\cu}{{\mathcal U}}
	\newcommand{\cv}{{\mathcal V}}
	\newcommand{\cn}{{\mathcal N}}
	\newcommand{\mcr}{{\mathcal R}}
	\newcommand{\ch}{{\mathcal H}}
	\newcommand{\ca}{{\mathcal A}}
	\newcommand{\cb}{{\mathcal B}}
	\newcommand{\cj}{{\mathcal J}}
	\newcommand{\cl}{{\mathcal L}}
	\newcommand{\cm}{{\mathcal M}}
	\newcommand{\cp}{{\mathcal P}}
	\newcommand{\cg}{{\mathcal G}}
	\newcommand{\cw}{{\mathcal W}}
	\newcommand{\co}{{\mathcal O}}
	\newcommand{\cq}{{\mathcal Q}}
	\newcommand{\cd}{{\mathcal D}}
	\newcommand{\ck}{{\mathcal K}}
	\newcommand{\calr}{{\mathcal R}}
	\newcommand{\cz}{{\mathcal Z}}
	\newcommand{\ol}{\overline}
	\newcommand{\ul}{\underline}
	\newcommand{\st}{[1]}
	\newcommand{\ow}{\widetilde}
	\renewcommand{\P}{\mathbf{P}}
	\newcommand{\pic}{\operatorname{Pic}\nolimits}
	\newcommand{\Spec}{\operatorname{Spec}\nolimits}
	\newcommand{\Fr}{\mathrm{Fr}}
	\newcommand{\Gp}{\mathrm{Gp}}

	%Theorem for the introduciton
	\newtheorem{innercustomthm}{{\bf Theorem}}
	\newenvironment{customthm}[1]
	{\renewcommand\theinnercustomthm{#1}\innercustomthm}
	{\endinnercustomthm}
	
	\newtheorem{innercustomcor}{{\bf Corollary}}
	\newenvironment{customcor}[1]
	{\renewcommand\theinnercustomcor{#1}\innercustomcor}
	{\endinnercustomthm}
	
	\newtheorem{innercustomprop}{{\bf Proposition}}
	\newenvironment{customprop}[1]
	{\renewcommand\theinnercustomprop{#1}\innercustomprop}
	{\endinnercustomthm}
	
	\newtheorem{theorem}{Theorem}[section]
	\newtheorem{acknowledgement}[theorem]{Acknowledgement}
	\newtheorem{algorithm}[theorem]{Algorithm}
	\newtheorem{axiom}[theorem]{Axiom}
	\newtheorem{case}[theorem]{Case}
	\newtheorem{claim}[theorem]{Claim}
	\newtheorem{conclusion}[theorem]{Conclusion}
	\newtheorem{condition}[theorem]{Condition}
	\newtheorem{conjecture}[theorem]{Conjecture}
	\newtheorem{construction}[theorem]{Construction}
	\newtheorem{corollary}[theorem]{Corollary}
	\newtheorem{criterion}[theorem]{Criterion}
	\newtheorem{definition}[theorem]{Definition}
	\newtheorem{example}[theorem]{Example}
	\newtheorem{assumption}[theorem]{Assumption}
	\newtheorem{lemma}[theorem]{Lemma}
	\newtheorem{notation}[theorem]{Notation}
	\newtheorem{problem}[theorem]{Problem}
	\newtheorem{proposition}[theorem]{Proposition}
	\newtheorem{solution}[theorem]{Solution}
	\newtheorem{summary}[theorem]{Summary}
	\newtheorem{hypothesis}[theorem]{Hypothesis}
	\newtheorem*{thm}{Theorem}
	
	\theoremstyle{remark}
	\newtheorem{remark}[theorem]{Remark}
	
	\def \Br{\mathrm{Br}}
	\newcommand{\tK}{K}
	
	\newcommand{\tk}{\widetilde{k}}
	\newcommand{\tU}{\widetilde{{\mathbf U}}}
	\newcommand{\Ui}{{\mathbf U}^\imath}
	\newcommand{\tUi}{\widetilde{{\mathbf U}}^\imath}
	\newcommand{\qbinom}[2]{\begin{bmatrix} #1\\#2 \end{bmatrix} }
	\newcommand{\ov}{\overline}
	\newcommand{\tMHg}{\operatorname{\widetilde{\ch}(Q,\btau)}\nolimits}
	\newcommand{\tMHgop}{\operatorname{\widetilde{\ch}(Q^{op},\btau)}\nolimits}
	
	\newcommand{\rMHg}{\operatorname{\ch_{\rm{red}}(Q,\btau)}\nolimits}
	\newcommand{\dg}{\operatorname{dg}\nolimits}
	\def \fu{{\mathfrak{u}}}
	\def \fv{{\mathfrak{v}}}
	\def \sqq{{\mathbbb{v}}}
	\def \bp{{\mathbf p}}
	\def \bv{{\mathbf v}}
	\def \bw{{\mathbf w}}
	\def \bA{{\mathbf A}}
	\def \bL{{\mathbf L}}
	\def \bF{{\mathbf F}}
	\def \bS{{\mathbf S}}
	\def \bC{{\mathbf C}}
	\def \bU{{\mathbf U}}
	\def \U{{\mathbf U}}
	\def \btau{\varpi}
	\def \La{\Lambda}
	\def \Res{\Delta}
	\newcommand{\ev}{\bar{0}}
	\newcommand{\odd}{\bar{1}}
	\def \fk{\mathfrak{k}}
	\def \ff{\mathfrak{f}}
	\def \fp{{\mathfrak{P}}}
	\def \fg{\mathfrak{g}}
	\def \fn{\mathfrak{n}}
	\def \gr{\mathfrak{gr}}
	\def \Z{\mathbb{Z}}
	\def \F{\mathbb{F}}
	\def \D{\mathbb{D}}
	\def \C{\mathbb{C}}
	\def \N{\mathbb{N}}
	\def \Q{\mathbb{Q}}
	\def \G{\mathbb{G}}
	\def \P{\mathbb{P}}
	\def \K{\mathbb{K}}
	\def \E{\mathbb{K}}
	\def \I{\mathbb{I}}
	
	\def \eps{\varepsilon}
	\def \BH{\mathbb{H}}
	\def \btau{\varrho}
	\def \cv{\varpi}
	
	\def \tR{\widetilde{\bf R}}
	\def \tRZ{\widetilde{\bf R}_\cz}
	\def \hR{\widehat{\bf R}}
	\def \hRZ{\widehat{\bf R}_\cz}
	\def\tRi{\widetilde{\bf R}^\imath}
	\def\hRi{\widehat{\bf R}^\imath}
	\def\tRiZ{\widetilde{\bf R}^\imath_\cz}
	\def\reg{\mathrm{reg}}
	
	\def\hRiZ{\widehat{\bf R}^\imath_\cz}
	\def \tTT{\widetilde{\mathbf{T}}}
	\def \TT{\mathbf{T}}
	\def \br{\mathbf{r}}
	\def \bp{{\mathbf p}}
	\def \tS{\texttt{S}}
	\def \bq{{\bm q}}
	\def \bvt{{v}}
	\def \bs{{ r}}
	\def \tt{{v}}
	\def \k{k}
	\def \bnu{\bm{\nu}}
	\def\bc{\mathbf{c}}
	\def \ts{\textup{\texttt{s}}}
	\def \tt{\textup{\texttt{t}}}
	\def \tr{\textup{\texttt{r}}}
	\def \tc{\textup{\texttt{c}}}
	\def \tg{\textup{\texttt{g}}}
	\def \bW{\mathbf{W}}
	\def \bV{\mathbf{V}}

	\newcommand{\browntext}[1]{\textcolor{brown}{#1}}
	\newcommand{\greentext}[1]{\textcolor{green}{#1}}
	\newcommand{\redtext}[1]{\textcolor{red}{#1}}
	\newcommand{\bluetext}[1]{\textcolor{blue}{#1}}
	\newcommand{\brown}[1]{\browntext{ #1}}
	\newcommand{\green}[1]{\greentext{ #1}}
	\newcommand{\red}[1]{\redtext{ #1}}
	\newcommand{\blue}[1]{\bluetext{ #1}}
	\numberwithin{equation}{section}
	\renewcommand{\theequation}{\thesection.\arabic{equation}}
	
	%todo
	\newcommand{\wtodo}{\rightarrowdo[inline,color=orange!20, caption={}]}
	\newcommand{\lutodo}{\rightarrowdo[inline,color=green!20, caption={}]}
	\def \tT{\widetilde{\mathcal T}}
	
	\def \tTL{\tT(\iLa)}
	\def \iH{\widetilde{\ch}}
	
	%\title[Dual canonical bases arising from quantum symmetric pairs]{Dual canonical bases arising from quantum symmetric pairs}
	
	\title[Dual canonical bases of quantum groups and $\imath$quantum groups I]{Dual canonical bases of quantum groups and $\imath$quantum groups I: Hall algebras}
	
	\author[Ming Lu]{Ming Lu}
	\address{Department of Mathematics, Sichuan University, Chengdu 610064, P.R.China}
	\email{luming@scu.edu.cn}

	\author[Xiaolong Pan]{Xiaolong Pan}
	\address{Department of Mathematics, Sichuan University, Chengdu 610064, P.R.China}
	\email{xiaolong\_pan@stu.scu.edu.cn}


	\begin{abstract}
		
		The $\imath$quantum groups have two realizations: one via the $\imath$Hall algebras and the other via the quantum Grothendieck rings of quiver varieties, as developed by the first author and Wang. The isoclasses of perverse sheaves provide the dual canonical bases for $\imath$quantum groups of type ADE with integral and positive structure constants. In this paper, we present a new construction of the dual canonical bases in the setting of $\imath$Hall algebras. We also introduce Fourier transforms for $\imath$Hall algebras, and prove the invariance of the dual canonical bases under braid group actions and Fourier transforms. Since quantum groups are $\imath$quantum groups of diagonal type, all results also apply to this classical case.
	\end{abstract}

	\maketitle
	\setcounter{tocdepth}{1}
	
	\tableofcontents
	
	%%%%%%%
	%%%%%%%
	\section{Introduction}

	%%%%%%%%
	\subsection{Background}

	Bridgeland \cite{Br13} constructed a Hall algebra from $\Z_2$-graded complexes of projective representations of a quiver to realize the Drinfeld double version $\tU$ of a quantum group $\U$, which was built on the Hall algebra realization of halves of quantum groups $\U^+$ in \cite{Rin90,Gr95}. 
	Inspired by Bridgeland's construction, Gorsky \cite{Gor18} introduced semi-derived Hall algebras for Frobenius categories while the first author and Peng \cite{LP21} formulated the  semi-derived Ringel-Hall algebras for the cateogries of $\Z_2$-graded complexes of hereditary abelain categories. In \cite{LW19}, the semi-derived Ringel-Hall algebras are further developed for 1-Gorenstein algebras. 
	The Hall algebra realization %\cite{Rin90,Gr95,Br13} 
	has led to fruitful ways of understanding the quantum groups, such as braid group actions and PBW bases; see \cite{Rin3, SV99,  LW21a}. 
	%. The braid group actions of $\U$ by Lusztig can also be understood from a categorical viewpoint via the BGP reflection functors.  
	%The reflection functors were adapted for the Hall algebra setting by Ringel  \cite{Rin3} to provide isomorphisms of subalgebras of $\U^+$. 
	%The reflection functors are also used to construct PBW bases of $\U^+$ \cite{Lus93} in the framework of Hall algebras \cite{Rin3}. Later the reflection functors were extended in \cite{LW21a} to automorphisms of semi-derived Ringel-Hall algebras, and they were identified with Lusztig's braid group actions on $\U$; also see  \cite{SV99, XY01} for braid group actions via automorphisms of Drinfeld double of Hall algebras. %Using the Hall algebras, one can apply BGP-reflection functors to construct Lusztig’s braid group actions and PBW bases of quantum groups; see \cite{Rin3, Lus98,SV99,XY01}. We refer to \cite{DDPW}.

	%, and Bridgeland's Theorem \cite{Br13} is reformulated by using the semi-derived Ringel-Hall algebras of framed double quiver algebras. 
	The universal $\imath$quantum group $\tUi =\langle B_i, \tk_i \mid i\in \I \rangle$ is by definition a subalgebra of $\tU$ associated to a Satake diagram, and $(\tU, \tUi)$ forms a quantum symmetric pair; cf. \cite{LW19}.
	%As a quantization of symmetric pairs $(\fg, \fg^\theta)$, the quantum symmetric pairs $(\U, \Ui)$ were formulated by Letzter \cite{Let99} (also cf. \cite{Ko14}) with Satake diagrams as inputs. 
	A central reduction of $\tUi$ produces Letzter's  $\imath$quantum group $\Ui =\Ui_\bvs$ with parameters $\bvs$ \cite{Let99,Ko14}. A breakthrough of $\imath$quantum groups is the discovery of $\imath$canonical bases \cite{BW18,BW18b}. In this paper, we focus on universal $\imath$quantum groups $\tUi$ of quasi-split type (Satake diagrams with only white nodes) and Drinfeld double version of quantum groups $\tU$. We view Letzter's $\imath$quantum groups as a vast generalization of Drinfeld--Jimbo quantum groups, and quantum groups can be viewed as $\imath$quantum groups of diagonal type (see Example \ref{ex:QGvsiQG}).

	%The first author and Wang \cite{LW19,LW20} used the $\imath$Hall algebra (i.e., the twisted semi-derived Ringel-Hall algebra) of an $\imath$quiver  algebra $\Lambda^\imath$ associated to an $\imath$quiver $(Q,\btau)$ to realize the  $\imath$quantum groups $\tUi$. 

	An $\imath$quiver $(Q, \btau)$ by definition consists of a  quiver $Q =(Q_0, Q_1)$ together with an involution $\btau$ on $Q$ (here $\btau =\Id$ is allowed).
	Associated to an $\imath$quiver $(Q, \btau)$, the first author and Wang  \cite{LW19} recently formulated a family of finite-dimensional 1-Gorenstein algebras, called $\imath$quiver algebras and denoted by $\iLa$. The semi-derived Ringel-Hall algebras for $\iLa$ (or simply the $\imath$Hall algebras) were then shown to provide a categorical realization of the quasi-split  universal $\imath$quantum groups $\tUi$; in the special case of $\imath$quivers of diagonal type, this construction reduces to a reformulation of Bridgeland's Hall algebra realization of quantum groups $\tU$. %the Drinfeld double $\tU$ is reproduced by using $\imath$Hall algebra associated to the $\imath$quiver of diagonal type. 
	In \cite{LW21b}, %they introduce a BGP type reflection functor for $\imath$quiver algebras $\Lambda^\imath$, and then 
	the $\imath$Hall algebra approach also provides a conceptual construction of braid group actions and  PBW bases on $\tUi$; see also \cite{KP11,WZ23} for the algebraic construction of the braid group actions on $\imath$quantum groups.

	%%%%%%%%%%%%
	%\subsection{Background on quiver varieties}
	
	Lusztig in \cite{Lus90,Lus91,Lus93} used perverse sheaves on the variety of representations of a quiver $Q$ to realize $\U^+$, and produced the canonical basis of $\U^+$ by simple perverse sheaves; cf. also \cite{Ka91} for the crystal basis. For type ADE, Lusztig \cite{Lus90} also used PBW bases to give an elementary construction of canonical bases of $\U^+$; also see \cite[\S11.6]{DDPW}. 
	The relation between Ringel’s realization and Lusztig’s categorification is given by sheaf-function correspondence; see \cite{Lus98}. 
	Nakajima \cite{Na01,Na04} further developed Lusztig's construction, and introduced (graded) Nakajima quiver varieties.
	
	Motivated by the works \cite{Na01,HL15,LeP13},
	%Hernander-Leclerc \cite[\S9]{HL15} and Leclerc-Plamondon \cite{LeP13}, 
	Keller and Scherotzke \cite{KS16,Sch19} formulated the notion of regular/singular Nakajima categories $\mcr, \cs$ from the mesh category of the repetition category $\Z Q$ for an acyclic quiver $Q$, called Nakajima-Keller-Scherotzke (NKS) categories in \cite{LW21b}. The NKS varieties are by definition the representation varieties of $\mcr$ and $\cs$. Via (dual) quantum Grothendieck rings of NKS varieties $\cs^\imath$, $\mcr^\imath$ of $\imath$quiver varieties, a geometric construction of $\tUi$ of type ADE is given in \cite{LW21b}. This result recovers Qin's construction \cite{Qin} of $\tU$ via such cyclic quiver varieties, since cyclic quiver varieties can be viewed as NKS varieties of $\cs^\imath$ of diagonal type; cf. \cite{LW19}.

	%%%%%%%%%%%%%%%
	\subsection{Goal}
	The geometric construction in \cite{Qin,LW21b} provides a favorable basis for $\tU$ and  $\tUi$  with positive integral structure constants, which are called {\em dual canonical basis}.
	The goal of this series of papers is to study in depth from both Hall algebra and quiver variety the dual canonical basis for the quantum group $\tU$ and the $\imath$quantum group $\tUi$ of type ADE. Unlike Lusztig's canonical basis for $\U^+$, the dual canonical basis is defined for the entire $\tU$ and $\tUi$; compare with canonical bases of modified quantum groups $\dot{\U}$ and $\imath$quantum groups $\dot{\U}^\imath$ (see \cite{Lus93,BW18b}).
	
	In this paper, to make the dual canonical bases tractable, we % use Lusztig's Lemma and a bar involution (different to Lusztig's bar involution) to
	construct the dual canonical basis for the $\imath$Hall algebra by using a rescaled Hall basis (natural basis of $\imath$Hall algebra).
	We prove that the dual canonical basis does not depend on the orientation of $\imath$quivers, and is invariant under braid group actions of $\tUi$.  For this purpose, we construct Fourier transforms for $\imath$Hall algebras.

	%%%%%%%%%%%%
	\subsection{Main results}
	From a Cartan matrix $C=(c_{ij})_{i,j\in\I}$, we can associate a variant of Drinfeld double quantum group $\hU$, the $\Q(v^{1/2})$-algebra generated by
	$E_i,F_i,K_i,K_i'$ ($i\in\I$)  subject to \eqref{eq:KK}--\eqref{eq:serre2}. The Drinfeld double $\tU$ is the $\Q(v^{1/2})$-algebra constructed from $\hU$ by making $K_i,K_i'$ ($i\in\I$) invertible. Given an involution $\btau\in\aut(C)$, we can define the (universal) $\imath$quantum group $\tUi$ to be the subalgebra of $\tU$ generated by $B_i= F_i +  E_{\btau i} \tK_i',
	\tk_i = \tK_i \tK_{\btau i}'$ ($i \in \I$), and the inverses of $\tk_i$. A variant of $\imath$quantum group $\hUi$ is defined to be subalgebra of $\hU$ similarly, but $\tk_i$ ($i\in\I$) are not invertible. Following \cite{BG17}, the bar-involution $\ov{\cdot}$ of $\tU$ (and $\hU$) is the $\Q$-anti-involution which fixes the generators $E_i,F_i,K_i,K_i'$, and $\ov{v^{1/2}}=v^{-1/2}$ (which is different to Lusztig's bar-involution). We also define the bar-involution $\ov{\cdot}$ on $\tUi$ (and $\hUi$) given by $\ov{v^{1/2}}=v^{-1/2}$, $\ov{B_i}=B_i$, $\ov{\K_i}=\K_i$, for $i\in\I$. Here, $\K_i:=-v\tk_i$, if $\varrho i=i$; $\K_j:=\tk_j$, otherwise.
	
	%%%%
	%\subsubsection{Dual canonical bases via $\imath$Hall algebras}
	
	For a Dynkin $\imath$quiver $(Q,\varrho)$, we can associate a generic $\imath$Hall algebra $\widetilde{\ch}(Q,\varrho)$ defined over $\Q(v^{1/2})$ \cite{LW19}. This algebra $\widetilde{\ch}(Q,\varrho)$ has a basis (called Hall basis) 
	\begin{align}
		\label{eq:iHallbasis}\{\K_\alpha*\fu_\lambda\mid \alpha\in\Z^\I,\lambda\in\fp\},
	\end{align}
	where $\fp$ is the set of $\N$-valued functions over the set $\Phi^+$ of positive roots, and $*$ is the (twisted) Hall product. We also define $\widehat{\ch}(Q,\btau)$ to be the generic $\imath$Hall algebra with the Hall basis $\{\K_\alpha*\fu_\lambda\mid \alpha\in\N^\I,\lambda\in\fp\}$. Let $\widehat{\ct}(Q,\varrho)$ be the subalgebra generated by $\K_\alpha$, $\alpha\in\N^\I$. Inspired by \cite{BG17}, we define an action $\diamond$ of $\widehat{\ct}(Q,\varrho)$ on $\widehat{\ch}(Q,\varrho)$; see \S\ref{subsec:dCB of iHA}. This construction can be carried out to $\widetilde{\ch}(Q,\btau)$. 
	Let $\cz:=\Z[v^{1/2},v^{-1/2}]$. The integral form of the $\imath$Hall algebra $\widetilde{\ch}(Q,\varrho)_\cz$ is the $\cz$-algebra defined over the (free) $\cz$-module spanned by \eqref{eq:iHallbasis}. Similarly for $\widehat{\ch}(Q,\varrho)_\cz$.

	%It is proved in \cite{LW19} that there exist two $\Q(v^{1/2})$-algebra isomorphisms (see \S\ref{subsec:iHall}--\S\ref{sub:generic}):
	%\begin{align}
	%\label{eq:2psi}
	%\widetilde{\psi}:\tUi\longrightarrow \widetilde{\ch}(Q,\btau),\qquad \widehat{\psi}:\hUi\longrightarrow \widehat{\ch}(Q,\btau).
	%\end{align}
	Using the isomorphism between $\imath$quantum groups and $\imath$Hall algebras (see Theorem \ref{lem:Hall-iQG}), the bar-involution of $\tUi$ (and $\hUi$) induces the bar-involution of $\widetilde{\ch}(Q,\btau)$ (and $\widehat{\ch}(Q,\btau)$) which satisfies $\ov{\K_\alpha\diamond \fu_\lambda}=\K_\alpha\diamond\ov{\fu_\lambda}$. The partial order on $\fp$ given by orbit closures \cite{Lus90} induces a partial order $\prec$ on $\N^\I\times\fp$. With this partial order, we can apply Lusztig's Lemma \cite{BZ14} to a rescaled Hall basis $\{\K_\alpha\diamond\mathfrak{U}_\lambda\mid \alpha\in\N^\I,\lambda\in\fp\}$ of $\widehat{\ch}(Q,\btau)$; see \S\ref{subsec:dCB of iHA}. This allows us to construct a bar-invariant integral basis on $\widehat{\ch}(Q,\btau)$. 
	
	\begin{customthm}{{\bf A}} [Theorem~\ref{iHA dCB theorem}]\label{thm B}
		For each $\alpha\in\N^\I$ and $\lambda\in\mathfrak{P}$, there exists a unique element $\mathfrak{L}_{\alpha,\lambda}\in\widehat{\ch}(Q,\varrho)$ such that $\ov{\mathfrak{L}_{\alpha,\lambda}}=\mathfrak{L}_{\alpha,\lambda}$ and
		\[
		\mathfrak{L}_{\alpha,\lambda}-\K_\alpha\diamond \mathfrak{U}_\lambda\in\sum_{(\beta,\mu)}v^{-1}\Z[v^{-1}]\cdot \K_\beta\diamond \mathfrak{U}_\mu.
		\]
		Moreover, $\mathfrak{L}_{\alpha,\lambda}$ satisfies 
		\[
		\mathfrak{L}_{\alpha,\lambda}-\K_\alpha\diamond \mathfrak{U}_\lambda\in\sum_{(\alpha,\lambda)\prec(\beta,\mu)}v^{-1}\Z[v^{-1}]\cdot \K_\beta\diamond \mathfrak{U}_\mu,
		\]
		and $\mathfrak{L}_{\alpha,\lambda}=\K_\alpha\diamond \mathfrak{L}_{0,\lambda}$.
	\end{customthm}
	
	If we write $\mathfrak{L}_\lambda:=\mathfrak{L}_{0,\lambda}$, then $\{\K_\alpha\diamond\mathfrak{L}_\lambda\mid\alpha\in\Z^\I,\lambda\in\fp\}$ is a basis of $\widetilde{\ch}(Q,\varrho)$, called the dual canonical basis. The dual canonical bases give rise to $\cz$-bases of $\widehat{\ch}(Q,\varrho)_\cz$ and $\widetilde{\ch}(Q,\varrho)_\cz$.
	
	For a sink $\ell\in Q_0$, the $\imath$quiver $(r_\ell Q,\btau)$ is constructed from the $\imath$quiver $(Q,\btau)$ by reversing all arrows ending at $\ell$ and $\btau(\ell)$. A BGP type reflection $F_\ell^+:\mod(\Lambda^\imath)\rightarrow \mod(r_\ell\Lambda^\imath)$ constructed in \cite{LW21a} gives an isomorphism $\Gamma_\ell:\widetilde{\ch}(Q,\btau)\rightarrow \widetilde{\ch}(r_\ell Q,\btau)$, which can be used to realize the braid group action $\tTT_\ell:\tUi\rightarrow \tUi$; see \eqref{eq:defT}. %with the help of the isomorphisms $\widetilde{\psi}_Q:\tUi\rightarrow \widetilde{\ch}(Q,\varrho)$ and  $\widetilde{\psi}_{r_\ell Q}:\tUi\rightarrow \widetilde{\ch}(r_\ell Q,\varrho)$.
	We show in Proposition \ref{Gamma_l on H_lambda} that $\Gamma_\ell$ maps the rescaled Hall basis of $\widetilde{\ch}(Q,\btau)$ to the one of $\widetilde{\ch}(r_\ell Q,\btau)$. Since $\Gamma_\ell$ commutes with the bar-involutions of $\imath$Hall algebras, this implies $\Gamma_\ell$ maps the dual canonical basis of $\widetilde{\ch}(Q,\btau)$ to the one of $\widetilde{\ch}(r_\ell Q,\btau)$;  see Corollary~\ref{iHA reflection functor of dCB}.

	%\subsubsection{Fourier transforms}
	%In order to prove that the Hall algebra $\widetilde{\ch}(\bfk Q)$ is independent of the orientation of $Q$, a Fourier transform of Hall algebras is constructed in \cite{Lus90,SV99}. Correspondingly, by using Fourier-Deligne transform, this can also be defined on the quantum Grothendieck rings of Lusztig's varieties \cite{Lus90}, and it is shown that the canonical basis is preserved by the Fourier transform. 
	
	%By realizing Hall algebra using functions on the representation variety of the quiver $Q$, in \cite{SV99} the authors defined a Fourier transform on Hall algebras of different orientations of the underlying diagram. Using Fourier-Deligne transform, this can be defined on the quantum Grothendieck rings \cite{Lus90}, and it is shown that the canonical basis is preserved. 
	
	%In order to construct Fourier transforms on $\imath$Hall algebras, following \cite{SV99}, we modify the construction in \S\ref{subsec:QV function iHA} of realizing $\widetilde{\ch}(Q,\varrho)$ via functions on $\rep(\bw,\mathcal{S}^\imath)$. 
	In order to prove that the Hall algebra $\widetilde{\ch}(\bfk Q)$ is independent of the orientation of $Q$, a Fourier transform of Hall algebras is constructed in \cite{Lus90,SV99}. Using the same idea we define a Fourier transform for $\imath$Hall algebras.
	
	\begin{customthm}{{\bf B}} [Theorem \ref{thm:FThall}, Corollary \ref{cor:FT-ihall}] \label{thm D}
		For a Dynkin $\imath$quiver $(Q,\varrho)$, if $(Q',\varrho)$ is constructed from $(Q,\varrho)$ by reversing some arrows, then there exist isomorphisms of algebras $\Phi_{Q',Q}:\widehat{\ch}(\bfk Q,\btau)\rightarrow \widehat{\ch}(\bfk Q',\btau)$ and $\Phi_{Q',Q}:\widetilde{\ch}(\bfk Q,\btau)\rightarrow \widetilde{\ch}(\bfk Q',\btau)$.
	\end{customthm}
	%The Fourier transforms descent to integral forms of $\imath$Hall algebras; see Corollary \ref{cor:FT-ihall-integral}.
	These Fourier transforms preserve the integral forms of $\imath$Hall algebras; see Corollary \ref{cor:FT-ihall-integral}. In a sequel \cite{LP26}, we shall prove that 
	the dual canonical bases of $\imath$Hall algebras coincide with those of the dual Grothendieck rings $\tRi$ established in \cite{Qin,LW21b}. The fourier transforms can be transferred to quantum Grothendieck rings $\tRi$ and $'\tRi$ of $(Q,\btau)$ and $(Q',\btau)$ respectively, that is, there is an isomorphism $\Psi=\Psi_{Q',Q}:\tRi\rightarrow {'\tRi}$. In this way, we shall prove in the sequel \cite{LP26} that the dual canonical bases of $\widetilde{\ch}(\bfk Q,\varrho)$ and $\tRi$ are invariant under Fourier transforms.

	The dual canonical basis of $\tUi$ is defined to be the image of $\{\mathfrak{L}_{\alpha,\lambda}\mid \alpha\in\Z^\I,\lambda\in\fp\}$ under the isomorphism given in Theorem~\ref{lem:Hall-iQG}. By the above, this basis is independent of the orientations of the $\imath$quivers. %We have also remarked that braid group actions of $\tUi$ can be realized using reflection functors on the $\imath$Hall algebra $\widetilde{\ch}(Q,\varrho)$. 
	Combining with Corollary~\ref{iHA reflection functor of dCB}, we get the invariance of the dual canonical basis of $\tUi$ under braid group actions. 
	
	\begin{customthm}{{\bf C}} [Theorem~\ref{iQG dCB braid invariant}] \label{thm G}
		The dual canonical basis of $\tUi$ is preserved by the action of the  braid group operator $\tTT_i$.
	\end{customthm}
	
	%%%%%%%
	%\subsection{Perspectives}
	
	In \cite{BG17} Berenstein-Greenstein constructed a double canonical basis for quantum group $\tU$. Since quantum groups can be viewed as $\imath$quantum groups of diagonal type, they also possess dual canonical basis. In a sequel to this paper, we shall show that the double canonical basis coincides with the dual canonical basis for quantum groups of type ADE. As applications, we shall solve several conjectures given in {\em loc. cit.} for type ADE, such as positivity and invariance under Lusztig's braid group actions.  
	
	The dual (double) canonical basis of $\tU$ (and also $\tUi$) of type ADE considered in this paper can be viewed as an ideal extension of Lusztig’s dual canonical basis for $\U^+$.
	We expect this basis (integral, bar-invariant, invariant under the braid group actions) also exists for other types, and moreover it is positive for symmetric  types.  
	%In addition, we shall develop an algebraic approach for the positive dual canonical bases for $\tUi$ of split type, 
	%as an analog of Berenstein-Greenstein's construction of double canonical bases for $\tU$ \cite{BG17}. 
	As $\imath$quantum groups have much more types than quantum groups, even for the quasi-split type, %(Satake diagrams with only white nodes), 
	it is a challenging and important open question to find an algebraic characterization and algorithm for the  dual canonical bases for $\tUi$ of all types (as an analog of Berenstein-Greenstein's construction of double canonical bases for $\tU$ \cite{BG17}).

	%%%%%%%%
	\subsection{Organization}
	
	The paper is organized as follows. In \S\ref{sec:QG and iQG} we review the basics of quantum groups and $\imath$quantum groups associated to a Dynkin ($\imath$)quiver. Here we will use different presentations following \cite{BG17}, and the braid group actions are modified accordingly. This presentation will simplify the realization by means of $\imath$Hall algebras, which is reviewed in \S\ref{sec:iHA}.
	
	In \S\ref{sec:dCB via Hall bases} we define the integral forms of $\imath$Hall algebras. Then we introduce the rescaled Hall bases for $\imath$Hall algebras, from which we construct the dual canonical bases using Lusztig's Lemma. It is also proved that this rescaled Hall basis and the dual canonical basis are preserved under the BGP type reflection functors formulated in \cite{LW21a}.

	Using the function realization, we define in \S\ref{sec:FT of iHA} the Fourier transforms of $\imath$Hall algebras, and prove the invariance of dual canonical basis under these isomorphisms.

	\vspace{2mm}
	
	%%%%
	\noindent{\bf Acknowledgments.}
	We thank Jiepeng Fang, Yixin Lan and Fan Qin for helpful discussions on quiver varieties.
	%ML thanks Liangang Peng and Weiqiang Wang for their guidance and continuing encouragement. 
	We thank Weiqiang Wang for his collaboration in related projects, also for his helpful comments and stimulating discussions. ML is partially supported by the National Natural Science Foundation of China (No. 12171333). 
	
	%%%%%%%%%%%%%
	\section{Quantum groups and $\imath$quantum groups}\label{sec:QG and iQG}

	\subsection{Quantum groups}
	\label{subsec:QG}
	
	%Let $Q=(Q_0,Q_1)$ be a Dynkin quiver with vertex set $Q_0= \I$.
	%Let $n_{ij}$ be the number of edges connecting vertex $i$ and $j$. 
	Let $\I=\{1,\dots,n\}$ be the index set. 
	Let $C=(c_{ij})_{i,j \in \I}$ be the Cartan matrix of of a simply-laced semi-simple Lie algebra $\fg$.
	Let $\Delta^+=\{\alpha_i\mid i\in\I\}$ be the set of simple roots of $\fg$, and denote the root lattice by $\Z^{\I}:=\Z\alpha_1\oplus\cdots\oplus\Z\alpha_n$. Let $\Phi^+$ be the set of positive roots. The simple reflection $s_i:\Z^{\I}\rightarrow\Z^{\I}$ is defined to be $s_i(\alpha_j)=\alpha_j-c_{ij}\alpha_i$, for $i,j\in \I$.
	Denote the Weyl group by $W =\langle s_i\mid i\in \I\rangle$.
	
	%Let $\btau$ be an involution of $Q$. We assume that $c_{i,\btau i}=0$ for all $i$, which always hold  for the {\em Dynkin} $\imath$quivers. We denote by $\bs_{i}$ the following element of order 2 in the Weyl group $W$, i.e.,
	%\begin{align}
	%\label{def:simple reflection}
	%\bs_i= \left\{
	%\begin{array}{ll}
	%s_{i}, & \text{ if } \btau i=i;
	%\\
	%s_is_{\btau i}, & \text{ if } \btau i\neq i.
	%\end{array}
	%\right.
	%\end{align}
	%It is well known (cf., e.g., \cite{KP11}) that the {\rm restricted Weyl group} associated to the quasi-split symmetric pair $(\fg, \fg^\theta)$ can be identified with the following subgroup $W_\btau$ of $W$:
	%\begin{align}
	%  \label{eq:Wtau}
	%W_{\btau} =\{w\in W\mid \btau w =w \btau\}
	%\end{align}
	%where $\btau$ is regarded as an automorphism of $\Aut(C)$. Moreover, it admits the following property.

	Let $v$ be an indeterminate. %Write $[A, B]=AB-BA$, and $[A,B]_v=AB-vBA$. 
	Denote, for $r,m \in \N$,
	\[
	[r]=\frac{v^r-v^{-r}}{v-v^{-1}},
	\quad
	[r]!=\prod_{i=1}^r [i], \quad \qbinom{m}{r} =\frac{[m][m-1]\ldots [m-r+1]}{[r]!}.
	\]
	Following \cite{BG17}, the Drinfeld double $\hU := \hU_v(\fg)$ is defined to be the $\Q(v^{1/2})$-algebra generated by $E_i,F_i, \tK_i,\tK_i'$, $i\in \I$, %\red{where $\tK_i, \tK_i'$ are invertible,?}
	subject to the following relations:  for $i, j \in \I$,
	\begin{align}
		[E_i,F_j]= \delta_{ij}(v^{-1}-v) (\tK_i-\tK_i'),  &\qquad [\tK_i,\tK_j]=[\tK_i,\tK_j']  =[\tK_i',\tK_j']=0,
		\label{eq:KK}
		\\
		\tK_i E_j=v^{c_{ij}} E_j \tK_i, & \qquad \tK_i F_j=v^{-c_{ij}} F_j \tK_i,
		\label{eq:EK}
		\\
		\tK_i' E_j=v^{-c_{ij}} E_j \tK_i', & \qquad \tK_i' F_j=v^{c_{ij}} F_j \tK_i',
		\label{eq:K2}
	\end{align}
	and for $i\neq j \in \I$,
	\begin{align}
		& \sum_{r=0}^{1-c_{ij}} (-1)^r \left[ \begin{array}{c} 1-c_{ij} \\r \end{array} \right]  E_i^r E_j  E_i^{1-c_{ij}-r}=0,
		\label{eq:serre1} \\
		& \sum_{r=0}^{1-c_{ij}} (-1)^r \left[ \begin{array}{c} 1-c_{ij} \\r \end{array} \right]  F_i^r F_j  F_i^{1-c_{ij}-r}=0.
		\label{eq:serre2}
	\end{align}

	We define $\tU=\tU_v(\fg)$ to be the $\Q(v^{1/2})$-algebra constructed from $\widehat{\U}$ by making $\tK_i,\tK_i'$ ($i\in\I$) invertible. Then $\tU$ and $\hU$ are $\Z^\I$-graded by setting $\deg E_i=\alpha_i$, $\deg F_i=-\alpha_i$, $\deg K_i=0=\deg K_i'$. %Moreover,
	%let $\tU_\mu$ be the homogeneous subspace of degree $\mu$.  $\tU=\oplus_{\mu\in\Z^\I} \tU_\mu$ and $\hU=\oplus_{\mu\in\Z^\I} \hU_\mu$. 
	
	\begin{remark}
		Note that $\tU$ defined here is different but isomorphic to the one given in \cite{LW19}. 
		In fact, denote 
		\begin{align}
			\label{eq:Udj-gen}
			\ce_i=\frac{E_i}{v^{-1}-v},\qquad \cf_i=\frac{F_i}{v-v^{-1}},\qquad \forall i\in\I.
		\end{align} 
		Then the presentation of $\tU$ given in \cite{LW19} is generated by $\ce_i,\cf_i,K_i,K_i'$ ($i\in\I$).
	\end{remark}

	The quantum group $\bU$ is defined to be quotient algebra of $\hU$ (also $\tU$) modulo the ideal generated by $K_iK_i'-1$ ($i\in\I$). In fact, $\U$ 
	is 
	the $\Q(v^{1/2})$-algebra generated by $E_i,F_i, K_i, K_i^{-1}$, $i\in \I$, subject to the  relations modified from \eqref{eq:KK}--\eqref{eq:serre2} with $\tK_i'$ replaced by $K_i^{-1}$. %The comultiplication $\Delta$ is obtained by modifying \eqref{eq:Delta} with $\tK_i'$ replaced by  $K_i^{-1}$. 
	
	Let $\hU^+$ be the subalgebra of $\hU$ generated by $E_i$ $(i\in \I)$, $\hU^0$ be the subalgebra of $\widehat{\bU}$ generated by $\tK_i, \tK_i'$ $(i\in \I)$, and $\hU^-$ be the subalgebra of $\widehat{\bU}$ generated by $F_i$ $(i\in \I)$, respectively.
	The subalgebras $\tU^+$, $\tU^0$ and $\tU^-$ of $\tU$, and subalgebras $\bU^+$, $\bU^0$ and $\bU^-$ of $\bU$ are defined similarly. Then the algebras $\hU$, $\widetilde{\bU}$ and $\bU$ have triangular decompositions:
	\begin{align*}
		\hU=\hU^+\otimes\hU^0\otimes \hU^-,\qquad 
		\widetilde{\bU} =\widetilde{\bU}^+\otimes \widetilde{\bU}^0\otimes\widetilde{\bU}^-,
		\qquad
		\bU &=\bU^+\otimes \bU^0\otimes\bU^-.
	\end{align*}
	Clearly, ${\bU}^+\cong\widetilde{\bU}^+\cong\hU^+$, ${\bU}^-\cong \widetilde{\bU}^-\cong\hU^-$, $\hU^0\cong\Q(v^{1/2})[\tK_i,\tK_i'\mid i\in\I]$, $\tU^0\cong \Q(v^{1/2})[\tK_i^{\pm1},(\tK_i')^{\pm1}\mid i\in\I]$ and $\bU^0\cong\Q(v^{1/2})[K_i^{\pm1}\mid i\in\I]$. Note that 
	${\bU}^0 \cong \hU^0/(\tK_i \tK_i' -1 \mid   i\in \I)$. 
	For any $\mu=\sum_{i\in\I}m_i\alpha_i\in\Z^\I$, we denote $K_\mu=\prod_{i\in\I} K_i^{m_i}$, $K_\mu'=\prod_{i\in\I} (K_i')^{m_i}$ in $\tU$ (or $\U$); we can view $\tK_\mu,\tK_\mu'$ in $\hU$ if $\mu\in\N^\I$.
	
	\begin{lemma}[cf. \cite{BG17}]\label{QG bar-involution def}
		There exists an anti-involution $u\mapsto \ov{u}$ on $\hU$ (also $\tU$, $\U$) given by $\ov{v^{1/2}}=v^{-1/2}$, $\ov{E_i}=E_i$, $\ov{F_i}=F_i$, and $\ov{K_i}=K_i$, $\ov{K_i'}=K_i'$, for $i\in\I$.
	\end{lemma}
	
	\begin{proof}
		It is enough to check that the bar-involution defined above preserves \eqref{eq:KK}--\eqref{eq:serre2}, which is obvious.
	\end{proof}
	
	\begin{remark}
		This bar-involution is different to the one considered by Lusztig \cite{Lus90}. 
	\end{remark}
	
	Let $\Br(W)$ be the braid group associated to the Weyl group $W$, generated by $t_i$ ($i\in\I$).
	Lusztig introduced braid group symmetries $T_{i,e}',T_{i,e}''$ for $i\in\I$ and $e=\pm1$, on the Drinfeld--Jimbo quantum group $\U$ \cite[\S37.1.3]{Lus93}. These braid group symmetries can be lifted to  the Drinfeld double $\tU$; see \cite[Propositions 6.20–6.21]{LW21a}, which are denoted by $\widetilde{T}_{i,e}',\widetilde{T}_{i,e}''$. We shall use the following modified version of $\widetilde{T}_i:=\widetilde{T}_{i,1}''$ (compare with \cite[Theorem 1.13]{BG17}).
	\begin{proposition}
		%[\cite{Lus90a}]
		\label{prop:braid1}
		There exists an automorphism $\widetilde{T}_{i}$, for $i\in \I$, on $\tU$ such that
		\begin{align*}
			&\widetilde{T}_{i}(K_\mu)= K_{s_i(\mu)},
			\qquad \widetilde{T}_{i}(K'_\mu)= K'_{s_i(\mu)},\;\;\forall \mu\in \Z^\I,\\
			&\widetilde{T}_{i}(E_i)=v(K_i')^{-1}F_i,\qquad \widetilde{T}_{i}(F_i)=v^{-1}E_iK_i^{-1},\\
			&\widetilde{T}_i(E_j)=E_j,\qquad \widetilde{T}_i(F_j)=F_j,\qquad \text{ if }c_{ij}=0,
			\\
			&\widetilde{T}_{i}(E_j)=\frac{v^{\frac{1}{2}}E_iE_j-v^{-\frac{1}{2}}E_jE_i}{v-v^{-1}},\qquad  \widetilde{T}_{i}(F_j)=\frac{v^{\frac{1}{2}}F_iF_j-v^{-\frac{1}{2}}F_jF_i}{v-v^{-1}}, \text{ if } c_{ij}=-1. %\qquad  j\neq i.
		\end{align*}
		Moreover, there exists a group homomorphism $\Br(W)\rightarrow \Aut(\tU)$, $t_i\mapsto \widetilde{T}_i$ for $i\in\I$.    
	\end{proposition}
	
	Hence, we can define
	\begin{align}\widetilde{T}_w 
		:= \widetilde{T}_{i_1}\cdots
		\widetilde{T}_{i_r} \in \Aut(\tU),
	\end{align}
	where $w = s_{i_1}
	\cdots s_{i_r}$
	is any reduced expression of $w \in W$.
	
	\begin{lemma}
		\label{lem:QGbraid-bar}
		The braid group actions $\widetilde{T}_{i}$ commute with the bar-involution, i.e., $\ov{\widetilde{T}_i(u)}=\widetilde{T}_i(\ov{u})$ for any $u\in\tU$.
	\end{lemma}
	
	\begin{proof}
		It is enough to check $\ov{\widetilde{T}_i(u)}=\widetilde{T}_i(\ov{u})$ for $u=E_i,F_i,K_i,K_i'$, which is obvious by Proposition \ref{prop:braid1}. 
	\end{proof}

	\begin{proposition}[\text{\cite[Proposition 1.10]{Lus90a}}]
		\label{prop:PBW-QG}
		Let $w_0$ be the longest element of $W$ and fix a reduced expression ${w_0}=s_{i_1} \cdots s_{i_l}$ for the longest element $\omega_0\in W$. Set
		\begin{align}
			\label{eq:F-root}
			F_{\beta_k}=\widetilde{T}_{i_1}^{-1}\widetilde{T}_{i_2}^{-1}\cdots \widetilde{T}_{i_{k-1}}^{-1}(F_{i_k}),\qquad 
			E_{\beta_k}=\widetilde{T}_{i_1}^{-1}\widetilde{T}_{i_2}^{-1}\cdots \widetilde{T}_{i_{k-1}}^{-1}(E_{i_k}),
		\end{align}
		for $1\leq k\leq l,a\in \N$. Then the monomials 
		$$F^{\ba}=F_{\beta_1}^{a_1} F_{\beta_2}^{a_2}\cdots F_{\beta_l}^{a_l}, \qquad \ba=(a_1,\ldots, a_l)\in \N^l$$ 
		form a $\Q(v^{1/2})$-basis for $\tU^-$, and  the monomials $$E^{\ba}=E_{\beta_1}^{a_1} E_{\beta_2}^{a_2}\cdots E_{\beta_l}^{a_l}, \qquad \ba=(a_1,\ldots, a_l)\in \N^l$$ 
		form a $\Q(v^{1/2})$-basis for $\tU^+$.
		Moreover, the monomials $$F^{\ba}E^{\bc}K_\mu K'_{\nu},\qquad \ba,\bc\in\N^l,\mu,\nu\in\Z^\I$$
		form a $\Q(v^{1/2})$-basis for $\tU$. Similar result holds for $\hU$. 
	\end{proposition}

	\subsection{The $\imath$quantum groups}
	
	For a  Cartan matrix $C=(c_{ij})$, let $\Aut(C)$ be the group of all permutations $\btau$ of the set $\I$ such that $c_{ij}=c_{\btau i,\btau j}$. An element $\btau\in\Aut(C)$ is called an \emph{involution} if $\btau^2=\Id$.

	In this paper, we always assume that $\btau\in\Aut(C)$ is an involution such that $c_{i,\btau i}=0$ for all $i\neq \btau i$. We denote by $\bs_{i}$ the following element of order 2 in the Weyl group $W$, i.e.,
	\begin{align}
		\label{def:simple reflection}
		\bs_i= \left\{
		\begin{array}{ll}
			s_{i}, & \text{ if } \btau i=i;
			\\
			s_is_{\btau i}, & \text{ if } \btau i\neq i.
		\end{array}
		\right.
	\end{align}
	It is well known (cf., e.g., \cite{KP11}) that the {\rm restricted Weyl group} associated to the quasi-split symmetric pair $(\fg, \fg^\theta)$ can be identified with the following subgroup $W_\btau$ of $W$:
	\begin{align}
		\label{eq:Wtau}
		W_{\btau} =\{w\in W\mid \btau w =w \btau\}
	\end{align}
	where $\btau$ is regarded as an automorphism of the root lattice $\Z^\I$. Moreover, the restricted Weyl group $W_{\btau}$ can be identified with a Weyl group with $\bs_i$ ($i\in \I_\btau$) as its simple reflections. %More precisely, associated to the Dynkin $\imath$quivers $(Q,\btau)$, we have
	%\begin{align*}
	%W_{\btau} =
	%\begin{cases}
	%W, & \text{ if } \btau =\Id,
	%\\
	% W(B_{r+1}), & \text{ if } \btau \neq \Id, \text{ for }Q \text{ of type } A_{2r+1} \text{ or }D_{2r},
	%\\
	%W(F_4), & \text{ if } \btau \neq \Id, \text{ for } Q \text{ of type } E_6.
	%\end{cases}
	%\end{align*} 

	For a  Cartan matrix $C=(c_{ij})$, 
	let $\btau$ be an involution in $\Aut(C)$. We define the universal $\imath$quantum groups ${\hU}^\imath:=\hU'_v(\fk)$ (resp. $\tUi:=\tU'_v(\fk)$) to be the $\Q(v^{1/2})$-subalgebra of $\hU$ (resp. $\tU$) generated by
	\begin{equation}
		\label{eq:Bi}
		B_i= F_i +  E_{\btau i} \tK_i',
		\qquad \tk_i = \tK_i \tK_{\btau i}', \quad \forall i \in \I,
	\end{equation}
	(with $\tk_i$ invertible in $\tUi$).
	Let $\hU^{\imath 0}$ be the $\Q(v^{1/2})$-subalgebra of $\hUi$ generated by $\tk_i$, for $i\in \I$. Similarly, let $\tU^{\imath 0}$ be the $\Q(v^{1/2})$-subalgebra of $\tUi$ generated by $\tk_i^{\pm1}$, for $i\in \I$. 
	%By \cite[Lemma 6.1]{LW19}, the elements $\tk_i$ (for $i= \btau i$) and $\tk_i \tk_{\btau i}$  (for $i\neq \btau i$) are central in $\hUi$ (and also $\tUi$). 

	A presentation for $\tUi$ can be found in \cite[Proposition~ 6.4]{LW19}. The following is a modified version based on the quantum groups given in \S\ref{subsec:QG} (this kind of presentation also holds for $\hUi$ by omitting that $\tk_i$ ($i\in\I$) are invertible).
	\begin{lemma}
		\label{prop:Serre}
		%Let $(Q, \btau)$ be a Dynkin $\imath$quiver. 
		The $\Q(v)$-algebra $\tUi$ has a presentation with generators $B_i, \tk_i$ $(i\in \I)$, where $\tk_i$ are invertible, subject to the relations \eqref{relation1}--\eqref{relation2}: for $\ell \in \I$, and $i\neq j \in \I$,
		\begin{align}
			\tk_i \tk_\ell =\tk_\ell \tk_i,
			\quad
			\tk_\ell B_i & = v^{c_{\btau \ell,i} -c_{\ell i}} B_i \tk_\ell,
			%K_\mu B_i-q_i^{-\langle \mu,\alpha_i\rangle} B_iK_\mu & =0,
			\label{relation1}
			\\
			B_iB_{j}-B_jB_i &=0, \quad \text{ if } c_{ij} =0 \text{ and }\btau i\neq j,
			\label{relationBB}
			\\
			\sum_{s=0}^{1-c_{ij}} (-1)^s \qbinom{1-c_{ij}}{s} B_i^{s}B_jB_i^{1-c_{ij}-s} &=0, \quad \text{ if } j \neq \btau i\neq i,
			\\
			B_{\btau i}B_i -B_i B_{\btau i}& =  (v^{-1}-v) (\tk_i -\tk_{\btau i}),
			%\sum_{s=0}^{1-c_{i,\DTr i}} (-1)^s B_i^{(s)}B_{\DTr i}B_i^{(1-c_{i,\DTr i}-s)}& =\frac{1}{q_i-q_i^{-1}}
			%  \\
			% \cdot \left(q_i^{c_{i,\DTr i}} (q_i^{-2};q_i^{-2})_{-c_{i,\DTr i}} \vs_{\DTr i}B_i^{(-c_{i,\DTr i})} \tK_i \tK_{\DTr i}^{-1} \right.  & \left. -(q_i^{2};q_i^{2})_{-c_{i,\tau i}}\vs_{i}B_i^{(-c_{i,\tau i})} \tK_{\tau i} \tK_i^{-1} \right),
			\quad \text{ if } \btau i \neq i,
			\label{relation5}
			\\
			B_i^2B_{j} - [2] B_iB_{j}B_i +B_{j}B_i^2 &= -v(v-v^{-1})^2 \tk_i B_{j},
			%\sum_{s=0}^{1-c_{ij}} (-1)^s  B_{i,\overline{c_{ij}}+\overline{p_i}}^{(s)}B_j B_{i,\overline{p}_i}^{(1-c_{ij}-s)} &=0,\quad
			%
			\quad \text{ if }  c_{ij} = -1 \text{ and }\btau i=i.
			\label{relation2}
		\end{align}
	\end{lemma}
	
	%Let $\bvs=(\vs_i)\in  (\Q(v^{1/2})^\times)^{\I}$ be such that $\vs_i=\vs_{\btau i}$ for each $i\in \I$ which satisfies $c_{i, \btau i}=0$. 	
	%Let $\Ui:=\Ui_{\bvs}$ be the $\Q(v)$-subalgebra of $\bU$ generated by
	%\[
	%B_i= F_i+\vs_i E_{\btau i}K_i^{-1},
	%\quad
	%k_j= K_jK_{\btau j}^{-1},
	%\qquad  \forall i \in \I, j \in \ci.
	%\]
	It is known \cite{Let99, Ko14,LW19} that the algebra $\widetilde{\bU}^\imath$ (resp. $\hUi$) is a right coideal subalgebra of $\widetilde{\bU}$ (resp. $\hU$); we call $(\widetilde{\bU}, \widetilde{\bU}^\imath)$ and $(\hU,\hUi)$ quantum symmetric pairs.
	%\end{proposition}
	
	We shall refer to $\hUi$ and $\tUi$ %and $\Ui$ 
	as the universal {\em (quasi-split) $\imath${}quantum groups} (cf. the $\imath$quantum groups defined in  \cite{Let99,Ko14}); they are called {\em split} if $\btau =\Id$.

	\begin{example}
		\label{ex:QGvsiQG}
		Let us explain quantum groups are $\imath$quantum group of diagonal type. 
		Consider the $\Q(v)$-subalgebra $\tUUi$ of $\tUU$
		generated by
		\[
		\ck_i:=\tK_{i} \tK_{i^{\diamond}}', \quad
		\ck_i':=\tK_{i^{\diamond}} \tK_{i}',  \quad
		\cb_{i}:= F_{i}+ E_{i^{\diamond}} \tK_{i}', \quad
		\cb_{i^{\diamond}}:=F_{i^{\diamond}}+ E_{i} \tK_{i^{\diamond}}',
		\qquad \forall i\in \I.
		\]
		Here we drop the tensor product notation and use instead $i^\diamond$ to index the generators of the second copy of $\tU$ in $\tUU$. There exists a $\Q(v)$-algebra isomorphism $\widetilde{\phi}: \tU \rightarrow \tUUi$ such that
		\[
		\widetilde{\phi}(E_i)= \cb_{i},\quad \widetilde{\phi}(F_i)= \cb_{i^{\diamond}}, \quad \widetilde{\phi}(\tK_i)= \ck_i', \quad \widetilde{\phi}(\tK_i')= \ck_i, \qquad \forall  i\in \I.
		\]
		% (consistent with the notation $Q^{\dbl}=Q\sqcup Q^\diamond$).
	\end{example}

	Choose one representative for each $\btau$-orbit on $\I$, and let
	\begin{align}\label{eq:ci}
		\ci = \{ \text{the chosen representatives of $\btau$-orbits in $\I$} \}.
	\end{align} 
	The braid group associated to the relative Weyl group for $(\fg, \fg^\theta)$ is of the form
	\begin{equation}
		\label{eq:braidCox}
		\brW =\langle \br_i \mid i\in \I_\btau \rangle
	\end{equation}
	where $\br_i$ satisfy the same braid relations as for $\bs_i$ in $W_{\varrho}$ (but no quadratic relations on $r_i$ are imposed). The following result was proved in \cite{LW21a} using $\imath$Hall algebra technique (cf. \cite{KP11, WZ23}), %and a purely algebraic proof was obtained subsequently in \cite{Z23},
	where for any $i\in\I$, we set
	\begin{align}
		\label{eq:bbKi}
		\K_i:=v\tk_i, \text{ if }\varrho i=i;
		\qquad
		\K_j:=\tk_j, \text{ otherwise.}
	\end{align}
	and $\K_\alpha:=\prod_{i\in\I}\K_i^{a_i}$ if $\alpha\in\Z^\I$. For elements $x,y\in\Ui$, we also define $[x,y]_v:=xy-vyx$.
	
	\begin{theorem} [\text{\cite[Theorem 6.8]{LW21a}}]
		\label{thm:Ti}
		\begin{enumerate}
			\item[(1)] For $i\in \I$ such that $\varrho i=i$, there exists an automorphism $\tTT_i$ of the $\Q(v)$-algebra $\tUi$ such that
			$\tTT_i(\K_\mu) =\K_{\bs_i\mu}$, and
			\[
			\tTT_i(B_j)= \begin{dcases}
				\K_i^{-1} B_i,  &\text{ if }j=i,\\
				B_j,  &\text{ if } c_{ij}=0, \\
				\frac{v^{\frac{1}{2}}B_iB_j-v^{-\frac{1}{2}}B_jB_i}{v-v^{-1}},  & \text{ if }c_{ij}=-1, %\\
				%{[}2]^{-1} \big(B_jB_i^{2} -v[2] B_i B_jB_i +v^2 B_i^{2} B_j \big) + B_j\K_i,  & \text{ if }c_{ij}=-2,
			\end{dcases}
			\]
			for $\mu\in \Z^\I$ and $j\in \I$.
			
			\item[(2)] For $i\in \I$ such that $c_{i,\varrho i}=0$, there exists an automorphism $\tTT_i$ of the $\Q(v)$-algebra $\tUi$ such that
			$\tTT_i(\K_\mu)=\K_{\bs_i\mu}$, and
			\[
			\tTT_i(B_j)= \begin{dcases}
				\frac{v^{\frac{1}{2}}B_iB_j-v^{-\frac{1}{2}}B_jB_i}{v-v^{-1}},  & \text{ if }c_{ij}=-1 \text{ and } c_{\varrho i,j}=0,
				\\
				\frac{v^{\frac{1}{2}}B_{\varrho i}B_j-v^{-\frac{1}{2}}B_jB_{\varrho i}}{v-v^{-1}},  & \text{ if } c_{ij}=0 \text{ and }c_{\varrho i,j}=-1 ,
				\\
				\frac{v^{-1}\big[[B_j,B_i]_v,B_{\varrho i}\big]_v}{(v-v^{-1})^2}+B_j\K_i  & \text{ if } c_{ij}=-1 \text{ and }c_{\varrho i,j}=-1 ,
				\\
				v\K_{i}^{-1} B_{\btau i},  & \text{ if }j=i,
				\\
				v\K_{\btau i}^{-1} B_i,  &\text{ if }j=\varrho i,
				\\
				B_j, & \text{ otherwise;}
			\end{dcases}
			\]
			for $\mu\in \Z^\I$ and $j\in \I$. 
		\end{enumerate}
		Moreover, there exists a homomorphism $\brW \rightarrow \Aut( \tUi)$, $\br_i\mapsto \tTT_i$, for all $i\in \I_\btau$.
	\end{theorem}
	
	Similar to Lemma \ref{QG bar-involution def} and Lemma \ref{lem:QGbraid-bar}, we have the following two lemmas.
	
	\begin{lemma}
		\label{iQG bar-involution def}
		There exists an anti-involution $u\mapsto \ov{u}$ on $\hUi$ (also $\tUi$) given by $\ov{v^{1/2}}=v^{-1/2}$, $\ov{B_i}=B_i$, $\ov{\K_i}=\K_i$, for $i\in\I$. In particular, $\ov{\tk_i}=\tk_i$ if $\varrho i\neq i$; $\ov{\tk_i}=v^2\tk_i$ if $\varrho i=i$.
		%  \begin{align}
		% \ov{\tk_i}=\begin{cases}
		%   \tk_i & \text{ if $\varrho i\neq i$},\\
		%   v^2\tk_i&\text{ if $\varrho i=i$}.
		%  \end{cases}
		%\end{align}
	\end{lemma}

	\begin{lemma}
		The braid group actions $\widetilde{\TT}_{i}$ commute with the bar-involution, i.e., $\ov{\widetilde{\TT}_i(u)}=\widetilde{\TT}_i(\ov{u})$ for any $u\in\tUi$.
	\end{lemma}
	
	Now we review the PBW basis of $\tUi$ (similar for $\bU^\imath$), cf. \cite[\S 8]{LW21a}. Let $\omega_0$ be the longest element of the Weyl group $W$ (note that $\omega_0\in W_\btau$). Let $\omega_0=r_{i_1}r_{i_2}\cdots r_{i_\ell}$ be a reduced expression and write
	\begin{align*}
		\beta_k=r_{i_1}r_{i_2}\cdots r_{i_{k-1}}(\alpha_{i_k}),\qquad \forall 1\leq k\leq \ell.
	\end{align*}
	Then we have
	$\{\beta_1,\btau(\beta_1), \beta_2,\btau(\beta_2),\cdots,\beta_\ell,\btau(\beta_\ell)\}=\Phi^+$, which also gives a total order of $\Phi^+$. Here $\Phi^+$ is the set of positive roots of $\fg$. For any $\gamma\in\Phi^+$, we define
	\begin{align*}
		B_\gamma=\begin{cases}
			\tTT_{i_1}^{-1}\tTT_{i_2}^{-1}\cdots \tTT_{i_{k-1}}^{-1}(B_{i_k}) & \text{ if }\gamma=\beta_k,
			\\
			\tTT_{i_1}^{-1}\tTT_{i_2}^{-1}\cdots \tTT_{i_{k-1}}^{-1}(B_{\varrho(i_k)}) & \text{ if }\gamma=\varrho(\beta_k).
		\end{cases}
	\end{align*}
	and for any $\ba\in\N^{|\Phi^+|}$, set
	\begin{equation}\label{eq:iQG PBW basis def}
		B^{\ba}:=\prod_{\gamma\in\Phi^+} B_{\gamma}^{a_\gamma}.
	\end{equation}
	
	\begin{proposition}[\text{\cite[Theorem 8.14]{LW21a}}]
		\label{prop:PBW-Ui}
		Given  a reduced expression $\omega_0=r_{i_1}r_{i_2}\cdots r_{i_\ell}$ of the longest element $\omega_0\in W_\btau$, the monomials $\{B^{\ba}\tk_\mu\mid \ba\in\N^{|\Phi^+|},\mu\in\Z^\I\}$
		form a $\Q(v^{1/2})$-basis of $\tUi$.   
	\end{proposition}

	%%%%%%%%%%%%%%%
	\section{$\imath$Quiver algebras and $\imath$Hall algebras}\label{sec:iHA}
	
	Throughout this paper, we denote $\bfk=\F_q$ the finite field of $q$ elements.

	\subsection{Hall algebras}
	%%%%%%%%%%%%
	
	Let $\ca$ be an essentially small exact category, linear over the finite field $\bfk=\F_q$.
	Assume that $\ca$ has finite morphism and extension spaces:
	$$|\Hom_\ca(A,B)|<\infty,\quad |\Ext^1_\ca(A,B)|<\infty,\,\,\forall A,B\in\ca.$$
	
	% Let $K_0(\ca)$ the Grothendieck group of $\mod(\ca)$. We denote by $\widehat{M}$ the class in $K_0(\ca)$ for $M\in\ca$.  
	%It is well-known $K_0(\bfk Q)$ is a free abelian group generated by $\widehat{S_i}$ ($\i\in\I$), where $S_i$ is the simple module corresponding to the vertex $i\in\I$. We identify $K_0(\bfk Q)$ with the root lattice $\Z^\I$; see \S\ref{subsec:QG}. 
	%Let $K^+_0(\ca)$ denote the subset of $K_0(\ca)$ corresponding to classes of objects in $\ca$.

	We denote by $\Iso(\ca)$ the set of isoclasses of objects of $\ca$.	Given objects $X,Y,Z\in\ca$, define $\Ext^1_\ca(X,Z)_Y\subseteq \Ext^1_\ca(X,Z)$ to be the subset parameterising extensions with the middle term isomorphic to $Y$. We define the Ringel-Hall algebra (also called Hall algebra) $\ch(\ca)$ to be the $\Q$-vector space whose basis is formed by the isomorphism classes $[X]$ of objects $X$ of $\ca$, with the multiplication
	defined by
	\begin{align}
		\label{eq:mult}
		[X]\cdot [Z]=\sum_{[Y]\in \Iso(\ca)}G_{XZ}^Y[Y],\text{ where }G_{XZ}^Y:=\frac{|\Ext_\ca^1(X,Z)_Y|}{|\Hom_\ca(X,Z)|}.
	\end{align}
	It is well known that
	the algebra $\ch(\ca)$ is associative and unital. The unit is given by $[0]$, where $0$ is the zero object of $\ca$; see \cite{Rin90,Br13}. % and also \cite{Rie,P,Hub,Br}.

	For any three objects $X,Y,Z$, let
	\begin{align}
		\label{eq:Fxyz}
		F_{XZ}^Y&:= \big |\{L\subseteq Y \mid L \cong Z,  Y/L\cong X\} \big |.\\
		\label{eq:P_xyz}
		P_{XZ}^Y&:=|\{(f,g)\mid 0\rightarrow Z\stackrel{f}{\rightarrow} Y\stackrel{g}{\rightarrow} X\rightarrow 0\text{ is short exact}\}|.
	\end{align}
	Let $\Aut(X)$ be the automorphism group of $X\in\ca$, and denote by $a_X=|\Aut(X)|$.  The Riedtmann-Peng formula states that
	\[
	F_{XZ}^Y= G_{XZ}^Y \cdot \frac{a_Y}{a_X a_Z}.
	\]
	It is clear that $P_{XZ}^Y=F_{XZ}^Ya_Xa_Z$.

	%%%%%
	
	\subsection{Notations for quivers and representations}
	\label{subsec:Quiver-represent}
	
	Let $Q=(Q_0,Q_1,\ts,\tt)$ be a quiver, where, $Q_0$ is the set of vertices (we sometimes write $\I =Q_0$), $Q_1$ is the set of arrows, and $\ts,\tt:Q_1\rightarrow Q_0$ are two maps such that for the arrow $\alpha$ starts at $\ts(\alpha)$, and ends at $\tt(\alpha)$. Let $I$ be an admissible ideal of the path algebra $\bfk Q$, and $A=\bfk Q/I$ the (finite-dimensional) quotient algebra. 
	
	Let $\rep_\bfk(Q)$ be the category of finite dimensional representations of $Q$ over $\bfk$, and we identify $\rep_\bfk(Q)$ with the category $\mod(\bfk Q)$ of finite dimensional left modules over the path algebra $\bfk Q$. Similarly, let $\rep(A)$ be the category of finite dimensional representations of $Q$ satisfying the relations in $I$, and identify $\rep(A)$ with $\mod(A)$ of finite dimensional left modules over $A$. Let $\proj(A)$ the subcategory of projective $A$-modules. Denote by ${\rm proj.dim}_AM$ (resp. ${\rm inj.dim}_AM$) the  projective (resp. injective) dimension of an $A$-module $M$.
	
	%	$\triangleright$ ${\rm inj.dim}_AM$ -- injective dimension of $M$

	For a representation $V=(V_i,V_h)_{i\in\I,h\in Q_1}$, the vector $\bd=\dimv V:=(\dim_\bfk V_i)_i\in\N^\I$ is called the dimension vector of $V$.  For a dimension vector $\bd$, let $\rep(\bd,A)$ be the variety of $A$-modules of dimension vector $\bd$. Let $G_\bd:=\prod_{i\in\I}\mathrm{GL}(d_i,\bfk)$. The group $G_\bd$ acts on $\rep(\bd,A)$ by conjugation: 
	\[
	(g_i)_{i\in\I}\cdot (V_i,V_h)=(V_i,g_{\tt(h)}V_\alpha g_{\ts(h)}^{-1}),\quad \forall (V_i,V_h)\in\rep(\bd,A). 
	\]
	Let $\mathfrak{O}_V$ be the $G_\bd$-orbit of $V$.  
	
	Let $S_i$ be the ($1$-dimensional) simple $A$-modules to each $i\in\I$. The Grothendieck group $K_0(\mod(A))$ of $\mod(A)$ is a free abelian group generated by $\widehat{S_i}$ $(i\in\I)$. Denote by $\widehat{M}$ the class of $M$ in $K_0(\mod(A))$. 
	We identify $K_0(\mod(A))$ (especially $K_0(\mod(\bfk Q))$) with the root lattice 
	$\Z^\I=\oplus_{i\in\I} \Z\alpha_i$ (see \S\ref{subsec:QG}) by identifying $\widehat{S_i}$ with the simple root $\alpha_i$, and $\widehat{M}$ with its dimension vector $\dimv M$ for any $M\in\mod(A)$. 
	Let $K^+_0(\mod(A))$ denote the subset of $K_0(\mod(A))$ corresponding to classes of objects in $\mod(A)$.
	In this way, we identify $K_0^+(\mod(A))$ with $\N^\I=\oplus_{i\in\I}\N\alpha_i$.

	Let $\langle-,- \rangle_Q$ be the Euler form of $Q$,  i.e.,
	\[
	\langle M,N\rangle_Q=\dim_\bfk\Hom_{\bfk Q}(M,N)-\dim_\bfk\Ext^1_{\bfk Q}(M,N)
	\]
	for $M,N\in\mod(\bfk Q)$. It is well-known that the Euler form descends to the Grothendieck group $K_0(\mod(\bfk Q))$. We denote by $(-,-)_Q:K_0(\mod(\bfk Q))\times K_0(\mod(\bfk Q))\rightarrow\Z$ the symmetrized Euler form defined by $(\alpha,\beta)_Q=\langle\alpha,\beta\rangle_Q+\langle \beta,\alpha\rangle_Q$.
	
	Let $\cd^b(A)$  be the bounded derived category with suspension functor $\Sigma$. The singularity category  $\cd_{sg}(A)$ of $\mod(A)$ is the Verdier quotient of $\cd^b(A)$ modulo the thick subcategory generated by modules with finite projective dimensions. %) be the singularity category of $\mod(A)$. 

	%%%
	\subsection{The $\imath$quivers and doubles}
	\label{subsec:i-quiver}

	An {\em involution} of $Q$ is defined to be an automorphism $\btau$ of the quiver $Q$ such that $\btau^2=\Id$ (we allow the {\em trivial} involution $\Id:Q\rightarrow Q$). An involution $\btau$ of $Q$ induces an involution of the path algebra $\bfk Q$, again denoted by $\btau$.
	A quiver together with an involution $\btau$, $(Q, \btau)$, will be called an {\em $\imath$quiver}. In this paper, we focus on Dynkin quivers.
	
	%Let $R_1$ denote the truncated polynomial algebra $\bfk[\varepsilon]/(\varepsilon^2)$.
	Let $R_2$ denote the radical square zero of the path algebra of $\xymatrix{1 \ar@<0.5ex>[r]^{\varepsilon} & 1' \ar@<0.5ex>[l]^{\varepsilon'}}$, i.e., $\varepsilon' \varepsilon =0 =\varepsilon\varepsilon '$. Define a $\bfk$-algebra
	\begin{equation}
		\label{eq:La}
		\Lambda=\bfk Q\otimes_\bfk R_2.
	\end{equation}
	
	Associated to the quiver $Q$, the {\em double framed quiver} $Q^\sharp$
	is the quiver such that
	\begin{itemize}
		\item the vertex set of $Q^{\sharp}$ consists of 2 copies of the vertex set $Q_0$, $\{i,i'\mid i\in Q_0\}$;
		\item the arrow set of $Q^{\sharp}$ is
		\[
		\{\alpha: i\rightarrow j,\alpha': i'\rightarrow j'\mid(\alpha:i\rightarrow j)\in Q_1\}\cup\{ \varepsilon_i: i\rightarrow i' ,\varepsilon'_i: i'\rightarrow i\mid i\in Q_0 \}.
		\]
	\end{itemize}
	Let $I^{\sharp}$ be the admissible ideal of $\bfk Q^{\sharp}$ generated by
	\begin{itemize}
		\item
		(Nilpotent relations) $\varepsilon_i \varepsilon_i'$, $\varepsilon_i'\varepsilon_i$ for any $i\in Q_0$;
		\item
		(Commutative relations) $\varepsilon_j' \alpha' -\alpha\varepsilon_i'$, $\varepsilon_j \alpha -\alpha'\varepsilon_i$ for any $(\alpha:i\rightarrow j)\in Q_1$.
	\end{itemize}
	Then the algebra $\La$ can be realized as $\Lambda\cong \bfk Q^{\sharp} \big/ I^{\sharp}$.
	
	Note $Q^\sharp$ admits a natural involution, $\swa$.
	The involution $\btau$ of a quiver $Q$ induces an involution ${\btau}^{\sharp}$ of $Q^{\sharp}$ which is the composition of $\swa$ and $\btau$ (on the two copies of subquivers $Q$ and $Q'$ of $Q^\sharp$).
	%defined by
	%\begin{itemize}
	%\item ${\btau}^{\sharp}(i)=(\btau i)'$, ${\btau}^{\sharp}(i') =\btau i$ for any $i\in Q_0$;
	%\item ${\btau}^{\sharp}(\varepsilon_i)= \varepsilon_{\btau i}'$, ${\btau}^{\sharp}(\varepsilon_i')= \varepsilon_{\btau i}$ for any $i\in Q_0$;
	%\item ${\btau}^{\sharp}(\alpha)= (\btau\alpha)'$, ${\btau}^{\sharp}(\alpha')=\btau\alpha$ for any $\alpha\in Q_1$.
	%\end{itemize}
	%The algebra $\La$ can be realized in terms of the quiver $Q^{\sharp}$ and a certain admissible ideal $I^{\sharp}$
	%so that $\Lambda\cong \bfk Q^{\sharp} \big/ I^{\sharp}$; see \cite[\S2.2]{LW19}.
	
	%The algebra $\Lambda$ can be described in terms of a quiver with relations. Let $I^{\sharp}$ be the admissible ideal of $\bfk Q^{\sharp}$ generated by
	%\begin{itemize}
	%\item
	%(Nilpotent relations) $\varepsilon_i \varepsilon_i'$, $\varepsilon_i'\varepsilon_i$ for any $i\in Q_0$;
	%\item
	%(Commutative relations) $\varepsilon_j' \alpha' -\alpha\varepsilon_i'$, $\varepsilon_j \alpha -\alpha'\varepsilon_i$ for any $(\alpha:i\rightarrow j)\in Q_1$.
	%\end{itemize}
	%
	%Let $Q$ (respectively, $Q'$) be the full subquiver of $Q^{\sharp}$ formed by all vertices $i$ (respectively, $i'$) for $i\in Q_0$. Then $Q\sqcup Q'$ is a subquiver of $Q^{\sharp}$.
	
	By \cite[Lemma~2.4]{LW19}, ${\btau}^{\sharp}$ on $Q^\sharp$ preserves $I^\sharp$ and hence induces an involution ${\btau}^{\sharp}$ on the algebra $\Lambda$. By \cite[Definition 2.5]{LW19}, the {\rm $\imath$quiver algebra} of $(Q, \btau)$ is the fixed point subalgebra of $\Lambda$ under ${\btau}^{\sharp}$,
	\begin{equation}
		\label{eq:iLa}
		\iLa
		= \{x\in \Lambda\mid {\btau}^{\sharp}(x) =x\}.
	\end{equation}
	The algebra $\iLa$ can be described in terms of a certain quiver $\ov{Q}$ and its admissible ideal $\ov{I}$ so that $\iLa \cong \bfk \ov{Q} / \ov{I}$; see \cite[Proposition 2.6]{LW19}.
	We recall $\ov{Q}$ and $\ov{I}$ as follows:
	\begin{itemize}
		\item[(i)] $\ov{Q}$ is constructed from $Q$ by adding a loop $\varepsilon_i$ at the vertex $i\in Q_0$ if $\btau i=i$, and adding an arrow $\varepsilon_i: i\rightarrow \btau i$ for each $i\in Q_0$ if $\btau i\neq i$;
		\item[(ii)] $\ov{I}$ is generated by
		\begin{itemize}
			\item[(1)] (Nilpotent relations) $\varepsilon_{i}\varepsilon_{\btau i}$ for any $i\in\I$;
			\item[(2)] (Commutative relations) $\varepsilon_i\alpha-\btau(\alpha)\varepsilon_j$ for any arrow $\alpha:j\rightarrow i$ in $Q_1$.
		\end{itemize}
	\end{itemize}
	Moreover, it follows by \cite[Proposition 3.5]{LW19} that $\Lambda^{\imath}$ is a $1$-Gorenstein algebra.
	
	The following quivers are examples of the quivers $\ov{Q}$ used to describe the $\imath$quiver algebras $\Lambda^\imath$
	associated to non-split $\imath$quivers of type ADE; cf. \cite{LW19}.

	%%%A
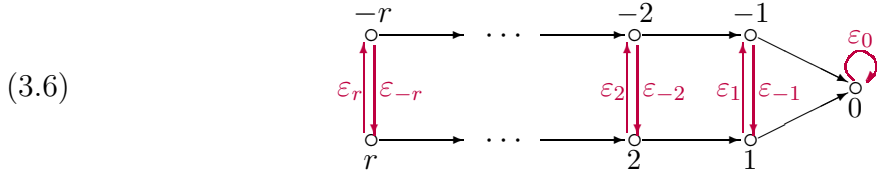
	\begin{center}\setlength{\unitlength}{0.7mm}
		\vspace{-1.5cm}
		\begin{equation}
			\label{diag: A}
			\begin{picture}(100,40)(0,20)
				\put(0,10){$\circ$}
				\put(0,30){$\circ$}
				\put(50,10){$\circ$}
				\put(50,30){$\circ$}
				\put(72,10){$\circ$}
				\put(72,30){$\circ$}
				\put(92,20){$\circ$}
				\put(0,6){$r$}
				\put(-2,34){${-r}$}
				\put(50,6){\small $2$}
				\put(48,34){\small ${-2}$}
				\put(72,6){\small $1$}
				\put(70,34){\small ${-1}$}
				\put(92,16){\small $0$}
				
				\put(3,11.5){\vector(1,0){16}}
				\put(3,31.5){\vector(1,0){16}}
				\put(23,10){$\cdots$}
				\put(23,30){$\cdots$}
				\put(33.5,11.5){\vector(1,0){16}}
				\put(33.5,31.5){\vector(1,0){16}}
				\put(53,11.5){\vector(1,0){18.5}}
				\put(53,31.5){\vector(1,0){18.5}}
				\put(75,12){\vector(2,1){17}}
				\put(75,31){\vector(2,-1){17}}
				\color{purple}
				\put(0,13){\vector(0,1){17}}
				\put(2,29.5){\vector(0,-1){17}}
				\put(50,13){\vector(0,1){17}}
				\put(52,29.5){\vector(0,-1){17}}
				\put(72,13){\vector(0,1){17}}
				\put(74,29.5){\vector(0,-1){17}}
				
				\put(-5,20){$\varepsilon_r$}
				\put(3,20){$\varepsilon_{-r}$}
				\put(45,20){\small $\varepsilon_2$}
				\put(53,20){\small $\varepsilon_{-2}$}
				\put(67,20){\small $\varepsilon_1$}
				\put(75,20){\small $\varepsilon_{-1}$}
				\put(92,30){\small $\varepsilon_0$}
				
				\qbezier(93,23)(90.5,25)(92,27)
				\qbezier(92,27)(94,30)(97,27)
				\qbezier(97,27)(98,24)(95.5,22.6)
				\put(95.6,23){\vector(-1,-1){0.3}}
			\end{picture}
		\end{equation}
		\vspace{-0.6cm}
	\end{center}

	%%%D
	\begin{center}\setlength{\unitlength}{0.8mm}
		
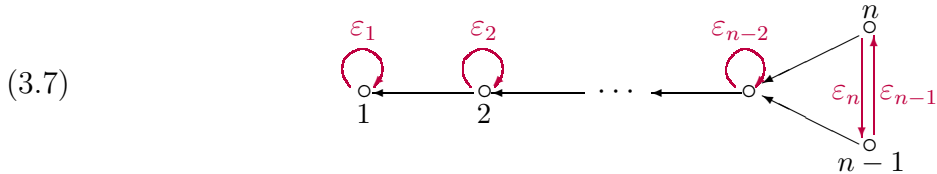
\begin{equation}
			\label{diag: D}
			\begin{picture}(100,25)(-5,0)
				\put(0,-1){$\circ$}
				\put(0,-5){\small$1$}
				\put(20,-1){$\circ$}
				\put(20,-5){\small$2$}
				\put(64,-1){$\circ$}
				\put(84,-10){$\circ$}
				\put(80,-13){\small${n-1}$}
				\put(84,9.5){$\circ$}
				\put(84,12.5){\small${n}$}

				\put(19.5,0){\vector(-1,0){16.8}}
				\put(38,0){\vector(-1,0){15.5}}
				\put(64,0){\vector(-1,0){15}}
				
				\put(40,-1){$\cdots$}
				\put(83.5,9.5){\vector(-2,-1){16}}
				\put(83.5,-8.5){\vector(-2,1){16}}
				\color{purple}
				\put(86,-7){\vector(0,1){16.5}}
				\put(84,9){\vector(0,-1){16.5}}
				
				\qbezier(63,1)(60.5,3)(62,5.5)
				\qbezier(62,5.5)(64.5,9)(67.5,5.5)
				\qbezier(67.5,5.5)(68.5,3)(66.4,1)
				\put(66.5,1.4){\vector(-1,-1){0.3}}
				\qbezier(-1,1)(-3,3)(-2,5.5)
				\qbezier(-2,5.5)(1,9)(4,5.5)
				\qbezier(4,5.5)(5,3)(3,1)
				\put(3.1,1.4){\vector(-1,-1){0.3}}
				\qbezier(19,1)(17,3)(18,5.5)
				\qbezier(18,5.5)(21,9)(24,5.5)
				\qbezier(24,5.5)(25,3)(23,1)
				\put(23.1,1.4){\vector(-1,-1){0.3}}
				
				\put(-1,9.5){$\varepsilon_1$}
				\put(19,9.5){$\varepsilon_2$}
				\put(59,9.5){$\varepsilon_{n-2}$}
				\put(79,-1){$\varepsilon_{n}$}
				\put(87,-1){$\varepsilon_{n-1}$}
			\end{picture}
		\end{equation}
		\vspace{.8cm}
	\end{center}

	%%%%E
	\begin{center}\setlength{\unitlength}{0.8mm}
		
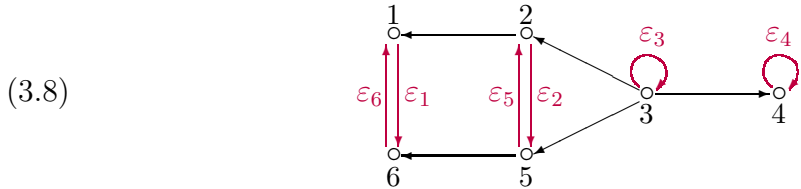
\vspace{-3cm}
		\begin{equation}
			\label{diag: E}
			\begin{picture}(100,40)(0,20)
				\put(10,6){\small${6}$}
				\put(10,10){$\circ$}
				\put(32,6){\small${5}$}
				\put(32,10){$\circ$}
				
				\put(10,30){$\circ$}
				\put(10,33){\small${1}$}
				\put(32,30){$\circ$}
				\put(32,33){\small${2}$}
				
				\put(31.5,11){\vector(-1,0){19}}
				\put(31.5,31){\vector(-1,0){19}}
				
				\put(52,22){\vector(-2,1){17.5}}
				\put(52,20){\vector(-2,-1){17.5}}
				
				\put(54.7,21.2){\vector(1,0){19}}
				
				\put(52,20){$\circ$}
				\put(52,16.5){\small$3$}
				\put(74,20){$\circ$}
				\put(74,16.5){\small$4$}
				\color{purple}
				\put(10,12.5){\vector(0,1){17}}
				\put(12,29.5){\vector(0,-1){17}}
				\put(32,12.5){\vector(0,1){17}}
				\put(34,29.5){\vector(0,-1){17}}
				
				\qbezier(52,22.5)(50,24)(51,26.5)
				\qbezier(51,26.5)(53,29)(56,26.5)
				\qbezier(56,26.5)(57.5,24)(55,22)
				\put(55.1,22.4){\vector(-1,-1){0.3}}
				\qbezier(74,22.5)(72,24)(73,26.5)
				\qbezier(73,26.5)(75,29)(78,26.5)
				\qbezier(78,26.5)(79,24)(77,22)
				\put(77.1,22.4){\vector(-1,-1){0.3}}
				
				\put(35,20){$\varepsilon_2$}
				\put(27,20){$\varepsilon_5$}
				\put(13,20){$\varepsilon_1$}
				\put(5,20){$\varepsilon_6$}
				\put(52,30){$\varepsilon_3$}
				\put(73,30){$\varepsilon_4$}
			\end{picture}
		\end{equation}
		\vspace{1cm}
	\end{center}

	\begin{example}[$\imath$quivers of diagonal type]
		\label{ex:diagquiver}
		Let $Q$ be an arbitrary quiver, and $Q^{\dbl} =Q\sqcup  Q^{\diamond}$,  where $Q^{\diamond}$ is an identical copy of $Q$  with a vertex set $\{i^{\diamond} \mid i\in Q_0\}$ and an arrow set $\{ \alpha^{\diamond} \mid \alpha \in Q_1\}$. We let $\rm{swap}$ be the involution of $Q^{\rm dbl}$ uniquely determined by $\swa(i)=i^\diamond$ for any $i\in Q_0$. 
		Then $(Q^{\rm dbl},\mathrm{swap})$ is an $\imath$quiver with $\Lambda$  as its $\imath$quiver algebra; see \cite[Example 2.10]{LW19}. 
	\end{example}

	%\begin{proposition}   \cite[Proposition 2.6]{LW19}
	%  \label{prop:invariant subalgebra}
	%We have $\iLa \cong \bfk \ov{Q} / \ov{I}$, where
	%\begin{itemize}
	%\item[(i)] $\ov{Q}$ is constructed from $Q$ by adding a loop $\varepsilon_i$ at the vertex $i\in Q_0$ if $\btau i=i$, and adding an arrow $\varepsilon_i: i\rightarrow \btau i$ for each $i\in Q_0$ if $\btau i\neq i$;
	%\item[(ii)] $\ov{I}$ is generated by
	%\begin{itemize}
	%\item[(1)] (Nilpotent relations) $\varepsilon_{i}\varepsilon_{\btau i}$ for any $i\in\I$;
	%\item[(2)] (Commutative relations) $\varepsilon_i\alpha-\btau(\alpha)\varepsilon_j$ for any arrow $\alpha:j\rightarrow i$ in $Q_1$.
	%\end{itemize}
	%\end{itemize}
	%\end{proposition}
	
	By \cite[Corollary 2.12]{LW19}, $\bfk Q$ is naturally a subalgebra and also a quotient algebra of $\Lambda^\imath$.
	Viewing $\bfk Q$ as a subalgebra of $\Lambda^{\imath}$, we have a restriction functor
	\[
	\res: \mod (\Lambda^{\imath})\longrightarrow \mod (\bfk Q).
	\]
	Viewing $\bfk Q$ as a quotient algebra of $\Lambda^{\imath}$, we obtain a pullback functor
	\begin{equation}\label{eqn:rigt adjoint}
		\iota:\mod(\bfk Q)\longrightarrow\mod(\Lambda^{\imath}).
	\end{equation}
	Hence a simple module $S_i (i\in Q_0)$ of $\bfk Q$ is naturally a simple $\iLa$-module.
	
	For each $i\in Q_0$, define a $\bfk$-algebra (which can be viewed as a subalgebra of $\iLa$)
	\begin{align}\label{dfn:Hi}
		\BH _i:=\left\{ \begin{array}{cc}  \bfk[\varepsilon_i]/(\varepsilon_i^2) & \text{ if }\btau i=i,
			\\
			\bfk(\xymatrix{i \ar@<0.5ex>[r]^{\varepsilon_i} & \btau i \ar@<0.5ex>[l]^{\varepsilon_{\btau i}}})/( \varepsilon_i\varepsilon_{\btau i},\varepsilon_{\btau i}\varepsilon_i)  &\text{ if } \btau i \neq i .\end{array}\right.
	\end{align}
	Note that $\BH _i=\BH _{\btau i}$ for any $i\in Q_0$. %Recall $\ci$ from \eqref{eq:ci}.
	%Choose one representative for each $\btau$-orbit on $\I=Q_0$, and let
	%\begin{align}   \label{eq:ci}
	%\ci = \{ \text{the chosen representatives of $\btau$-orbits in $\I$} \}.
	%\end{align}
	
	Define the following subalgebra of $\Lambda^{\imath}$:
	\begin{equation}  \label{eq:H}
		\BH =\bigoplus_{i\in \ci }\BH _i.
	\end{equation}
	Note that $\BH $ is a radical square zero selfinjective algebra. Denote by
	\begin{align}
		\res_\BH :\mod(\iLa)\longrightarrow \mod(\BH )
	\end{align}
	the natural restriction functor.
	On the other hand, as $\BH $ is a quotient algebra of $\iLa$ (cf. \cite[proof of Proposition~ 2.15]{LW19}), every $\BH $-module can be viewed as a $\iLa$-module.
	
	%Recall the algebra $\BH _i$ for $i \in \ci$ from \eqref{dfn:Hi}.
	For $i\in \I$, define the indecomposable module over $\BH _i$ (if $i\in \ci$) or over $\BH_{\btau i}$ (if $i\not \in \ci$)
	\begin{align}
		\label{eq:E}
		\E_i =\begin{cases}
			\bfk[\varepsilon_i]/(\varepsilon_i^2), & \text{ if }\btau i=i;
			\\
			\xymatrix{\bfk\ar@<0.5ex>[r]^1 & \bfk\ar@<0.5ex>[l]^0} \text{ on the quiver } \xymatrix{i\ar@<0.5ex>[r]^{\varepsilon_i} & \btau i\ar@<0.5ex>[l]^{\varepsilon_{\btau i}} }, & \text{ if } \btau i\neq i.
		\end{cases}
	\end{align}
	Then $\E_i$, for $i\in Q_0$, can be viewed as a $\iLa$-module and will be called a {\em generalized simple} $\iLa$-module.
	
	\begin{lemma}[\text{\cite[Proposition 3.9]{LW19}}]
		\label{cor: res proj}
		For any $M\in\mod(\Lambda^{\imath})$ the following are equivalent: (i) $\pd M<\infty$; 
		(ii) $\ind M<\infty$;
		(iii) $\pd M\leq1$;
		(iv) $\ind M\leq1$;
		(v) $\res_\BH (M)$ is projective as an $\BH $-module.
		%\end{itemize}
	\end{lemma}
	%As a corollary, we know $\Lambda^\imath$ is a $1$-Gorenstein algebra (see \cite[Proposition 3.6]{LW19}).  
	
	%Let $\langle\cdot,\cdot \rangle_Q$ be the Euler form of $Q$,  i.e.,
	%\[
	%\langle M,N\rangle_Q=\dim_\bfk\Hom_{\bfk Q}(M,N)-\dim_\bfk\Ext^1_{\bfk Q}(M,N)
	%\]
	%for $M,N\in\mod(\bfk Q)$. It is well-known that the Euler form descends to the Grothendieck group $K_0(\bfk Q)$. We denote by $(-,-)_Q:K_0(\bfk Q)\times K_0(\bfk Q)\rightarrow\Z$ the symmetrised Euler form defined by $(\alpha,\beta)_Q=\langle\alpha,\beta\rangle_Q+\langle \beta,\alpha\rangle_Q$.
	
	We denote by $\cp^{\leq 1}(\Lambda^\imath)$ the subcategory of $\Lambda^\imath$-modules with finite projective dimension (equivalently, with projective diemnsion at most $1$). 
	For $K,M\in \mod(\iLa)$, if $K\in\cp^{\leq 1}(\Lambda^\imath)$, we define the Euler forms
	\begin{align}\label{left Euler form}
		\langle K,M\rangle=\sum_{i=0}^{+\infty}(-1)^i \dim_\bfk\Ext^i(K,M)=\dim_\bfk\Hom(K,M)-\dim_\bfk\Ext^1(K,M),\\
		\label{right Euler form}
		\langle M,K\rangle=\sum_{i=0}^{+\infty}(-1)^i \dim_\bfk\Ext^i(M,K)=\dim_\bfk\Hom(M,K)-\dim_\bfk\Ext^1(M,K).
	\end{align}
	As in \cite[(A.1)-(A.2)]{LW19}, these forms descend to bilinear Euler forms on the Grothendieck groups $K_0(\cp^{\leq 1}(\Lambda^\imath))$ and $K_0(\mod(\iLa))$:
	\begin{align}
		\label{eq:Euler1}
		\langle\cdot,\cdot\rangle: K_0(\cp^{\leq 1}(\Lambda^\imath))\times K_0(\mod(\Lambda^\imath))\longrightarrow \Z,
		\\
		\label{eq:Euler2}
		\langle\cdot,\cdot\rangle: K_0(\mod(\Lambda^\imath))\times K_0(\cp^{\leq 1}(\Lambda^\imath))\longrightarrow \Z.
	\end{align}
	
	\begin{lemma}[\text{\cite[Lemma 4.3]{LW19}}]\label{lemma compatible of Euler form}
		We have
		\begin{itemize}
			\item[(i)]
			$\langle \E_i, M\rangle = \langle S_i,\res (M) \rangle_Q$ and $\langle M,\E_i\rangle =\langle \res(M), S_{\btau i} \rangle_Q$, for any $i\in Q_0$, $M\in\mod(\Lambda^{\imath})$;
			\item[(ii)] $\langle M,N\rangle=\frac{1}{2}\langle \res(M),\res(N)\rangle_Q$, for any $M,N\in\cp^{\leq 1}(\Lambda^\imath)$.
		\end{itemize}
	\end{lemma}

	Denote by $\widehat{\varrho}$ the exact functor of $\cd^b(\bfk Q)$ induced by $\varrho:\mod(\bfk Q)\rightarrow \mod(\bfk Q)$. 
	Let $\Gproj(\Lambda^\imath)$ be the subcategory of Gorenstein projective $\Lambda^\imath$-modules $G$ (that is, such that $\Ext^1_{\Lambda^\imath}(G,\Lambda^\imath)=0$). Then $\Gproj(\Lambda^\imath)$ is a Frobenius category with its stable category denoted by $\ul{\Gproj}(\Lambda^\imath)$; see \cite{Ha3}.
	
	\begin{theorem}[\text{\cite[Theorem 3.18]{LW19}}]
		\label{thm:sigma}
		Let $(Q, \btau)$ be an $\imath$quiver. Then $\cd^b(\bfk Q)/\Sigma \circ \widehat{\btau}$ is a triangulated orbit category \`a la Keller \cite{Ke05}, and we have the following triangulated equivalence
		\[
		\ul{\Gproj}(\Lambda^\imath)\simeq \cd_{sg}(\Lambda^{\imath})\simeq \cd^b(\bfk Q)/\Sigma \circ \widehat{\btau}.
		\]
	\end{theorem}
	
	As a corollary,	for any $M\in \cd_{sg}(\Lambda^{\imath})$, there exists a unique (up to isomorphisms) module $N\in \mod(\bfk Q)\subseteq \mod(\Lambda^{\imath})$ such that
	$M\cong N$ in $\cd_{sg}(\Lambda^{\imath})$; see  \cite[Corollary 3.21]{LW19}.

	%%%
	\subsection{$\imath$Hall algebras}
	\label{subsec:iHall}
	
	Let $\ch(\Lambda^\imath)$ be the  Hall algebra of $\mod(\Lambda^\imath)$. Let $\sqq$ be a fixed square root of $q$ in $\C$. We define the twisted Hall algebra $\widetilde{\ch}(\Lambda^\imath)$ to be the $\Q(\sqq^{1/2})$-algebra on the same vector space as $\ch(\Lambda^\imath)$ with twisted multiplication given by
	$$[M]*[N]=\sqq^{\langle \res M,\res N\rangle_{Q}}[M]\cdot [N].$$
	
	Let $\mathcal{I}$ be the subspace of $\widetilde{\ch}(\Lambda^\imath)$ spanned by all differences 
	\begin{align}
		\label{def:I}
		[M]-[N], \text{if $\res_\BH(M)=\res_\BH(N)$ and $M\cong N$ in $\cd_{sg}(\Lambda^\imath)$}.
	\end{align}
	Using the proof of \cite[Proposition 3.8]{LP24} (cf. \cite[Lemmas A11\& A12]{LW19}), we know that $\mathcal{I}$ is a two-sided ideal of $\widetilde{\ch}(\Lambda^\imath)$, and then we obtain a quotient algebra $\widetilde{\ch}(\Lambda^\imath)/\mathcal{I}$, and denote it by $\widehat{\ch}(\bfk Q,\varrho)$.

	Let	\begin{equation}
		\label{eq:Sca}
		\cs := \{ a[K] \in \widehat{\ch}(\bfk Q,\varrho) \mid a\in \Q(\sqq^{1/2})^\times, \pd K\leq1\}.
	\end{equation}
	By \cite{LW19,LW20,LP24}, the right localization of
	$\ch(\Lambda^\imath)/\mathcal{I}$ with respect to $\cs$ exists, and will be denoted by $\widetilde{\ch}(\bfk Q,\btau)$ or $\cs\cd\widetilde{\ch}(\Lambda^\imath)$, called the $\imath$Hall algebra (also called the twisted semi-derived Ringel-Hall algebra) of $\mod(\Lambda^\imath)$.

	For any $M\in\mod(\Lambda^\imath)$, we denote 
	$$\mathfrak{S}_M=\{[N]\mid \text{$\res_\BH(M)=\res_\BH(N)$ and $M\cong N$ in $\cd_{sg}(\Lambda^\imath)$}\}.$$
	Let $\mathcal{B}$ be the set of representatives of the sets $\mathfrak{S}_M$. Then $\mathcal{B}$ is a basis of $\widehat{\ch}(\bfk Q,\varrho)$. 
	
	%Let $S_i$ be the simple module of $\bfk Q$ corresponding to $i\in\I$. 
	
	We have $[\K_i]*[\K_j]=[\K_j]*[\K_i]=[\K_i\oplus \K_j]$ in $\widehat{\ch}(\bfk Q,\varrho)$ for any $i,j\in\I$; cf. \cite[Lemma 4.7]{LW19}. 
	For any $\alpha=(a_i)_{i\in\I}\in\N^{\I}$, we define in $\widehat{\ch}(\bfk Q,\btau)$: 
	$$[\K_\alpha]=[\oplus_{i\in\I}\K_i^{\oplus a_i}]=\prod_{i\in\I}[\K_i]^{a_i}.$$
	Similarly, one can define $[\K_\alpha]$ in $\widetilde{\ch}(\bfk Q,\btau)$ for $\alpha\in\Z^\I$. 
	
	\begin{lemma}
		\label{lem:KM-MK}
		For any $\alpha\in \N^\I$ and $M\in\mod(\Lambda^\imath)$, we have 
		$$[\K_\alpha]*[M]=\sqq^{\langle \varrho\alpha-\alpha,\widehat{\res (M)}\rangle_Q}[\K_\alpha\oplus M],\qquad [M]*[\K_\alpha]=\sqq^{\langle \widehat{\res (M)}, \alpha-\varrho\alpha\rangle_Q}[\K_\alpha\oplus M]$$ in $\widehat{\ch}(\bfk Q,\btau)$.
	\end{lemma}
	
	\begin{proof}
		It is enough to check it for $\alpha=\alpha_i$, that is $\K_\alpha=\K_i$, which follows from Lemma \ref{lemma compatible of Euler form} and the proof of \cite[Lemma A8]{LW19}.
	\end{proof}
	Moreover, we have
	\begin{align}
		\label{eq:KX=XK}
		[\K_\alpha]*[X]=\sqq^{( \btau\alpha, \widehat{X}  )_Q-( \alpha, \widehat{X}  )_Q}[X]*[\K_\alpha],\quad \forall X\in\mod(\bfk Q), \alpha\in\N^{\I}. 
	\end{align}

	\begin{lemma}[\text{cf. \cite[Proposition 4.9]{LW19}}]
		\label{basis-iHall}
		
		(1) The algebra $\widetilde{\ch}(\bfk Q,\btau)$ has a (Hall) basis given by
		\begin{align}
			\label{eq:Hallbasis-hat}
			\{[X]*[\K_\alpha]\mid X\in\mod(\bfk Q)\subseteq \mod(\Lambda^\imath), \alpha\in\Z^{\I}\}.
		\end{align}
		
		(2) The algebra $\widehat{\ch}(\bfk Q,\btau)$ has a (Hall) basis given by
		\begin{align}
			\label{eq:Hallbasis-tilde}
			\{[X]*[\K_\alpha]\mid X\in\mod(\bfk Q)\subseteq \mod(\Lambda^\imath), \alpha\in\N^{\I}\}.
		\end{align}
	\end{lemma}
	
	\begin{proof}
		The first statement is proved in \cite[Proposition 4.9]{LW19}. Let us prove the second one. For any $M\in\mod(\Lambda^\imath)$, by Theorem \ref{thm:sigma}, we know there exists a (unique up to isomorphism) $N\in\mod(\bfk Q)$, such that $M\cong N$ in $\cd_{sg}(\Lambda^\imath)$.
		Note that $\res_\BH(M)\cong L\oplus K$ with $L$ semisimple and $K$ projective (also injective). By definition, we assume $K=\K_\alpha$ for some $\alpha\in\N^\I$.   
		
		{\bf Claim}: $\res_\BH(N)=L$. 
		
		If the claim holds, then  $[M]=[N\oplus K]=[N\oplus \K_\alpha]=\sqq^{\langle \widehat{N},\varrho \alpha-\alpha \rangle_Q}[N]*[\K_\alpha]$ in $\widehat{\ch}(\bfk Q,\btau)$ by Lemma \ref{lem:KM-MK}. So $\{[X]*[\K_\alpha]\mid X\in\mod(\bfk Q)\subseteq \mod(\Lambda^\imath), \alpha\in\Z^{\I}\}$ spans $\widehat{\ch}(\bfk Q,\btau)$. 
		By using the natural projection $\widetilde{\ch}(\bfk Q,\btau)\rightarrow\widetilde{\ch}(\bfk Q,\btau)$ and (1), we know $\{[X]*[\K_\alpha]\mid X\in\mod(\bfk Q)\subseteq \mod(\Lambda^\imath), \alpha\in\Z^{\I}\}$ is linearly independent, and then the desired result follows. 
		
		Let us prove the claim. For $M$ and $N$, since $M\cong N$ in $\cd_{sg}(\Lambda^\imath)$, there exist short exact sequences $0\rightarrow P_M\rightarrow G\rightarrow M\rightarrow 0$ and $0\rightarrow P_N\rightarrow G\rightarrow N\rightarrow 0$  for some $G\in\Gproj(\Lambda^\imath)$ and $P_M,P_N\in\proj(\Lambda^\imath)$; see \cite[Theorem 11.5.1]{EJ}. By applying $\res_\BH$, one can easily see that $\res_\BH(N)\cong L$ by Lemma \ref{cor: res proj}. The claim is proved. 
	\end{proof}
	
	Let $\widetilde{\ct}(\bfk Q,\btau)$ be the subalgebra of $\widetilde{\ch}(\bfk Q,\btau)$ generated by $[\K_\alpha]$, $\alpha\in K_0(\mod(\bfk Q))$, which is a Laurent polynomial algebra in $[\E_i]$, for $i\in \I$. Similarly, one can define the subalgebra $\widehat{\ct}(\bfk Q,\btau)$ of $\widehat{\ch}(\bfk Q,\btau)$, which is a polynomial algebra in $[\E_i]$, for $i\in \I$.

	\begin{lemma}
		\label{lem:independent}
		For any $X,X',Y,Y',Z\in\mod(\Lambda^\imath)$, if $X'\in\mathfrak{S}_{X}$ and $Y'\in\mathfrak{S}_Y$, then we have 
		\begin{align}
			\sum_{[W]\in\mathfrak{S}_Z}G^W_{XY}=\sum_{[W]\in\mathfrak{S}_Z}G^W_{X'Y'}.
		\end{align}
	\end{lemma}
	
	\begin{proof}
		In the quotient algebra $\widehat{\ch}(\bfk Q,\varrho)=\widetilde{\ch}(\Lambda^\imath)/\mathcal{I}$, we have
		$[X]=[X']$, and $[Y]=[Y']$. Choose the basis $\mathcal{B}$ such that it includes $[Z]$. 
		Note that
		\begin{align}
			\label{eq:prod}
			[X]*[Y]=&\sqq^{\langle X,Y\rangle_Q}\sum_{[W]\in\Iso(\mod(\Lambda^\imath))} G_{XY}^W[W]=\sqq^{\langle X,Y\rangle_Q}\sum_{[B]\in\mathcal{B} } \big(\sum_{[W]\in\mathfrak{S}_B} G_{XY}^W\big)[B],
			\\\notag
			[X']*[Y']=&\sqq^{\langle X',Y'\rangle_Q}\sum_{[W]\in\Iso(\mod(\Lambda^\imath))} G_{X'Y'}^W[W]=\sqq^{\langle X',Y'\rangle_Q}\sum_{[B]\in\mathcal{B} } \big(\sum_{[W]\in\mathfrak{S}_B} G_{X'Y'}^W\big)[B],
		\end{align}
		which are equal. 
		By definition, we know $\langle X,Y\rangle_Q=\langle X',Y'\rangle_Q$. 
		Since $\mathcal{B}$ is a basis of $\widehat{\ch}(\bfk Q,\varrho)$, we have
		$\sum_{[W]\in\mathfrak{S}_B}G^W_{XY}=\sum_{[W]\in\mathfrak{S}_B}G^W_{X'Y'}$ 
		for any $[B]\in \mathcal{B}$, especially for $[Z]$.    
	\end{proof}
	
	\begin{proposition}[\text{\cite[Proposition 3.3]{LR24}}]
		\label{prop:iHallmult}
		For any $A,B\in\mod(\bfk Q)\subseteq \mod(\Lambda^\imath)$, we have in $\widetilde{\ch}(\bfk Q,\btau)$ (and also in $\widehat{\ch}(\bfk Q,\btau)$)
		\begin{align*}
			[A]*[B]=&
			\sum_{[L],[M],[N],[X]} \sqq^{\langle X,M\rangle-\langle \varrho X,M\rangle-\langle A,B\rangle} q^{\langle N,L\rangle} F_{N,L}^M F_{X,N}^AF_{ L,\varrho X}^{B} 
			\frac{a_L a_N a_X}{a_M} \cdot[K_{\widehat{X}}]*[M]
		\end{align*}
		where the sum is over $ [L],[M],[N],[X]\in\Iso(\mod(\bfk Q))$.
	\end{proposition}
	
	\begin{proof}
		The proof given in \cite[Proposition 3.3]{LR24} is for $\widetilde{\ch}(\bfk Q,\btau)$, which also works for $\widehat{\ch}(\bfk Q,\btau)$. 
	\end{proof}
	
	%Let $(Q,\btau)$ be a Dynkin $\imath$quiver ane $\tUi$ be its $\imath$quantum group.
	
	\begin{lemma}[cf. \text{\cite[Theorem 7.7]{LW19}}]
		\label{lem:Hall-iQG}
		Let $(Q, \btau)$ be a Dynkin $\imath$quiver. Then we have the following isomorphism $\widetilde{\psi}:\tUi|_{v=\sqq}\stackrel{\simeq}{\rightarrow} \widetilde{\ch}(\bfk Q,\btau)$ of $\Q({\sqq^{1/2}})$-algebras, which sends
		\begin{align}
			\label{eq:psi}
			B_i \mapsto \sqq^{-\frac{1}{2}}[S_i], \qquad 
			\tilde{k}_i \mapsto
			\begin{cases}
				[\K_i]&\text{if $\varrho i\neq i$},\\
				\sqq^{-1}[\K_{i}]&\text{if $\varrho i=i$}.
			\end{cases}
		\end{align}
	\end{lemma}
	
	\begin{proof}
		The result follows from \cite[Theorem 7.7]{LW19} by noting that the $\imath$quantum group $\tUi$ is different but isomorphic to the one given there. 
	\end{proof}
	
	Similarly, we have the isomorphism of $\Q(\sqq^{1/2})$-algebras
	$$\widehat{\psi}:\hUi|_{v=\sqq}\longrightarrow \widehat{\ch}(\bfk Q,\btau).$$

	As Example \ref{ex:QGvsiQG} shows, we can view quantum group $\tU$ to be an $\imath$quantum group; correspondingly, 
	as Example \ref{ex:diagquiver} shows, we can view the algebra $\Lambda$ in \eqref{eq:La} to be the $\imath$quiver algebra of $(Q^{\rm dbl},\swa)$.
	From Lemma \ref{lem:Hall-iQG}, we have the following result.
	\begin{lemma}[Bridgeland's Theorem reformulated] 
		\label{lem:bridgeland}
		Let $Q$ be a Dynkin quiver. Then we have the following isomorphism of $\Q(\sqq^{1/2})$-algebras
		\begin{align*}
			&\widetilde{\psi}:\tU|_{v=\sqq}\stackrel{\simeq}{\longrightarrow} \widetilde{\ch}(\bfk Q^{\rm dbl},\swa),
			\\
			E_i \mapsto \sqq^{-\frac{1}{2}}[S_i],&\qquad F_i\mapsto \sqq^{-\frac{1}{2}}[S_{i^\diamond}],
			\qquad
			K_i\mapsto [\K_{i^\diamond}],\qquad K_i'\mapsto [\K_i].
		\end{align*}
	\end{lemma}
	
	Similarly, we have the isomorphism of $\Q(\sqq^{1/2})$-algebras
	\[
	\widehat{\psi}:\hU|_{v=\sqq}\longrightarrow \widehat{\ch}(\bfk Q^{\rm dbl},\swa).
	\]

	%Then we have the following so called $\imath$Hall basis.
	%\begin{lemma}
	%\label{lem: basis of i-hall}
	%$\tMH$ is free as a right (respectively, left) $\tTL$-module, with a basis given by $[M]$ in $\Iso(\mod(\bfk Q))\subseteq \Iso(\mod(\Lambda^{\imath}))$.
	%\end{lemma}
	
	%\begin{definition}  \cite[Definition 4.11]{LW19}
	%  \label{def:reducedHall}
	%\red{Let $\bvs=(\vs_i)\in   (\Q(\sqq)^\times)^{\I}$ be such that $\vs_i=\vs_{\btau i}$ for each $i\in \I$. The \emph{reduced $\imath$Hall algebra associated to $(Q,\btau)$} \cite[Definition 4.11]{LW19}, denoted by $\ch(\bfk Q,\btau)$, is defined to be the quotient $\Q(\sqq^{1/2})$-algebra of $\widetilde{\ch}(\bfk Q,\btau)$ by the ideal generated by the central elements
	%\begin{align}
	%\label{eqn: reduce1}
	%[\E_i] +q \vs_i \; (\forall i\in \I \text{ with } \btau i=i), \text{ and }\; [\E_i]*[\E_{\btau i}] -\vs_i^2\; (\forall i\in \I \text{ with }\btau i\neq i).
	%\end{align}
	%}
	%\end{definition}
	
	%%
	\subsection{Generic $\imath$Hall algebras}
	\label{sub:generic}
	
	For a Dynkin $\imath$quiver $(Q,\btau)$,
	we recall the generic $\imath$Hall algebras defined in \cite[\S9.3]{LW19}.
	Recall that $\Phi^+$ is the set of positive roots.
	For any $\beta\in\Phi^+$, denote by $M_q(\beta)$ its corresponding indecomposable $\bfk Q$-module, i.e., $\dimv M_q(\beta)=\beta$.
	Let $\mathfrak{P}:=\mathfrak{P}(Q)$ be the set of functions $\lambda: \Phi^+\rightarrow \N$.
	Then the modules
	\begin{align}
		\label{def:Mlambda}
		M_q(\lambda):= \bigoplus_{\beta\in\Phi^+}\lambda(\beta) M_q(\beta),\quad \text{ for } \lambda\in\mathfrak{P},
	\end{align}
	provide a complete set of isoclasses of $\bfk Q$-modules.
	
	For $(\alpha,\nu),(\beta,\mu)\in\Z^\I\times\fp$, % (or $\N^\I\times\fp$), 
	there exists a polynomial $\boldsymbol{\varphi}^{\lambda,\gamma}_{\mu,\alpha;\nu,\beta}(v)\in\Z[v,v^{-1}]$ such that
	\[\big([\K_\alpha]\ast[M_q(\mu)]\big)\ast\big([\K_\beta]\ast[M_q(\nu)]\big)=\sum_{\lambda\in\fp,\gamma\in\Z^\I}\boldsymbol{\varphi}^{\lambda,\gamma}_{\mu,\alpha;\nu,\beta}({\sqq})[\K_\gamma]\ast[M_q(\lambda)]
	\]
	in $\widetilde{\ch}(\bfk Q,\btau)$. % (or $\widehat{\ch}(\bfk Q,\btau)$). 
	The generic $\imath$Hall algebra $\tMHg$ %(respectively, $\widehat{\ch}(Q,\btau)$) 
	is defined to be the  $\Q(v^{1/2})$-space with a basis $\{\K_\alpha*\fu_\lambda\mid \alpha\in\Z^\I,\lambda\in\fp\}$ % (respectively, $\{\K_\alpha*\fu_\lambda\mid \alpha\in\N^\I\lambda\in\fp$) 
	with multiplication
	\begin{align}
		\label{eq:generic-mult}
		(\K_\alpha*\fu_\mu)*(\K_\beta*\fu_\nu)=\sum_{\lambda\in\fp,\gamma\in\Z^\I}\boldsymbol{\varphi}^{\lambda,\gamma}_{\mu,\alpha;\nu,\beta}(v)\K_\gamma*\fu_\lambda.
	\end{align}
	
	For $\widehat{\ch}(\bfk Q,\btau)$, we can construct its generic version $\widehat{\ch}(Q,\btau)$ with a basis $\{\K_\alpha*\fu_\lambda\mid \alpha\in\N^\I,\lambda\in\fp\}$ similarly.
	
	From Lemma \ref{lem:Hall-iQG}, see also \cite[Theorem 9.8]{LW19}, we obtain 
	the isomorphisms of $\Q(v^{1/2})$-algebras
	\begin{alignat*}{2}
		&\widetilde{\psi}:\tUi\longrightarrow \widetilde{\ch}(Q,\btau),&\qquad &\widehat{\psi}:\hUi\longrightarrow \widehat{\ch}(Q,\btau),\\
		&\widetilde{\psi}:\tU\longrightarrow \widetilde{\ch}(Q^{\rm dbl},\swa),&\qquad &\widehat{\psi}:\hU\longrightarrow \widehat{\ch}(Q^{\rm dbl},\swa).
	\end{alignat*}

	\begin{remark}
		We abuse the notation $\K_\alpha$. To distinguish them, we use $[\K_\alpha]$ in $\widetilde{\ch}(\bfk Q,\btau)$,  and use $\K_\alpha$ in $\widetilde{\ch}( Q,\btau)$. The notation $\K_\alpha$ is also defined in $\tUi$, see \eqref{eq:bbKi}. However, they coincide with each other when we use $\imath$Hall algebra to realize the $\imath$quantum group; see Lemma \ref{lem:Hall-iQG}.
	\end{remark}

	\section{Dual canonical bases via Hall bases}\label{sec:dCB via Hall bases}
	
	%%%%%%%%%%%%%%%
	In this section, we shall use the natural bases (Hall bases) of $\imath$Hall algebras to construct dual canonical bases.
	
	\subsection{Dual canonical bases of Hall algebras}\label{dCB of HA subsec}
	
	In this subsection, we shall review the Lusztig's (dual) canonical basis of $\U^-$ via the generic twisted Hall algebra $\widetilde{\ch}(Q)$. 
	
	Recall the generators $\ce_i$, $\cf_i$ in \eqref{eq:Udj-gen}. %Let $$\xi^\pm: \U^{\pm}\rightarrow \U^{\pm},\qquad E_i\mapsto \ce_i,\quad F_i\mapsto \cf_i$$
	%be the automorphism of $\U^{\pm}$. 
	Lusztig constructed canonical basis of $\U^{\pm}$ via the generators $\ce_i,\cf_i$ ($i\in\I$). 
	%Recall that under the isomorphism $\psi:\hU\rightarrow\hU$ defined by $E_i\mapsto(v^{-1}-v)^{-1}E_i$, $F_i\mapsto (v-v^{-1})^{-1}F_i$, $K_i\mapsto K_i$, $K_i'\mapsto K_i'$, we obtain the standard presentation of $\hU$. 
	Let $\mathbf{B}^-$ be %the preimage of 
	Lusztig's canonical basis of $\U^-$ (via the generators $\cf_i$ ($i\in\I$)). % under the isomorphism $\xi^-$. 
	Following \cite{Ka91}, we define a Hopf pairing on $\U^-\otimes\U^+$ by
	\begin{align} 
		\label{eq:hopf}(F_i,E_j)_{K}=\delta_{ij}(v-v^{-1}).
	\end{align}
	
	For $b\in\mathbf{B}^-$, we denote by $\delta_b\in\U^+$ the dual basis of $b$ under the pairing $(\cdot,\cdot)_K$. We define a norm function $N:\Z^\I\rightarrow\Z$ by
	\[N(\alpha)=\frac{1}{2}(\alpha,\alpha)_Q-\eta(\alpha).\]
	where $\eta:\Z^\I\rightarrow\Z$ is the augmentation map defined by $\eta(\sum_ia_i\alpha_i)=\sum_ia_i$ (cf. \cite{GLS13,HL15}). The rescaled dual canonical basis of $\U^+$ is then defined to be
	\[\widetilde{\mathbf{B}}^+:=\{v^{\frac{1}{2}N(-\deg(b))}\delta_b\mid b\in\mathbf{B}^-\}.\]
	The rescaled dual canonical basis of $\U^-$ is defined as $\varphi(\widetilde{\mathbf{B}}^+)$, where $\varphi$ is the isomorphism $\U^+\rightarrow\U^-$ defined by $E_i\mapsto F_i$.

	\begin{example}
		For $\tU=\tU_v(\mathfrak{sl}_2)$, we have $\widetilde{\mathbf{B}}^+=\{E_1^n\mid n\in\N\}$ and $\widetilde{\mathbf{B}}^-=\{F_1^n\mid n\in\N\}$.
	\end{example}

	For a Dynkin quiver $Q$, we let $\widetilde{\ch}(\bfk Q)$ be the twisted Hall algebra of $\mod(\bfk Q)$ over $\Q(\sqq^{1/2})$, that is, 
	$$[M]\cdot [N]=\sqq^{\langle M,N\rangle_Q}\sum_{[L]}\frac{|\Ext^1(M,N)_L|}{|\Hom(M,N)|}[L],\quad \forall M,N\in\mod(\bfk Q).$$
	Recall that $\mathfrak{P}:=\mathfrak{P}(Q)$ be the set of functions $\lambda: \Phi^+\rightarrow \N$. 
	We let $\widetilde{\ch}(Q)$ be the generic Hall algebra of $Q$. It is a $\Q(v^{\frac{1}{2}})$-vector space with basis $\{\fu_\lambda\mid\lambda\in\mathfrak{P}\}$, and the multiplication is defined by 
	\[\fu_\mu\cdot\fu_\nu=\sum_{\lambda}g^{\lambda}_{\mu,\nu}(v)\fu_\lambda\]
	where $g^{\lambda}_{\mu,\nu}(v)\in\Z[v,v^{-1}]$ is the polynomial such that
	\[
	[M_q(\mu)]\cdot[M_q(\lambda)]=\sum_{\lambda\in\mathfrak{P}}g^{\lambda}_{\mu,\nu}(\sqq)[M_q(\lambda)]\in\widetilde{\ch}(\bfk Q).
	\]
	The following result is well-known (note the coefficient on $\fu_{\alpha_i}$).
	
	\begin{proposition}[\cite{Rin90}]
		There exists an isomorphism of $\Q(v^{1/2})$-algebras
		\begin{align*}
			\psi^+:\U^+\longrightarrow \widetilde{\ch}(Q),\quad E_i\mapsto v^{-\frac{1}{2}}\fu_{\alpha_i},\quad\forall i\in\I.
		\end{align*}
	\end{proposition}
	
	Similarly, there exists an isomorphism of $\Q(v^{1/2})$-algebras 
	\[
	\psi^-:\U^-\longrightarrow \widetilde{\ch}(Q),\quad F_i\mapsto v^{-\frac{1}{2}}\fu_{\alpha_i},\quad\forall i\in\I.
	\] 
	The Hopf pairing in \eqref{eq:hopf} can then be transferred to $\widetilde{\ch}(\bfk Q)\otimes\widetilde{\ch}(\bfk Q)$: it is given by
	\begin{align}
		\label{eq:hopf-Hall}
		([M],[N])_K=\delta_{M,N}|\aut(M)|.
	\end{align}
	We can also define a similar paring on $\widetilde{\ch}(Q)\otimes \widetilde{\ch}(Q)$.

	\begin{example}
		For $\tU_v(\mathfrak{sl}_2)$, the dual canonical basis of $\widetilde{\ch}(Q)$ is $\{v^{-\frac{m^2}{2}}\fu_{m\alpha_1}\mid m\in\N\}$.
	\end{example}
	
	Recall from Proposition~\ref{QG bar-involution def} that there is a bar-involution on $\U^+$ defined by $v^{\frac{1}{2}}\mapsto v^{-\frac{1}{2}}$ and $E_i\mapsto E_i$, which is an anti-automorphism. Using the isomorphism $\psi^+$, this can be transferred to the generic Hall algebra $\widetilde{\ch}(Q)$: it is characterized by $\ov{\fu_{\alpha_i}}=v^{-1}\fu_{\alpha_i}$. 
	
	For $\lambda\in\mathfrak{P}$, we define a rescaled Hall basis of $\widetilde{\ch}(Q)$ by
	%\red{$\mathfrak{U}_\lambda$}
	\begin{equation}\label{eq:HA element H_lambda}
		\mathfrak{U}_\lambda=v^{-\dim\End_{\bfk Q}(M_q(\lambda))+\frac{1}{2}\langle M_q(\lambda),M_q(\lambda)\rangle_Q}\fu_\lambda.
	\end{equation}
	
	The definition of $\mathfrak{U}_\lambda$ is justified by the fact that $\mathfrak{U}_\beta$ is bar-invariant for each $\beta\in\Phi^+$. To see this, let us give a useful lemma.
	
	\begin{lemma}[cf. \cite{Rin3}]
		\label{lem:short exact sequence}
		Let $k$ be any field, and $Q$ be of type ADE.
		For any indecomposable $kQ$-module $M$, there exist two indecomposable $kQ$-modules $M_1$ and $M_2$ such that the following hold:
		\begin{align}
			\label{eqn: split cond 1}
			&\Hom(M_1,M_2)=0=\Hom(M_2,M_1),\qquad \dim\Ext^1_{kQ}(M_2,M_1)=1;
			\\\label{eqn: split cond 2}
			&\text{there exists a short exact sequence
			} 0\rightarrow M_1\rightarrow M\rightarrow M_2\rightarrow0.
		\end{align}
		Furthermore, in this case, we have $\Ext^1_{kQ}(M_1,M_2)=0$.
	\end{lemma}
	
	\begin{proof}
		From \cite[Proposition 5]{Rin3}, there exist two indecomposable $kQ$-modules $M_1$ and $M_2$ such that $\Hom(M_1,M_2)=0=\Hom(M_2,M_1)$, and \eqref{eqn: split cond 2} holds. Note that $M_1,M_2,M$ are indecomposable. Then $0\rightarrow M_1\rightarrow M\rightarrow M_2\rightarrow0$ is non-split, and $\dim \Ext^1_{kQ}(M_2,M_1)\neq0$.  Recall the symmetric bilinear form defined in \S\ref{subsec:Quiver-represent}. Since $Q$ is of type ADE, we have $(\dimv M_1,\dimv M_2)_Q\geq-1$. %\lutodo{($=0,1,-1$)}. 
		Since $\dim \Ext^1_{kQ}(M_2,M_1)\neq0$, we have $\Ext^1_{kQ}(M_1,M_2)=0$ and $\dim\Ext^1_{kQ}(M_2,M_1)=1$. Therefore, we have the desired properties of $M_1$ and $M_2$.
	\end{proof}
	
	\begin{proposition}
		\label{HA bar of H_beta}
		For each $\beta\in\Phi^+$, we have $\ov{\mathfrak{U}_\beta}=\mathfrak{U}_\beta$.  
	\end{proposition}
	
	\begin{proof}
		For $\beta\in\Phi^+$ we have $\dim\End_{\bfk Q}(M_q(\beta))=1$ and $\langle M_q(\beta),M_q(\beta)\rangle_Q=1$, so by definition 
		$\mathfrak{U}_\beta=v^{-\frac{1}{2}}\fu_\beta$. 
		We prove $\ov{\mathfrak{U}_\beta}=\mathfrak{U}_\beta$ by induction on dimension of $M_q(\beta)$.
		
		If $\beta=\alpha_i$, it is obvious. 
		In general, if $M=M_q(\beta)$ is an indecomposable $\bfk Q$-module, then we can find $\bfk Q$-modules $M_1$ and $M_2$ as in Lemma~\ref{lem:short exact sequence}. Then we have 
		\begin{align*}
			[M_2]\cdot[M_1]&=\sqq^{-1}[M_1\oplus M_2]+(\sqq-\sqq^{-1})[M],\\
			[M_1]\cdot[M_2]&=[M_1\oplus M_2].
		\end{align*}
		Let us write $M_1=M_q(\beta_1)$ and $M_2=M_q(\beta_2)$ for some $\beta_1,\beta_2\in\Phi^+$. In $\widetilde{\ch}(Q)$ we get
		\begin{align*}
			\mathfrak{U}_\beta&=\frac{v^{\frac{1}{2}}\mathfrak{U}_{\beta_2}\cdot \mathfrak{U}_{\beta_1}-v^{-\frac{1}{2}}\mathfrak{U}_{\beta_1}\cdot \mathfrak{U}_{\beta_2}}{v-v^{-1}}.
		\end{align*}
		As both $\mathfrak{U}_{\beta_1}$ and $\mathfrak{U}_{\beta_2}$ are bar-invariant due to our induction assumption, it is easy to see that $\mathfrak{U}_\beta$ is also bar-invariant.
	\end{proof}
	
	Let $\prec$ be the partial order on $\mathfrak{P}$ defined by orbit closure: we say that $\lambda\prec\mu$ if the orbit $\mathfrak{O}_{M_q(\lambda)}$ is contained in the closure of $\mathfrak{O}_{M_q(\mu)}$; see \cite{Lus90} or \cite[\S1.6]{DDPW}. The result of Proposition~\ref{HA bar of H_beta} can be generalized as follows. 
	
	\begin{proposition}
		\label{HA bar of H_lambda}
		For each $\lambda\in\mathfrak{P}$, one can write
		\[\ov{\mathfrak{U}_\lambda}-\mathfrak{U}_\lambda\in\sum_{\lambda\prec\mu}\Z[v,v^{-1}]\mathfrak{U}_\mu.\]
	\end{proposition}
	\begin{proof}
		To prove this we can specialize to the Hall algebra $\widetilde{\ch}(\bfk Q)$. The bar-involution of $\widetilde{\ch}(\bfk Q)$ is similarly defined by setting $\ov{[S_i]}=\sqq^{-1}[S_i]$. Now choose an ordering $\beta_1,\dots,\beta_l$ of positive roots such that
		\begin{align*}
			\Hom_{\bfk Q}(M_q(\beta_a),M_q(\beta_b))&=0,\quad \text{for $a>b$},\\
			\Ext^1_{\bfk Q}(M_q(\beta_a),M_q(\beta_b))&=0,\quad \text{for $a\leq b$}.
		\end{align*}
		For any $\lambda\in\mathfrak{P}$, we can write
		\begin{align}
			\label{eq:lambda-decop}
			[M_q(\lambda)]&=\sqq^{\sum_{a<b}\lambda_a\lambda_b\langle\beta_a,\beta_b\rangle_Q}[M_q(\beta_1)^{\oplus\lambda_1}]\cdot[M_q(\beta_2)^{\oplus\lambda_2}]\cdots[M_q(\beta_l)^{\oplus\lambda_l}]
		\end{align}
		where $\lambda_k=\lambda(\beta_k)$. Note that $\ov{[M_q(\beta_k)^{\oplus\lambda_k}]}=\sqq^{-\lambda_k^2}[M_q(\beta_k)^{\oplus\lambda_k}]$ for each $k$. 
		Applying the bar-involution to \eqref{eq:lambda-decop} gives \begin{equation}\label{HA bar of H_lambda-1}
			\begin{aligned}
				\ov{[M_q(\lambda)]}&=\sqq^{-\sum_{a<b}\lambda_a\lambda_b\langle\beta_a,\beta_b\rangle_Q-\sum_k\lambda_k^2}[M_q(\beta_l)^{\oplus\lambda_l}]\cdot[M_q(\beta_{l-1})^{\oplus\lambda_{l-1}}]\cdots[M_q(\beta_1)^{\oplus\lambda_1}]\\
				&=\sqq^{\sum_{a>b}\lambda_a\lambda_b\langle\beta_a,\beta_b\rangle_Q-\sum_{a<b}\lambda_a\lambda_b\langle\beta_a,\beta_b\rangle_Q-\sum_k\lambda_k^2}[M_q(\lambda)]+\sum_{[N]}f_N[N],
			\end{aligned}
		\end{equation}
		for some $f_N\in\Z[\sqq,\sqq^{-1}]$. 
		By the definition of multiplication in $\widetilde{\ch}(\bfk Q)$, we see that each class $[N]$ appearing in the RHS of (\ref{HA bar of H_lambda-1})  with $f_N\neq0$ satisfies $\dimv N=\dimv M_q(\lambda)$ and $\mathfrak{O}_{M_q(\lambda)}\subseteq\ov{\mathfrak{O}}_N$. The result then follows by noting
		\begin{align*}
			&-\dim\End_{\bfk Q}(M_q(\lambda))+\frac{1}{2}\langle M_q(\lambda),M_q(\lambda)\rangle_Q\\
			&=-\frac{1}{2}\dim\End_{\bfk Q}(M_q(\lambda))-\frac{1}{2}\Ext^1_{\bfk Q}(M_q(\lambda),M_q(\lambda))
			\\
			&=-\frac{1}{2}\sum_k\lambda_k^2-\frac{1}{2}\sum_{a<b}\lambda_a\lambda_b\langle\beta_a,\beta_b\rangle_Q+\frac{1}{2}\sum_{a>b}\lambda_a\lambda_b\langle\beta_a,\beta_b\rangle_Q
		\end{align*}
		and the fact that $\langle M_q(\lambda),M_q(\lambda)\rangle_Q=\langle N,N\rangle_Q$ for any $N$ such that $f_N\neq 0$.
	\end{proof}
	
	Since for each $\lambda\in\mathfrak{P}$ there are only finitely many elements $\mu$ such that $\lambda\prec\mu$, we can apply Lusztig's Lemma (\cite[Theorem 1.1]{BZ14}) to extract a bar-invariant basis from $\mathfrak{U}_\lambda$. More precisely, we have the following theorem.
	
	\begin{theorem}
		\label{thm:dualCB-U+}
		For each $\lambda\in\mathfrak{P}$, there exists a unique element $\mathfrak{C}_\lambda\in\widetilde{\ch}(Q)$ such that $\ov{\mathfrak{C}_\lambda}=\mathfrak{C}_\lambda$ and 
		\[\mathfrak{C}_\lambda-\mathfrak{U}_\lambda\in\sum_{\mu\in\mathfrak{P}}v^{-1}\Z[v^{-1}]\mathfrak{U}_\mu.\]
		Moreover, $\mathfrak{C}_\lambda$ satisfies
		\[\mathfrak{C}_\lambda-\mathfrak{U}_\lambda\in\sum_{\lambda\prec\mu}v^{-1}\Z[v^{-1}]\mathfrak{U}_\mu.\]
	\end{theorem}
	
	We shall see that $\mathfrak{C}_\lambda$ coincides with the dual canonical basis of $\widetilde{\ch}(Q)$. For this, let us recall Lusztig's elementary algebraic construction of canonical basis of $\U^{\pm}$ by using the PBW basis \cite{Lus90}. Using the isomorphisms $\psi^{\pm}$, we transfer his construction to $\widetilde{\ch}(Q)$ (see also \cite[\S11.6]{DDPW}). 
	
	%\red{Fix a $Q$-admissible ordering $\beta_1,\dots,\beta_N$ of the roots in $\Phi^+$, that is,
	%$$\Hom_{\bfk Q}(M_q(\beta_i),M_q(\beta_j))\neq0 \text{ only if }i\leq j.$$}
	
	For any $\lambda\in\fp$, we define a PBW basis of $\widetilde{\ch}(Q)$ by (compare with \cite[eq. (11.4.1)]{DDPW}, the differences come from the Hall multiplication and the map $\psi^{\pm}$)
	\[
	\mathfrak{E}_\lambda=v^{\dim \End(M_q(\lambda))-\frac{1}{2}\dim M_q(\lambda)}\frac{\fu_\lambda}{a_\lambda(v)},
	\]
	here $a_\lambda(v)\in\Z[v,v^{-1}]$ is the polynomial such that $a_\lambda(\sqq)=|\aut(M_q(\lambda))|$; see e.g. \cite[Lemma 10.19]{DDPW}. The corresponding canonical basis $\mathfrak{B}_\lambda\in \U^+$ then satisfies
	\[
	\psi^+(\mathfrak{B}_\lambda)-\mathfrak{E}_\lambda\in\sum_{\mu\prec\lambda}v^{-1}\Z[v^{-1}]\mathfrak{E}_\mu.
	\]
	Let $\mathfrak{E}_\lambda^*$ ($\lambda\in\fp$) be the dual PBW basis with respect to the paring \eqref{eq:hopf-Hall}. Then one can check that %(using Theorem~\ref{dCB of U^+ by L})
	\begin{align}
		\label{eq:TE}
		\mathfrak{U}_\lambda=v^{\frac{1}{2}n_\lambda}\mathfrak{E}_\lambda^*, \text{ where $n_\lambda=N(\sum_k\lambda(\beta_k)\beta_k)$}.
	\end{align}
	
	Let $\mathfrak{B}_\lambda^*$ ($\lambda\in\fp$) be the dual canonical basis of $\U^+$ with respect to the paring \eqref{eq:hopf}. Then we have 
	\begin{align*}
		\psi^+(\mathfrak{B}_\lambda^*)\in \mathfrak{E}_\lambda^*+\sum_{\lambda\prec\mu}v^{-1}\Z[v^{-1}]\mathfrak{E}_\mu^*,\quad \ov{\mathfrak{B}_\lambda^*}=v^{n_\lambda}\mathfrak{B}_\lambda^*.
	\end{align*}
	In view of \eqref{eq:TE} and Theorem~\ref{thm:dualCB-U+}, this means $\psi^+(v^{\frac{1}{2}n_\lambda}\mathfrak{B}_\lambda^*)=\mathfrak{C}_\lambda$. %The following proposition then follows Theorem~\ref{dCB of U^+ by L}.

	\subsection{Integral forms}
	
	Set $\cz=\Z[v^{1/2},v^{-1/2}]$ and $\cz_\sqq=\Z[\sqq^{1/2},\sqq^{-1/2}]$. We define the integral form of Hall algebras and $\imath$Hall algebras in the following. 
	
	By definition of Hall algebras, we know $\widetilde{\ch}(\Lambda^\imath)$ can be defined over $\cz_\sqq$, which is a free $\cz_\sqq$-module  with a basis given by isoclasses $[M]$ of all modules $M$. This is denoted by $\widetilde{\ch}(\Lambda^\imath)_{\cz_\sqq}$. By considering $\mathcal{I}$ in \eqref{def:I} the submodule (ideal) of $\widetilde{\ch}(\Lambda^\imath)_{\cz_\sqq}$, we can define $\widehat{\ch}(\bfk Q,\varrho)_{\cz_\sqq}$ over $\cz_\sqq$. Similarly, using \eqref{eq:Sca}, we can define the $\imath$Hall algebra $\widetilde{\ch}(\bfk Q,\varrho)_{\cz_\sqq}$ over $\cz_\sqq$. Using the proof of Lemma \ref{basis-iHall}, we know $\widehat{\ch}(\bfk Q,\varrho)_{\cz_\sqq}$ and $\widetilde{\ch}(\bfk Q,\varrho)_{\cz_\sqq}$ are free $\cz_\sqq$-module with their bases given by \eqref{eq:Hallbasis-hat} and \eqref{eq:Hallbasis-tilde} respectively. Furthermore, the algebra $\tMHg_{\cz}$ (respectively, $\widehat{\ch}(Q,\btau)_{\cz}$) is defined to be the free $\cz$-module with a basis $\{\K_\alpha*\fu_\lambda\mid \alpha\in\Z^\I, \lambda\in\fp\}$ (respectively, $\{\K_\alpha*\fu_\lambda\mid \alpha\in\N^\I, \lambda\in\fp\}$ ), and the multiplication as in \eqref{eq:generic-mult}.  
	
	Using the  isomorphism of $\Q(v^{1/2})$-algebras
	$$\widetilde{\psi}:\tUi\rightarrow \widetilde{\ch}(Q,\btau),\qquad \widehat{\psi}:\hUi\rightarrow \widehat{\ch}(Q,\btau),$$
	we can define the integral form of $\tUi$ (resp. $\hUi$) as the preimage of $\widehat{\ch}(Q,\varrho)_\cz$ (resp. $\widehat{\ch}(Q,\varrho)_\cz$), which is denoted by $\tUi_\cz$ (resp. $\hUi_\cz$). Both of them are free $\cz$-modules. 
	We shall see that $\tUi_\cz$ and $\hUi_\cz$ are independent of orientations; see Corollary \ref{cor:FT-ihall-integral}. Note that $\hUi_\cz=\tUi_\cz\cap\hUi$.

	\begin{remark}
		By definition of quantum groups, we know $\tU$ can be defined as $\cz$-algebra generated by $E_i,F_i,K_i,K_i'$ subject to relations \eqref{eq:KK}--\eqref{eq:serre2}. Similarly, $\tUi$ can be defined as the subalgebra of $\tU$ over $\cz$. However, from Proposition \ref{prop:braid1} and Theorem \ref{thm:Ti}, we can see that these naive integral forms are different to the integral forms $\tU_\cz$ and $\tUi_\cz$ defined above.   
	\end{remark}

	Recall that an ordering $\beta_1,\dots,\beta_N$ of the positive roots in $\Phi^+$ is called to be {\em $Q$-admissible} if $\Ext^1_{\bfk Q}(M_q(\beta_r),M_q(\beta_t))\neq 0$ implies $r>t$ (or $\Hom_{\bfk Q}(M_q(\beta_r),M_q(\beta_t))\neq0$ only if $r\leq t$). 
	A sequence $i_1,\dots,i_t$ in $\ci$ is called {\em $\imath$-admissible} if $i_j$ is a sink in $r_{i_{j-1}}\cdots r_{i_1}(Q)$, for each $1\leq j\leq t$.
	
	\begin{lemma}[{\cite[Theorem 8.3]{LW21a}}]
		\label{thm:i-seq}
		Let $(Q,\btau)$ be a Dynkin $\imath$quiver. Then there is an $\imath$-admissible sequence of length $N_\imath$, $\{i_1,\dots,i_{\ell}\}$, such that
		\[
		\beta_{1}, \btau(\beta_1), \beta_{2}, \btau(\beta_2), \dots,\beta_{\ell },\btau(\beta_{\ell})
		\]
		is a $Q$-admissible sequence of the roots in $\Phi^+$. (By convention the redundant $\btau(\beta_j)$ is omitted here and below whenever $\btau(\beta_j)=\beta_j$.) Moreover, $\bs_{i_1}\cdots \bs_{i_{\ell}}$ is a reduced expression of $\omega_0$ in $W_\btau$ (and it becomes a reduced expression of $\omega_0 \in W$ if each $r_{i_k}$ is replaced by products of simple reflections as in \eqref{def:simple reflection}). In particular, $\ell$ is the length of $\omega_0$ in $W_\btau$.
	\end{lemma}
	
	The $\imath$-admissible sequence in Lemma~\ref{thm:i-seq} is called {\em complete}.
	
	\begin{proposition}\label{prop:iQG integral form by PBW}
		For a Dynkin $\imath$quiver $(Q,\btau)$, we consider the complete $\imath$-admissible sequence $i_1,\cdots,i_\ell$ to give a reduced expression of $\omega_0$ in $W_\varrho$. Then $\tUi_\cz$ is the (free) $\cz$-module generated by the PBW basis obtained in Proposition \ref{prop:PBW-Ui} by using this $\imath$-admissible sequence. 
	\end{proposition}
	
	\begin{proof}
		Denote by ${'\tUi_\cz}$ the (free) $\cz$-module generated by the PBW basis obtained in Proposition \ref{prop:PBW-Ui}. Using the proof of \cite[Lemma 8.3]{LW21a}, we have $\widetilde{\psi}(B_\beta)= v^{-1/2}\fu_\beta$  for any $\beta\in\Phi^+$. In general, for $B^\ba=\prod_{\beta\in\Phi^+}B_\beta^{a_\beta}$ (with the order given by the $\imath$-admissible sequence $i_1,\dots,i_\ell$), we have $\widetilde{\psi}(B^{\ba})\in(v^{-1/2})^{f_\ba}\fu_{\lambda}+\sum_{\tilde{\lambda},\mu} \Z[v,v^{-1}]\K_{\mu}*\fu_{\tilde{\lambda}}$, where $f_\ba\in\Z$, $\lambda,\tilde{\lambda}:\Phi^+\rightarrow\N$, such that $\lambda(\beta)=a_\beta$, and $\mu\in\N^\I$ which is nonzero. By induction on the $|\ba|=\sum_\beta a_\beta$, we can see that  $B^\ba\in\tUi_\cz$, and  then $'\tUi_\cz\subseteq \tUi_\cz$. Conversely, by the leading term of $\widetilde{\psi}(B^\ba)$, we know $\fu_\lambda\in \widetilde{\psi}('\tUi_\cz)$ by induction, and then $\tUi_\cz\subseteq {}'\tUi_\cz$ by definition. Then  $\tUi_\cz={'\tUi_\cz}$.   
	\end{proof}
	
	\begin{corollary}
		\label{cor:bar-integral-Ui}
		The bar-involution of $\tUi$ (resp. $\hUi$) preserves $\tUi_\cz$ (resp. $\hUi_\cz$).
	\end{corollary}
	
	\begin{proof}
		By Proposition \ref{prop:iQG integral form by PBW}, we can see that $\tUi_\cz$ is the $\cz$-algebra generated by $\{\K_\alpha,B_\beta\mid \alpha\in\Z^\I,\beta\in\Phi^+\}$. It is enough to prove $\ov{B_\beta}\in\tUi_\cz$ for all $\beta\in\Phi^+$, which follows from Lemma \ref{iQG bar-involution def} immediately.
	\end{proof}

	The bar-involution of $\widehat{\ch}(Q,\varrho)$ (and also $\widetilde{\ch}(Q,\varrho)$) can be defined similarly to $\widetilde{\ch}(Q)$ by setting
	\[\ov{\fu_{\alpha_i}}=v^{-1}\fu_{\alpha_i},\quad \ov{\K_{\alpha_i}}=\K_{\alpha_i}.\]
	
	By Corollary \ref{cor:bar-integral-Ui}, we have the following lemma. 
	\begin{lemma}\label{lem:bar on iHall integral}
		The bar-involution of $\widetilde{\ch}(Q,\varrho)$ (resp. $\widehat{\ch}(Q,\varrho)$) preserves $\widetilde{\ch}(Q,\varrho)_\cz$ (resp. $\widehat{\ch}(Q,\btau)_\cz$).
	\end{lemma}

	As quantum groups $\tU$ can be viewed as $\imath$quantum groups, the integral forms $\widetilde{\ch}(Q^{\rm dbl},\swa)_\cz$ and $\tU_\cz$ can be defined. Then $\tU_\cz$ is the free $\cz$-module generated by the PBW basis obtained in Proposition \ref{prop:PBW-QG}, which is preserved by the bar-involution of $\tU$. 
	
	\subsection{Dual canonical bases of $\imath$Hall algebras}\label{subsec:dCB of iHA}
	
	Inspired by the results of \S\ref{dCB of HA subsec}, we construct a bar-invariant basis in $\widehat{\ch}(Q,\varrho)$. Recall that there is a natural inclusion of $\Q(v^{\frac{1}{2}})$-vector spaces
	\begin{equation}\label{eq:inclusion of HA to iHA}
		\iota:\widetilde{\ch}(Q)\hookrightarrow \widehat{\ch}(Q,\varrho),\quad \fu_\lambda\mapsto\fu_\lambda,\forall\lambda\in\mathfrak{P}.
	\end{equation}
	The element $\mathfrak{U}_\lambda$ defined by \eqref{eq:HA element H_lambda} can therefore be regarded as an element in $\widehat{\ch}(Q,\varrho)$. Let $\widehat{\ct}(Q,\varrho)$ be the subalgebra generated by $\K_\alpha$, $\alpha\in\N^\I$. Then $\widehat{\ch}(Q,\varrho)$ has a basis given by 
	\[
	\{\K_\alpha\ast\fu_\lambda\mid\alpha\in\N^\I,\lambda\in\mathfrak{P}\}.
	\]
	
	Consider the quotient maps 
	\begin{align}
		\label{eq:iHA quotient}
		\pi:&\widehat{\ch}(Q,\varrho)\longrightarrow \widetilde{\ch}(Q),\quad \K_\alpha\ast\fu_\lambda\mapsto \begin{cases}
			\fu_\lambda&\text{if $\alpha=0$},\\
			0&\text{otherwise}.
		\end{cases}
	\end{align}
	Note that it is a homomorphism of $\Q(v^{1/2})$-algebras; cf. Proposition \ref{prop:iHallmult}.
	
	Following \cite{BG17}, we define an action of $\widehat{\ct}(Q,\varrho)$ on $\widehat{\ch}(Q,\varrho)$ by
	\[\K_\alpha\diamond\fu_\lambda=v^{\frac{1}{2}(\alpha-\varrho\alpha,\,\dimv{M_q(\lambda)})_Q}\K_\alpha\ast\fu_\lambda.\]
	This definition has the advantage that for any $\alpha\in\N^\I$ and $\lambda\in\mathfrak{P}$,
	\begin{equation}\label{eq:diamond action and bar}
		\ov{\K_\alpha\diamond \fu_\lambda}=\K_\alpha\diamond\ov{\fu_\lambda}.
	\end{equation}
	
	Similar to the case of Hall algebras, we define a partial order on $\N^\I\times\mathfrak{P}$: we say $(\alpha,\lambda)\prec(\beta,\mu)$ if $\alpha+\btau(\alpha)+\dimv M_q(\lambda)=\beta+\btau(\beta)+\dimv M_q(\mu)$ and either $\alpha\prec\beta$ (i.e. $\alpha\neq\beta$ and $\beta-\alpha\in\N^\I$) or $\alpha=\beta$ and $\lambda\prec\mu$. The following result is an analogue of Proposition~\ref{HA bar of H_lambda} for $\imath$Hall algebras.
	
	\begin{proposition}\label{iHA bar of H_lambda}
		For $\alpha\in\N^\I$ and $\lambda\in\mathfrak{P}$, we have 
		\[\ov{\K_\alpha\diamond \mathfrak{U}_\lambda}-\K_\alpha\diamond \mathfrak{U}_\lambda\in\sum_{(\alpha,\lambda)\prec(\beta,\mu)}\Z[v,v^{-1}]\cdot \K_\beta\diamond \mathfrak{U}_\mu.\]
	\end{proposition}
	\begin{proof}
		In view of \eqref{eq:diamond action and bar} and the definition of the partial order, it suffices to prove the proposition for $\alpha=0$. In this case we can write (by Lemma~\ref{lem:bar on iHall integral})
		\begin{align}\label{eq:iHA bar of H_lambda-1}
			\ov{\fu_\lambda}=\sum_{\beta,\mu}f_{\beta,\mu}(v)\cdot \K_\beta\ast \fu_\mu
		\end{align}
		where $f_{\beta,\mu}\in\Z[v,v^{-1}]$ and $f_{\beta,\mu}\neq 0$ only if $\dimv M_q(\lambda)=\beta+\varrho(\beta)+\dimv M_q(\mu)$. 
		
		We first note that for $f_{\beta,\mu}\neq 0$ only if $(0,\lambda)\prec(\beta,\mu)$. To see this, consider the quotient map \eqref{eq:iHA quotient}. It is easy to see that $\pi$ commutes with the bar-involutions of $\widehat{\ch}(Q,\varrho)$ and $\widetilde{\ch}(Q)$, so by Proposition~\ref{HA bar of H_lambda}, we find
		\[
		\pi(\ov{\fu_\lambda}-\fu_\lambda)=\ov{\pi(\fu_\lambda)}-\pi(\fu_\lambda)\in\sum_{\lambda\prec\mu}\Z[v,v^{-1}]\cdot \fu_\mu.
		\]
		On the other hand, by our convention $\pi(\ov{\fu_\lambda})=\sum_\mu f_{0,\mu}\cdot\pi(\fu_\mu)$, so our assertion follows.
		
		From \eqref{eq:iHA bar of H_lambda-1} we then conclude that (set $\bd=\dimv M_q(\lambda)$ and $\bd'=\dimv M_q(\mu)$)
		\begin{align*}
			\ov{\mathfrak{U}_\lambda}=\sum_{(0,\lambda)\prec(\beta,\mu)}f_{\beta,\mu}(v)\cdot v^{\dim\End_{\bfk Q}(M_q(\lambda))+\dim\End_{\bfk Q}(M_q(\mu))-\frac{1}{2}\langle \bd,\bd\rangle_Q-\frac{1}{2}\langle\bd',\bd'\rangle_Q-\frac{1}{2}(\beta-\varrho\beta,\bd')_Q} \K_\beta\diamond \mathfrak{U}_\mu
		\end{align*}
		It remains to show that
		\begin{align}\label{eq:iHA bar of H_lambda-2}
			-\frac{1}{2}\langle \bd,\bd\rangle_Q-\frac{1}{2}\langle\bd',\bd'\rangle_Q-\frac{1}{2}(\beta-\varrho\beta,\bd')_Q\in\Z.
		\end{align}
		For this we note that $\bd=\bd'+\beta+\varrho\beta$, so 
		\begin{align*}
			\text{LHS of \eqref{eq:iHA bar of H_lambda-2}}&=-\frac{1}{2}\langle \bd,\bd\rangle_Q-\frac{1}{2}\langle\bd-\beta-\varrho\beta,\bd-\beta-\varrho\beta\rangle_Q-\frac{1}{2}(\beta-\varrho\beta,\bd-\beta-\varrho\beta)_Q\\
			&=-\langle\bd,\bd\rangle_Q+\langle\varrho\beta,\bd-\beta-\varrho\beta\rangle_Q+\langle\bd,\varrho\beta\rangle_Q+\frac{1}{2}\langle\beta,\varrho\beta\rangle_Q-\frac{1}{2}\langle\varrho\beta,\beta\rangle_Q\\
			&=-\langle\bd,\bd\rangle_Q+\langle\varrho\beta,\bd-\beta-\varrho\beta\rangle_Q+\langle\bd,\varrho\beta\rangle_Q
		\end{align*}
		In the last equation we have used the fact that $\varrho$ is an involution of the quiver $Q$, so $\langle\beta,\varrho\beta\rangle_Q=\langle\varrho\beta,\beta\rangle_Q$. This completes the proof.
	\end{proof}
	
	Again, as for a fixed pair $(\alpha,\lambda)$ there are only finitely many pairs $(\beta,\mu)$ such that $(\alpha,\lambda)\prec(\beta,\mu)$, Lusztig's Lemma \cite[Theorem 1.1]{BZ14} is applicable and we obtain a bar-invariant basis for $\widehat{\ch}(Q,\varrho)$. We also note the $\diamond$-action preserves this basis, more precisely, we have the following theorem. 
	
	\begin{theorem}\label{iHA dCB theorem}
		For each $\alpha\in\N^\I$ and $\lambda\in\mathfrak{P}$, there exists a unique element $\mathfrak{L}_{\alpha,\lambda}\in\widehat{\ch}(Q,\varrho)$ such that $\ov{\mathfrak{L}_{\alpha,\lambda}}=\mathfrak{L}_{\alpha,\lambda}$ and
		\[
		\mathfrak{L}_{\alpha,\lambda}-\K_\alpha\diamond \mathfrak{U}_\lambda\in\sum_{(\beta,\mu)}v^{-1}\Z[v^{-1}]\cdot \K_\beta\diamond \mathfrak{U}_\mu.
		\]
		Moreover, $\mathfrak{L}_{\alpha,\lambda}$ satisfies 
		\[
		\mathfrak{L}_{\alpha,\lambda}-\K_\alpha\diamond \mathfrak{U}_\lambda\in\sum_{(\alpha,\lambda)\prec(\beta,\mu)}v^{-1}\Z[v^{-1}]\cdot \K_\beta\diamond \mathfrak{U}_\mu,
		\]
		and $\mathfrak{L}_{\alpha,\lambda}=\K_\alpha\diamond \mathfrak{L}_{0,\lambda}$.
	\end{theorem}

	Because of the last property, we often write $\mathfrak{L}_\lambda:=\mathfrak{L}_{0,\lambda}$ for $\lambda\in\mathfrak{P}$ and use $\{\K_\alpha\diamond \mathfrak{L}_\lambda\mid \alpha\in\N^\I,\lambda\in\mathfrak{P}\}$ to denote the basis constructed in Theorem~\ref{iHA dCB theorem}. This is called the dual canonical basis of $\widehat{\ch}(Q,\varrho)$. The dual canonical basis of $\widetilde{\ch}(Q,\varrho)$ is defined to be $\{\K_\alpha\diamond \mathfrak{L}_\lambda\mid \alpha\in\Z^\I,\lambda\in\mathfrak{P}\}$.
	
	As the  diagonal type of $\imath$quantum groups, the dual canonical basis of $\widetilde{\ch}(Q^{\rm dbl},\swa)$ is also defined. 
	
	%%%%%%%%%%%%%%%%%%%%%%%%%%%%%%%%
	\subsection{Reflection functors}\label{subsec:iHA reflection}
	Let $\ell\in Q_0$ be a sink. Denote by $s_\ell Q$ the quiver obtained from $Q$ by reversing all arrows ending at $\ell$. Let $\mathscr{R}_\ell^+:\mod(\bfk Q)\rightarrow \mod(\bfk s_\ell Q)$ be the BGP reflection. 
	
	For an $\imath$quiver $(Q,\btau)$, define for any $\ell\in Q_0$
	\begin{align*}
		r_\ell=\begin{cases}
			s_\ell s_{\btau \ell} & \text{ if }\btau \ell\neq \ell,
			\\
			s_\ell &\text{ if }\btau \ell=\ell.
		\end{cases}
	\end{align*}
	In this way, we define $r_\ell Q$ for any sink $\ell\in Q_0$, and then $(r_\ell Q,\varrho)$ is also an $\imath$quiver. Denote by $r_\ell\Lambda^\imath$ its $\imath$quiver algebra.
	In \cite{LW21a}, we define a BGP type reflection 
	$F_\ell^+:\mod(\Lambda^\imath)\rightarrow\mod(r_\ell \Lambda^\imath)$, which restricts to $\mod(\bfk Q)$ is 
	\begin{align*}
		\begin{cases}
			\mathscr{R}_\ell^+  :\mod(\bfk Q)\rightarrow\mod(\bfk r_\ell Q)&\text{ if }\btau\ell=\ell,
			\\
			\mathscr{R}_\ell^+ \mathscr{R}_{\varrho\ell}^+=\mathscr{R}_{\varrho\ell}^+ \mathscr{R}_\ell^+ :\mod(\bfk Q)\rightarrow\mod(\bfk r_\ell Q)&\text{ if }\btau\ell\neq\ell.
		\end{cases} 
	\end{align*}

	For any $i\in \I$, denote by $S_i$ (respectively, $S_i'$) the simple $\bfk Q$-module (respectively, $\bfk r_\ell Q$-module), denote by $\K_i$ (respectively, $\K_i'$) the
	generalized simple $\iLa$-module (respectively, $\bs_\ell\iLa$-module). We similarly define $[\K_\alpha],[\K_\beta']$ for $\alpha\in K_0(\mod(\bfk Q))$, $\beta\in K_0(\mod(\bfk(r_\ell Q)))$
	in the $\imath$Hall algebras (where $\ell$ is a sink of $Q$).
	
	Recall the root lattice $\Z^{\I}=\Z\alpha_1\oplus\cdots\oplus\Z\alpha_n$, and we have an isomorphism of abelian groups $\Z^\I\rightarrow K_0(\mod(\bfk Q))$, $\alpha_i\mapsto \widehat{S_i}$. This isomorphism induces the action of the reflection $s_i$ on $K_0(\mod(\bfk Q))$.
	Thus for $\alpha\in K_0(\mod (\bfk Q))$ and $i\in\I$, $[\K_{s_i\alpha}] \in \iH(\bfk Q,\btau)$ is well defined. Similarly, we have $[\K'_{s_i\alpha}]\in \iH(\bfk r_\ell Q,\btau)$.
	
	Let $\ct=\{X\in \mod(\Lambda^\imath) \mid \Hom(X,S_\ell\oplus S_{\btau \ell})=0\}$. Then any $\bfk Q$-module $M$ is of the form 
	\begin{align}
		\label{eq:M-form}
		M\cong \begin{cases}
			X_M\oplus S_\ell^{\oplus a}\oplus S_{\btau \ell}^{\oplus b}& \text{ if }\varrho\ell\neq \ell,
			\\
			X_M\oplus S_\ell^{\oplus a} &\text{ if }\btau \ell=\ell,
		\end{cases}
	\end{align}
	where $X_M\in \ct\cap \mod(\bfk Q)$. By \cite[Proposition 4.4]{LW21a}, there exists an isomorphism 
	\[\Gamma_\ell:\iH(\bfk Q,\btau)\rightarrow\iH(\bfk r_\ell Q,\btau),\]
	which induces an isomorphism
	\[
	\Gamma_\ell: \iH( Q,\btau)\rightarrow\iH( r_\ell Q,\btau).
	\]
	By \cite{LW21a}, we have the following commutative diagram:
	\begin{equation}\label{eq:defT}
		\begin{tikzcd}
			\tUi  \ar[r,"\tTT_\ell"] \ar[d,swap,"\widetilde{\psi}_Q"] & \tUi \ar[d,"\widetilde{\psi}_{r_\ell Q}"]\\
			\tMHg \ar[r,"\Gamma_\ell"] &  \widetilde{\ch}(r_\ell Q,\btau)
		\end{tikzcd}
	\end{equation}
	
	Since $\widetilde{\TT}_\ell$ commutes with the bar-involution of $\tUi$ by Lemma \ref{lem:QGbraid-bar}, we also have $\Gamma_\ell$ commutes with the bar-involution of $\widetilde{\ch}(Q,\varrho)$. 
	
	\begin{lemma}\label{Gamma_l on [M]}
		Let $M$ be a $\bfk Q$-module of the form \eqref{eq:M-form}. Let $\bd=(d_i)_{i\in\I}$ be the dimension vector of $M$.  Then we have     
		\begin{align*}
			\Gamma_\ell([M])&=\begin{cases}
				\sqq^{(a-b)(d_\ell-d_{\varrho \ell})} [(\K'_{\ell})^{\oplus a}\oplus(\K'_{\varrho\ell})^{\oplus b}]^{-1} *[F_{\ell}^+(X_M)\oplus (S_{\varrho \ell}')^{\oplus a}\oplus (S_{ \ell}')^{\oplus b}] &\text{ if }\varrho \ell\neq \ell,
				\\
				[\K_\ell^{\oplus a}]^{-1}*[F_\ell^{+}(X_M)\oplus (S_{\ell}')^{\oplus a}]&\text{ if }\varrho\ell=\ell.
			\end{cases}
		\end{align*}
	\end{lemma}
	
	\begin{proof}
		We assume $Q$ is connected in the following. By \cite[Proposition 4.4]{LW21a}, we know
		\begin{align}
			\Gamma_{\ell}([M])&= [F_{\ell}^+(M)], \quad \forall M\in\ct,
			\label{eqn:reflection 1}\\
			\Gamma_{\ell}([S_\ell])&=\begin{cases}
				\sqq [\K_\ell']^{-1}\ast [S_{\varrho \ell}'] &\text{if $\varrho \ell\neq\ell$},\\
				[\K_\ell']^{-1}\ast [S_{\varrho \ell}']      &\text{if $\varrho \ell=\ell$},
			\end{cases}
			\label{eqn:reflection 2}\\
			\Gamma_{\ell}([S_{\btau \ell}])&=\sqq [\K'_{\btau \ell}]^{-1} \ast [S_{\ell}'], \quad \text{if $\varrho \ell\neq \ell$}.
			\label{eqn:reflection 3}
		\end{align}
		If $Q$ has only one $\varrho$-orbit, that is, $Q_0=\{\ell,\varrho\ell\}$, then one can prove the result by induction and the following formulae: 
		\[
		[S_\ell] \ast [S_{\ell}^{\oplus a}] = \sqq^{-a} [S_\ell^{\oplus(a+1)}] + (\sqq^a-\sqq^{-a})[S_\ell^{\oplus(a-1)}]\ast[\K_\ell],\quad \text{if $\varrho\ell=\ell$};
		\]
		and
		\begin{align*}
			[S_j]*[S_\ell^{\oplus a}\oplus S_{\varrho\ell}^{\oplus b}]=\begin{cases}
				\sqq^{-a} [S_\ell^{\oplus (a+1)}\oplus S_{\varrho\ell}^{\oplus b}]+\sqq^{1-b}(q^b-1) [\K_\ell]*[S_\ell^{\oplus a}\oplus S_{\varrho\ell}^{\oplus (b-1)}]       &\text{ if $j=\ell$},\\
				\sqq^{-b} [S_\ell^{\oplus a}\oplus S_{\varrho\ell}^{\oplus (b+1)}]+\sqq^{1-a}(q^a-1) [\K_{\varrho\ell}]*[S_\ell^{\oplus (a-1)}\oplus S_{\varrho\ell}^{\oplus b}] &\text{ if $j=\varrho\ell$}.
			\end{cases}
		\end{align*}
		
		In general, since $\ell$ is a sink, there exists at least one arrow $\alpha:j\rightarrow\ell$ in $Q_0$, so there exists a string module $Y$ with its string $\ell\xleftarrow{\alpha} j\xrightarrow{\varepsilon_j}\btau j$. Then $Y$ and $\K_j$ are in $\mathcal{T}$, which admit short exact sequences
		\[
		0\rightarrow S_\ell\rightarrow Y\rightarrow \K_j\rightarrow 0,\quad  
		0\rightarrow S_{\varrho\ell}\rightarrow \varrho Y\rightarrow \K_{\varrho j}\rightarrow 0. 
		\]
		
		We only prove the case $\varrho \ell\neq \ell$. For $M=X_M\oplus S_\ell^{\oplus a}\oplus S_{\varrho \ell}^{\oplus b}$, we have a short exact sequence
		\[
		0\rightarrow M\rightarrow X_M\oplus Y^{\oplus a}\oplus \varrho Y^{\oplus b}\rightarrow \K_j^{\oplus a}\oplus \K_{\varrho j}^{\oplus b}\rightarrow 0.
		\]
		By \cite[Theorem 4.3]{LW21a}, we know that
		\begin{align*}
			\Gamma_{\ell}([M])=&\sqq^{\langle \res(\K_j^{\oplus a}\oplus \K_{\varrho j}^{\oplus b}),\res(M)\rangle_Q} q^{-\langle \K_j^{\oplus a}\oplus \K_{\varrho j}^{\oplus b}),M\rangle} [F_{\ell}^+(\K_j^{\oplus a}\oplus \K_{\varrho j}^{\oplus b}))]^{-1}\\
			&\ast [F_{\ell}^+(X_M)\oplus F_\ell^+(Y^{\oplus a})\oplus F_\ell^+(\varrho Y^{\oplus b})]
			\\
			=&\sqq^{(a+b)\langle S_j \oplus S_{\varrho j},\res(M)\rangle_Q}q^{-a\langle S_j, \res(M)\rangle_Q-b\langle S_{\varrho j},\res(M)\rangle_Q} [F_{\ell}^+(\K_j^{\oplus a}\oplus \K_{\varrho j}^{\oplus b}))]^{-1}
			\\
			&\ast [F_{\ell}^+(X_M)\oplus F_\ell^+(Y^{\oplus a})\oplus F_\ell^+(\varrho Y^{\oplus b})]
			\\
			=&\sqq^{(b-a)\langle S_j ,\res(M)\rangle_Q+(a-b)\langle S_{\varrho j},\res(M)\rangle_Q} [F_{\ell}^+(\K_j^{\oplus a}\oplus \K_{\varrho j}^{\oplus b}))]^{-1}
			\\
			&\ast [F_{\ell}^+(X_M)\oplus F_\ell^+(Y^{\oplus a})\oplus F_\ell^+(\varrho Y^{\oplus b})].
		\end{align*}
		Similar to \cite[Proposition 4.4]{LW21a}, we have 
		\begin{align*}
			[F_{\ell}^+(\K_j^{\oplus a}\oplus \K_{\varrho j}^{\oplus b}))]^{-1}&=[(\K'_{r_\ell(\alpha_j)})^{\oplus a}\oplus (\K'_{r_\ell(\alpha_{\varrho j})})^{\oplus b}]^{-1},
		\end{align*}
		and 
		\begin{align*}
			&[F_{\ell}^+(X_M)\oplus F_\ell^+(Y^{\oplus a})\oplus F_\ell^+(\varrho Y^{\oplus b})]
			\\
			=&[F_{\ell}^+(X_M)\oplus (\K'_{r_\ell(\alpha_j)-\alpha_\ell})^{\oplus a}\oplus (S'_{\varrho \ell})^{\oplus a}\oplus(\K'_{r_\ell(\varrho \alpha_j)-\alpha_{\varrho\ell}})^{\oplus b}\oplus (S'_{ \ell})^{\oplus b}]
			\\
			=&\sqq^{(a-b)\langle r_\ell(\alpha_j)-\alpha_\ell,\res(F_\ell^+(X_M)) \oplus (S'_{\varrho \ell})^{\oplus a}\oplus (S'_{ \ell})^{\oplus b}\rangle_{Q'}+(b-a)\langle r_\ell(\varrho \alpha_j)-\alpha_{\varrho \ell},\res(F_\ell^+(X_M)) \oplus (S'_{\varrho \ell})^{\oplus a}\oplus (S'_{ \ell})^{\oplus b}\rangle_{Q'}}
			\\
			&[(\K'_{r_\ell(\alpha_j)-\alpha_\ell})^{\oplus a}\oplus (\K'_{r_\ell(\varrho \alpha_j)-\alpha_{\varrho\ell}})^{\oplus b}]*[F_{\ell}^+(X_M)\oplus (S'_{\varrho \ell})^{\oplus a}\oplus (S'_{ \ell})^{\oplus b}].
		\end{align*}
		By \cite[Lemma 4.2]{LW21a}, we have
		\begin{align*}
			& \langle r_{\ell}(\alpha_j),\res(F_\ell^+(X_M))\rangle_{Q'}=\langle S_j,X_M\rangle_Q,\qquad \langle r_{\ell}(\varrho \alpha_j),\res(F_\ell^+(X_M))\rangle_{Q'}=\langle S_{\varrho j},X_M\rangle_Q,
			\\
			&\langle S'_{\ell},\res(F_\ell^+(X_M))\rangle_{Q'}=-\langle S_\ell,\res(X_M)\rangle_{Q},\qquad \langle S'_{\varrho\ell},\res(F_\ell^+(X_M))\rangle_{Q'}=-\langle S_{\varrho \ell},\res(X_M)\rangle_{Q},
			\\
			&\langle r_{\ell}(\alpha_j), (S'_{\varrho \ell})^{\oplus a}\oplus(S'_{ \ell})^{\oplus b} \rangle_{Q'}=0= \langle r_{\ell}(\varrho \alpha_j),(S'_{\varrho \ell})^{\oplus a}\oplus(S'_{ \ell})^{\oplus b} \rangle_{Q'},
		\end{align*}
		and these imply
		\begin{align*}
			\Gamma_\ell([M])
			=&\sqq^{(b-a)\langle S_\ell',\res(F_\ell^+(X_M)) \oplus (S'_{\varrho \ell})^{\oplus a}\oplus (S'_{ \ell})^{\oplus b} \rangle_{Q'}+(a-b)\langle S_{\varrho\ell}',\res(F_\ell^+(X_M)) \oplus (S'_{\varrho \ell})^{\oplus a}\oplus (S'_{ \ell})^{\oplus b} \rangle_{Q'}}
			\\
			&[\K'_{-a\alpha_\ell-b\alpha_{\varrho\ell}}] *[F_{\ell}^+(X_M)\oplus (S'_{\varrho \ell})^{\oplus a}\oplus (S'_{ \ell})^{\oplus b}]
			\\
			=&\sqq^{(b-a)\langle S_\ell',\res(F_\ell^+(X_M))  \rangle_{Q'}+(a-b)\langle S_{\varrho\ell}',\res(F_\ell^+(X_M)) \rangle_{Q'}+(a-b)^2}[\K'_{-a\alpha_\ell-b\alpha_{\varrho\ell}}] 
			\\
			&*[F_{\ell}^+(X_M)\oplus (S'_{\varrho \ell})^{\oplus a}\oplus (S'_{ \ell})^{\oplus b}]
			\\
			=&\sqq^{(a-b)\langle \widehat{S_\ell}-\widehat{S_{\varrho\ell}},\widehat{\res(X_M)} \rangle_Q+(a-b)^2} [\K'_{-a\alpha_\ell-b\alpha_{\varrho\ell}}] *[F_{\ell}^+(X_M)\oplus (S'_{\varrho \ell})^{\oplus a}\oplus (S'_{ \ell})^{\oplus b}].
		\end{align*} 
		Let $\bd'=(d'_i)_{i\in\I}$ be the dimension vector of $X_M$. Then we have $d_\ell=d'_\ell+a$, $d_{\varrho\ell}=d'_{\varrho\ell}+b$, and $\langle S_\ell,\res(X_M)\rangle_Q=d'_\ell$, $\langle S_{\varrho\ell},\res(X_M)\rangle_Q=d'_{\varrho\ell}$.
		From these the desired formula follows.
	\end{proof}
	
	For $M\in\mod(\bfk Q)$, let us write
	\[\mathfrak{U}_{[M]}=\sqq^{-\dim\End_{\bfk Q}(M)+\frac{1}{2}\langle M,M\rangle_Q}[M].\]
	This is the image of $\mathfrak{U}_\lambda$ under the specialization map $v\mapsto\sqq$, where $\lambda\in\mathfrak{P}$ is such that $M_q(\lambda)=M$. We can also define the $\diamond$-action of $\widetilde{\ct}(\bfk Q,\varrho)$ on $\widetilde{\ch}(\bfk Q,\varrho)$ by 
	\[
	\K_\alpha\diamond[M]=v^{\frac{1}{2}(\alpha-\varrho\alpha,\,\dimv M)_Q}\K_\alpha\ast[M].
	\]
	
	\begin{proposition}
		\label{Gamma_l on H_lambda}
		Let $M\in\mod(\bfk Q)$ be of the form \eqref{eq:M-form}, and define
		\[N=\begin{cases}
			F_{\ell}^+(X_M)\oplus (S'_{\varrho \ell})^{\oplus a}\oplus (S'_{ \ell})^{\oplus b}& \text{if $\varrho \ell\neq\ell$},\\
			F_\ell^{+}(X_M)\oplus (S_{\ell}')^{\oplus a}& \text{if $\varrho \ell=\ell$}.
		\end{cases}\]
		Then we have 
		\begin{align}
			\Gamma_\ell(\mathfrak{U}_{[M]})=\begin{cases}
				[\K_{a\alpha_\ell+b\alpha_{\varrho\ell}}]^{-1}\diamond \mathfrak{U}_{[N]} &\text{ if }\varrho \ell\neq \ell,
				\\
				[\K_{a\alpha_\ell}]^{-1}\diamond \mathfrak{U}_{[N]}&\text{ if }\varrho\ell=\ell.
			\end{cases}
		\end{align}
	\end{proposition}
	
	\begin{proof}
		Let $\bd=(d_i)_{i\in\I}$ be the dimension vector of $M$. Let $N=F_\ell^+(X_M)$ and $\bd'=(d'_i)_{i\in\I}$ be its dimension vector. Then $\dimv X_M=r_\ell(\bd')$.
		
		We first prove the case $\varrho\ell=\ell$. It is enough to show that 
		\begin{align*}
			&-\dim \End_{\bfk Q}(X_M\oplus S_\ell^{\oplus a})+\frac{1}{4}(X_M\oplus S_\ell^{\oplus a},X_M\oplus S_\ell^{\oplus a})_{Q}
			\\
			=&-\dim\End_{\bfk Q'}(F_\ell^{+}(X_M)\oplus (S_{\ell}')^{\oplus a})+\frac{1}{4}(F_\ell^{+}(X_M)\oplus (S_{\ell}')^{\oplus a},F_\ell^{+}(X_M)\oplus (S_{\ell}')^{\oplus a})_{Q'}.
		\end{align*}
		
		First,
		\begin{align*}
			\dim\End_{\bfk Q}(X_M\oplus S_\ell^{\oplus a})&=\dim\End_{\bfk Q}(X_M)+a^2+a\langle S_\ell,\res(X_M)\rangle_Q,\\
			\dim\End_{\bfk Q'}(F_\ell^{+}(X_M)\oplus (S_{\ell}')^{\oplus a})&=\dim\End_{\bfk Q'}(F_\ell^+(X_M))+a^2+a\langle\res(F_\ell^+(X_M)),S_\ell\rangle_{Q'}
		\end{align*}
		while
		\begin{align*}
			&(X_M\oplus S_\ell^{\oplus a},X_M\oplus S_\ell^{\oplus a})_{Q}-(F_\ell^{+}(X_M)\oplus (S_{\ell}')^{\oplus a},F_\ell^{+}(X_M)\oplus (S_{\ell}')^{\oplus a})_{Q'}
			\\
			=&(r_\ell(\bd')+a\alpha_\ell,r_\ell(\bd')+a\alpha_\ell)_{Q}-(\bd'+a\alpha_\ell,\bd'+a\alpha_\ell)_{Q'}
			\\
			=&4a(\alpha_\ell,\dimv(X_M))_Q.
		\end{align*}
		It then suffices to note that $\langle S_\ell,\res(X_M)\rangle_Q=-\langle S_\ell,\res(F_\ell^+(X_M))\rangle_{Q'}$ and $\End_{\bfk Q}(X_M)\cong\End_{\bfk Q'}(F_\ell^+(X_M))$.
		
		Next consider the case $\varrho\ell\neq\ell$. We have 
		\begin{align*}
			&-\dim\End_{\bfk Q}(X_M\oplus S_\ell^{\oplus a}\oplus S_{\varrho\ell}^{\oplus b})+\dim\End_{\bfk Q'}(F_\ell^+(X_M)\oplus (S_{\varrho\ell}')^{\oplus a}\oplus (S_{\ell}')^{\oplus b})\\
			&=-\langle a\alpha_\ell+b\alpha_{\varrho\ell},\dimv(X_M)\rangle_Q+\langle \dimv(F_\ell^+(X_M)), a\alpha_{\varrho\ell}+b\alpha_\ell\rangle_{Q'}\\
			&=\langle a\alpha_\ell+b\alpha_{\varrho\ell},\dimv(F_\ell^+(X_M))\rangle_{Q'}+\langle \dimv(F_\ell^+(X_M)), a\alpha_{\varrho\ell}+b\alpha_\ell\rangle_{Q'}
		\end{align*}
		and
		\begin{align*}
			&(X_M\oplus S_\ell^{\oplus a}\oplus S_{\varrho\ell}^{\oplus b},X_M\oplus S_\ell^{\oplus a}\oplus S_{\varrho\ell}^{\oplus b})_{Q}-(F_\ell^+(X_M)\oplus (S_{\varrho\ell}')^{\oplus a}\oplus (S_{\ell}')^{\oplus b})_{Q'}\\
			&=2(a+b)(r_\ell(\bd'),\alpha_\ell+\alpha_{\varrho\ell})_Q=-2(a+b)(\bd',\alpha_\ell+\alpha_{\varrho\ell})_{Q'}.
		\end{align*}   
		In view of Lemma~\ref{Gamma_l on [M]}, these give us
		\[\Gamma_\ell(\mathfrak{U}_{[M]})=\sqq^{(a-b)(\bd_\ell-\bd_{\varrho\ell}+\frac{1}{2}(\langle\bd',\alpha_{\varrho\ell}-\alpha_\ell\rangle_{Q'}-\langle\alpha_{\varrho\ell}-\alpha_\ell,\bd'\rangle_{Q'}))}[\K_{a\alpha_\ell+b\alpha_{\varrho\ell}}]^{-1} \ast \mathfrak{U}_{[N]}.\]
		It then remains to prove that
		\begin{align*}
			&(a-b)(\bd_\ell-\bd_{\varrho\ell}+\frac{1}{2}(\langle\bd',\alpha_{\varrho\ell}-\alpha_\ell\rangle_{Q'}-\langle\alpha_{\varrho\ell}-\alpha_\ell,\bd'\rangle_{Q'}))
			\\
			&=\frac{1}{2}(-a\alpha_\ell-b\alpha_{\varrho\ell}+a\alpha_{\varrho\ell}+b\alpha_\ell,\bd'+a\alpha_{\varrho\ell}+b\alpha_\ell)_{Q'},
		\end{align*}
		or equivalently
		\[(d_\ell-d_{\varrho\ell})+\langle\alpha_\ell-\alpha_{\varrho\ell},\bd'\rangle_{Q'}=a-b.\]
		This follows easily from the fact that $\bd'=r_\ell(\bd-a\alpha_\ell-b\alpha_{\varrho\ell})$.
	\end{proof}
	
	The computation in Proposition~\ref{Gamma_l on H_lambda} can be made generically, and this implies $\Gamma_\ell$ preserves the rescaled Hall bases, that is, it
	maps the basis $\{\K_\alpha\diamond \mathfrak{U}_\lambda\mid\alpha\in\Z^\I,\lambda\in\mathfrak{P}\}$ for $\widetilde{\ch}(Q,\varrho)$ to the one for $ \widetilde{\ch}(r_\ell Q,\varrho)$.
	
	\begin{corollary}\label{iHA reflection functor of dCB}
		The isomorphism $\Gamma_\ell:\widetilde{\ch}(Q,\varrho)\rightarrow \widetilde{\ch}(r_\ell Q,\varrho)$ maps the dual canonical basis of $\widetilde{\ch}(Q,\varrho)$ to the dual canonical basis of $\widetilde{\ch}(r_\ell Q,\varrho)$.
	\end{corollary}
	\begin{proof}
		This follows immediately from Proposition~\ref{Gamma_l on H_lambda} and the fact that $\Gamma_\ell$ commutes with the bar-involution of $\widetilde{\ch}(Q,\varrho)$.
	\end{proof}
	
	\begin{corollary}
		For any $\beta\in\Phi^+$, the element $\mathfrak{U}_\beta$ belongs to the dual canonical basis of $\widetilde{\ch}(Q,\varrho)$.
	\end{corollary}
	\begin{proof}
		This follows from Corollary~\ref{iHA reflection functor of dCB} and the fact that $\fu_\beta$ can be obtained from some $\fu_{\alpha}$, $\alpha\in\Delta^+$, by successive actions of reflection functors (cf. \cite[\S 8]{LW21a}).
	\end{proof}

	%%%%%%%%%
	\section{Fourier transforms of $\imath$Hall algebras}\label{sec:FT of iHA}
	
	In this section, all Hall algebras and $\imath$Hall algebras are defined over the field of complex numbers $\C$. We shall establish the Fourier transforms of $\imath$Hall algebras. First, we shall use the machinery of \cite{SV99} to realize Hall algebras and $\imath$Hall algebras via functions. % 

	%%%%%%%%%%%%%%%%%%%%
	\subsection{$\imath$Hall algebras via functions}\label{iHA by function subsec}
	
	%We introduce a slight modification based on the construction given in \S\ref{subsec:QV function iHA}. 
	For each $\bd=(d_i)_i\in\N^{\I}$, let $V=\oplus_{i\in\I} V_i$ be an $\I$-graded vector space of dimension vector $\bd$ over $\bfk$. We impose and fix $V_{\varepsilon_i}:V_i\rightarrow V_{\varrho i}$ for each $i\in\I$ such that $V_{\varepsilon_i}\circ V_{\varepsilon_{\varrho i}}=0$ for any $i\in\I$. (That means, $V^{i}=V^{\varrho i}:=(V_i,V_{\varrho i},V_{\varepsilon_i},V_{\varepsilon_{\varrho i}})$ is an $\BH_i$-module.) In this way, $V$ is an $\BH$-module, which is fixed in the following. We then define vector spaces
	\begin{align*}
		E_\bd=&E_{Q,\bd}=E_V=E_{Q,V}:=\prod_{\alpha\in Q_1} \Hom_{\bfk}(V_{\ts(\alpha)},V_{\tt(\alpha)}),
		\\
		\ov{E}_V=&\ov{E}_{Q,V}:=\{(x_\alpha:V_{\ts(\alpha)}\rightarrow V_{\tt(\alpha)})_{\alpha\in Q_1}\in E_\bd\mid x_{\varrho \alpha}\circ V_{\varepsilon_{\ts(\alpha)}}=V_{\varepsilon_{\tt(\alpha)}}\circ x_{\alpha},\forall \alpha\in Q_1\}.
	\end{align*}
	It is obvious that $\ov{E}_V$ is a linear subspace of $E_\bd$.
	Therefore, a point $x=V(x)=(x_\alpha)_\alpha$ of $\ov{E}_V$ determines a representation of $\Lambda^\imath$. We denote
	\[a_x:=|\Aut(V(x))|.\] 
	The group $G_V:=\prod_{i\in\I_\varrho} \Aut_{\BH_i}(V^i)\subset \prod_{i\in\I}\mathrm{GL}(V_i)$ acts on $\ov{E}_V$ by conjugation
	\[
	(g_i)_{i\in\I}\cdot (V_\alpha)_\alpha=(g_{\tt(\alpha)}V_\alpha g_{\ts(\alpha)}^{-1})_\alpha,
	\]
	and then the $G_V$-orbits $\mathfrak{O}_x$ correspond bijectively to the isoclasses $[V(x)]$ of representations of $\Lambda^\imath$ with the ground $\BH$-module $V$. 
	
	Let $\cf(\Lambda^\imath)_V$ denote the set of all $G_V$-invariant functions on $\ov{E}_V$ with value in $\C$. It has a basis given by the characteristic functions $\chi_{\mathfrak{O}_x}$ on $G_V$-orbits $\mathfrak{O}_x$. We define
	\begin{align}
		\cf(\Lambda^\imath)=\bigoplus_{[V]\in\Iso(\mod(\BH))}\cf(\Lambda^\imath)_V.
	\end{align}
	For convenience, we set for any dimension vector $\bd$,
	\begin{align}
		\cf(\Lambda^\imath)_\bd=\bigoplus_{[V]\in\Iso(\mod(\BH)),\dimv V=\bd}\cf(\Lambda^\imath)_V, 
	\end{align}
	and then $$\cf(\Lambda^\imath)=\bigoplus_{\bd\in\N^{\I}}\cf(\Lambda^\imath)_\bd.$$
	
	Let $V',V''$ be $\BH$-modules. Denote $\bd',\bd''$ the dimension vectors of $V',V''$ respectively, and $\bd=\bd'+\bd''$. Then for $f\in\cf(\Lambda^\imath)_{V'}$ and $g\in\cf(\Lambda^\imath)_{V''}$, define their convolution product $f\star g\in\cf(\Lambda^\imath)_\bd$ by setting
	\begin{align*}
		(f\star g)(z)=\frac{\sqq^{\langle \bd',\bd''\rangle_Q}}{|G_{V'}|\cdot |G_{V''}|}\sum_{x,y,\alpha,\beta} \frac{a_xa_y}{a_z}f(x)g(y), \quad \forall z\in \ov{E}_V,
	\end{align*}
	where $V$ is an $\BH$-module with dimension vector $\bd$, and the sum is over all exact sequence
	\begin{align}
		\label{eq:ses}
		0\longrightarrow V(y)\stackrel{\alpha}{\longrightarrow} V(z)\stackrel{\beta}{\longrightarrow} V(x)\longrightarrow0.
	\end{align}
	\begin{lemma}[\text{cf. \cite{Lus90,SV99}}]
		$\cf(\Lambda^\imath)$ is an associative algebra under the product $\star$, and it is isomorphic to $\widetilde{\ch}(\Lambda^\imath)$.
	\end{lemma}
	
	\begin{proof}
		%Keep notations as above. Similarly, we define a new product on $\cf(\Lambda^\imath)$, which is associative (see \cite[\S6]{SV99}): for $f\in\ov{E}_{V'}$, $g\in\ov{E}_{V''}$ with $\dimv V'=\bd'$, $\dimv V''=\bd''$, 
		%\[
		%(f\cdot g)(z)=\frac{\sqq^{\langle \bd',\bd''\rangle_Q}}{|G_{V'}|\cdot |G_{V''}|}\sum_{x,y,\phi,\theta}f(x)g(y), \quad \forall z\in \ov{E}_V,
		%\]
		%where the sum is over all exact sequence \eqref{eq:ses}. These two product coincide up to a linear isomorphism $\varphi:\cf(\Lambda^\imath)\longrightarrow \cf(\Lambda^\imath)$ by $\varphi(f)(x)=a_xf(x)$. %Then we have $\varphi(f\star g)=\varphi(f)\varphi(g)$. For any $f,g,h\in\cf(\Lambda^\imath)$, we have 
		%\begin{align*}
		%  \varphi((f\star g)\star h)&=\varphi(f\star g)\cdot\varphi(h)=\big(\varphi(f)\cdot\varphi(g)\big)\cdot\varphi(h)\\
		%&=\varphi(f)\cdot\big(\varphi(g)\cdot\varphi(h)\big)=\varphi(f\star(g\star h)),
		%\end{align*}
		%   and then $(f\star g)\star h=f\star (g\star h)$. 
		%So $\cf(\Lambda^\imath)$ is an associative algebra under the product $\star$. 
		
		We prove that $\cf(\Lambda^\imath)$ is isomorphic to $\widetilde{\ch}(\Lambda^\imath)$, by sending $\chi_{\mathfrak{O}_x}$ to $[V(x)]$. In fact, for $x\in\ov{E}_{V'}$, $y\in\ov{E}_{V''}$, 
		\begin{align*}
			\chi_{\mathfrak{O}_x}\star 
			\chi_{\mathfrak{O}_y}(z)=&\frac{\sqq^{\langle \bd',\bd''\rangle_Q}}{|G_{V'}|\cdot|G_{V''}|}\frac{a_xa_y}{a_z}P_{V(x),V(y)}^{V(z)}\cdot|\mathfrak{O}_x|\cdot|\mathfrak{O}_y|
			\\
			=&\sqq^{\langle \bd',\bd''\rangle_Q}\frac{P_{V(x),V(y)}^{V(z)}}{a_z}=\sqq^{\langle \bd',\bd''\rangle_Q}G_{V(x),V(y)}^{V(z)}.
		\end{align*}
		This completes the proof. 
	\end{proof}

	%%%%%%%%%
	Next we define $\cg\cf(\Lambda)_V\subseteq \cf(\Lambda^\imath)_V$ to be the set of all $G_V$-invariant functions $f:\ov{E}_V\rightarrow \C$ such that $f(x)=f(y)$ if $V(x)\cong V(y)$ in $\cd_{sg}(\Lambda^\imath)$, and let
	\[
	\cg\cf(\Lambda^\imath)=\bigoplus_{[V]\in\Iso(\mod(\BH))}\cg\cf(\Lambda^\imath)_V,
	\]
	which is a subspace (rather than subalgebra) of $\cf(\Lambda^\imath)$. As in \S\ref{subsec:iHall}, for any $x\in\ov{E}_V$, we set 
	\begin{align}
		\label{eq:Px}
		\mathfrak{S}_x:=\{\mathfrak{O}_y\mid y\in\ov{E}_V, V(y)\cong V(x) \text{ in }\cd_{sg}(\Lambda^\imath)\}.
	\end{align}
	
	Let us endow the subspace  $\cg\cf(\Lambda^\imath)$ an algebra structure. For any $f,g\in\cg\cf(\Lambda^\imath)$, we define
	\begin{align}
		(f\ast g)(z)=\frac{1}{|\mathfrak{S}_z|}\sum_{\mathfrak{O}_{\tilde{z}}\in\mathfrak{S}_z}(f\star g)(\tilde{z}). 
	\end{align}
	It is clear that $f\ast g\in\cg\cf(\Lambda^\imath)$.

	%In fact, for any $x\in\ov{E}_V$, this exist (unique up to isomorphism) $H(V(x))\in\mod(\bfk Q)$ and $r=(r_i)\in\N^{\I}$ such that $V(x)\cong H(V(x))$ in $\cd_{sg}(\Lambda^\imath)$ and $y\in\ov{E}_V$, where $y$ corresponds to $H(V(x))\oplus \bigoplus_{i\in\I}\K_i^{r_i}$.
	
	\begin{proposition}
		\label{prop:quotient}
		$\cg\cf(\Lambda^\imath)$ is an associative algebra under $\ast$. Moreover, it is isomorphic to $\widehat{\ch}(\bfk Q,\varrho)$.
	\end{proposition}
	
	\begin{proof}
		It is enough to prove the second statement. For any $\mathfrak{O}_x$ in $\ov{E}_V$, we define 
		\[
		\ov{\chi}_{\mathfrak{O}_x}:=\sum_{\mathfrak{O}_y\in\mathfrak{S}_x} \chi_{\mathfrak{O}_y}.
		\]
		Then $\ov{\chi}_{\mathfrak{O}_x}\in\cg\cf(\Lambda^\imath)_V$. Let $\mathcal{B}$ be the set of representatives $\mathfrak{O}_x$ of the sets $\mathfrak{S}_x$. Then $\{\ov{\chi}_{\mathfrak{O}_x}\mid\mathfrak{O}_x\in\mathcal{B}\}$ forms a basis of $\cg\cf(\Lambda^\imath)$. We show in the following that the linear map 
		\begin{align*}
			\Psi: \cg\cf(\Lambda^\imath)\longrightarrow \widehat{\ch}(\bfk Q,\varrho),
			\quad \frac{\ov{\chi}_{\mathfrak{O}_x}}{|\mathfrak{S}_x|}\mapsto [V(x)]   
		\end{align*}
		is a bijection, which satisfies $\Psi(\ov{\chi}_{\mathfrak{O}_x}\ast \ov{\chi}_{\mathfrak{O}_y})=[V(x)]\ast[V(y)]$. From the basis of $\widehat{\ch}(\bfk Q,\varrho)$ (see \S\ref{subsec:iHall}), we know that $\Psi$ is a linear  isomorphism. 
		
		%Then we must show that if $w,z\in\ov{E}_V$ such that $H(V(w))\cong H(V(z))$ in $\cd_{sg}(\Lambda^\imath)$, then
		%\begin{align}
		%\ov{\chi}_{\mathfrak{O}_x}*\ov{\chi}_{\mathfrak{O}_y}(w)= \ov{\chi}_{\mathfrak{O}_x}*\ov{\chi}_{\mathfrak{O}_{y}}(z).
		%\end{align}
		
		By definition, we have
		\begin{align*}
			(\ov{\chi}_{\mathfrak{O}_x}\star\ov{\chi}_{\mathfrak{O}_y})(w)=&\frac{\sqq^{\langle\bd',\bd''\rangle_Q}}{|G_{V'}|\cdot|G_{V''}|}\sum_{x',y',\phi,\theta}\frac{a_{x'}a_{y'}}{a_w}\ov{\chi}_{\mathfrak{O}_x}(x') \ov{\chi}_{\mathfrak{O}_y}(y'),
		\end{align*}
		where the sum is over all exact sequences $0\rightarrow V(y')\rightarrow V(w)\rightarrow V(x')\rightarrow 0$. Then
		\begin{align*}
			(\ov{\chi}_{\mathfrak{O}_x}\star\ov{\chi}_{\mathfrak{O}_y})(w)=&\frac{\sqq^{\langle\bd',\bd''\rangle_Q}}{|G_{V'}|\cdot|G_{V''}|}\sum_{\mathfrak{O}_{x'}\in\mathfrak{S}_x,\mathfrak{O}_{y'}\in\mathfrak{S}_y} \frac{a_{x'}a_{y'}}{a_w}F_{x',y'}^wa_{x'}a_{y'}|\mathfrak{O}_{x'}|\cdot|\mathfrak{O}_{y'}|
			\\
			=&\sqq^{\langle\bd',\bd''\rangle_Q}\sum_{\mathfrak{O}_{x'}\in\mathfrak{S}_x,\mathfrak{O}_{y'}\in\mathfrak{S}_y} G_{x',y'}^w,
		\end{align*}
		and so
		\begin{align*}
			\ov{\chi}_{\mathfrak{O}_x}\ast\ov{\chi}_{\mathfrak{O}_y}(z)=&\sqq^{\langle\bd',\bd''\rangle_Q}\frac{1}{|\mathfrak{S}_z|}\sum_{\mathfrak{O}_{w}\in\mathfrak{S}_z} \sum_{\mathfrak{O}_{x'}\in\mathfrak{S}_x,\mathfrak{O}_{y'}\in\mathfrak{S}_y} G_{x',y'}^w
			=\frac{\sqq^{\langle\bd',\bd''\rangle_Q}}{|\mathfrak{S}_z|}\sum_{\mathfrak{O}_{x'}\in\mathfrak{S}_x,\mathfrak{O}_{y'}\in\mathfrak{S}_y} \big(\sum_{\mathfrak{O}_{w}\in\mathfrak{S}_z} G_{x',y'}^w\big).
		\end{align*}
		
		By Lemma \ref{lem:independent}, we know that  all the $\sum_{\mathfrak{O}_{w}\in\mathfrak{S}_z} G_{x',y'}^w$ are equal for any $\mathfrak{O}_{x'}\in\mathfrak{S}_x,\mathfrak{O}_{y'}\in\mathfrak{S}_y$. So
		\begin{align*}
			\ov{\chi}_{\mathfrak{O}_x}\ast\ov{\chi}_{\mathfrak{O}_y}(z)=&\frac{\sqq^{\langle\bd',\bd''\rangle_Q}}{|\mathfrak{S}_z|}|\mathfrak{S}_x|\cdot|\mathfrak{S}_y|\big(\sum_{\mathfrak{O}_{w}\in\mathfrak{S}_z} G_{x,y}^w\big).
		\end{align*}
		%and then
		%\begin{align*}
		%\ov{\chi}_{\mathfrak{O}_x}\ast\ov{\chi}_{\mathfrak{O}_y}=\sum_{\mathfrak{O}_z\in\mathfrak{B}} \frac{\sqq^{\langle\bd',\bd''\rangle_Q}}{|\mathfrak{S}_z|}|\mathfrak{S}_x|\cdot|\mathfrak{S}_y|\big(\sum_{\mathfrak{O}_{w}\in\mathfrak{S}_z} G_{x,y}^w\big)\ov{\chi}_{\mathfrak{O}_z}.
		%\end{align*}
		Together with \eqref{eq:prod}, we get that $\Psi(\ov{\chi}_{\mathfrak{O}_x}\ast \ov{\chi}_{\mathfrak{O}_y})=[V(x)]*[V(y)]$. 
	\end{proof}
	
	\begin{corollary}
		\label{cor:decomposition}
		The algebra $\cg\cf(\Lambda^\imath)$ has a basis given by
		$$\{\ov{\chi}_{\mathfrak{O}_y}\ast\ov{\chi}_{\mathfrak{O}_x}\mid V(x)\in\mod(\bfk Q)\subseteq\mod(\Lambda^\imath), V(y)\in\proj(\BH)\subseteq \mod(\Lambda^\imath)\}.$$ 
	\end{corollary}
	
	\begin{proof}
		%We only need to consider $f=\ov{\chi}_{\mathfrak{O}_x}$ for $x\in\ov{E}_V$. 
		By Proposition \ref{prop:quotient}, it is equivalent to consider $\widehat{\ch}(\bfk Q,\varrho)$, and then the desired result follows from Lemma \ref{basis-iHall}.  
	\end{proof}
	
	For a $\Lambda^\imath$-module $M$, we also write $\chi_{[M]}$ and $\ov{\chi}_{[M]}$ instead of $\chi_{\mathfrak{O}_M}$ and $\ov{\chi}_{\mathfrak{O}_M}$ for simplicity. 
	For convenience, we define
	$$\K_X:=\bigoplus_{i\in\I}\K_i^{\oplus a_i}\in\proj(\BH)\subset\mod(\Lambda^\imath)$$
	for any semisimple $\BH$-module $X=\oplus_{i\in\I}S_i^{\oplus a_i}$.
	
	\begin{corollary}
		\label{cor:KX=XK}
		Let $W\in\mod(\BH)$ be semisimple, $K=\K_X\in\proj(\BH)\subseteq\mod(\Lambda^\imath)$. For any  $g\in\cg\cf(\Lambda^\imath)_{K}$, $h\in\cg\cf(\Lambda^\imath)_{W}$, we have
		$$g\ast h=\sqq^{(\btau X, W)_Q-(X,W)_Q}h\ast g.$$
	\end{corollary}
	
	\begin{proof}
		We can assume $g=\ov{\chi}_{\mathfrak{O}_y}$, $h=\ov{\chi}_{\mathfrak{O}_x}$. It is equivalent to consider $\widehat{\ch}(\bfk Q,\varrho)$, and then the result follows from \eqref{eq:KX=XK}.
	\end{proof}
	
	By Proposition \ref{prop:iHallmult}, we have the following corollary immediately.
	\begin{corollary}
		\label{cor:hallmult}
		For any $x\in\ov{E}_{V'},y\in\ov{E}_{V''}$ with $V',V''\in\mod(\BH)$ semisimple, we have
		\begin{align*}
			&\ov{\chi}_{\mathfrak{O}_x}\ast\ov{\chi}_{\mathfrak{O}_y}
			%  \\
			%    &=\sum_{[L],[M],[N],[X]} \sqq^{\langle X,M\rangle_Q-\langle \varrho X,M\rangle_Q-\langle V(x),V(y)\rangle_Q} q^{\langle N,L\rangle_Q} F_{N,L}^M F_{X,N}^{V(x)}F_{ L,\varrho (X)}^{V(y)} 
			%		 \frac{a_L a_N a_X}{a_M} \cdot\frac{\ov{\chi}_{[\K_{X}]}}{|\mathfrak{S}_{\K_X}|}\ast\ov{\chi}_{[M]}
			=\sum_{[L],[M],[N],[X]} \sqq^{\langle X,M\rangle_Q-\langle \varrho X,M\rangle_Q-\langle A,B\rangle_Q} q^{\langle N,L\rangle_Q} \frac{P_{N,L}^M P_{X,N}^{V(x)}P_{ L,\varrho (X)}^{V(y)} 
			}{a_L a_N a_Xa_M} \cdot\frac{\ov{\chi}_{[\K_{X}]}}{|\mathfrak{S}_{\K_X}|}\ast\ov{\chi}_{[M]}.
		\end{align*}
	\end{corollary}
	
	%%%%%%%%%%%%%%%
	
	\subsection{Fourier transforms of $\imath$Hall algebras}\label{FT of iHall subsec}
	Let $\ce$ be a subset of $Q_1$ closed under the action of $\btau$ and denote by $Q'$ the quiver obtained from $Q$ by reversing all the arrows in $\ce$. Then $(Q',\varrho)$ is an $\imath$quiver, and we denote ${'\Lambda^{\imath}}$ its $\imath$quiver algebra. Let $V=\oplus_{i\in\I}V_i$ be an $\I$-graded $\F_q$-vector space of dimension vector $\bd$, which is a semi-simple $\BH$-module. Note that $E_{Q,V}=\ov{E}_{Q,V}=E_{Q,\bd}$. In this case, we also denote $G_V=\prod_{i\in\I}{\rm GL}(V_i)$ by $G_\bd$. Set
	
	\begin{align*}
		X_{\bd}=X_V=&\prod_{\alpha\in Q_1\setminus \ce}\Hom_\bfk(V_{\ts(\alpha)},V_{\tt(\alpha)}),
		\qquad
		Y_{\bd}=Y_V=\prod_{\alpha\in \ce}\Hom_\bfk(V_{\ts(\alpha)},V_{\tt(\alpha)}),\\
		&Z_{\bd}=Z_V=\prod_{\alpha\in \ce}\Hom_\bfk(V_{\tt(\alpha)},V_{\ts(\alpha)}).
	\end{align*}
	Then there is a decomposition
	\begin{align}
		E_{Q,V}=X_V\times Y_V,\qquad E_{Q',V}=X_V\times Z_V.
	\end{align}
	
	Define a non-degenerate bilinear form
	$$\{-,-\}:Y_V\times Z_V\longrightarrow \bfk,\quad\big((x_\alpha),  (y_\alpha)\big)\mapsto \sum_{\alpha\in E}\mathrm{tr}(x_\alpha y_\alpha),$$
	which extends to a blinear form 
	\[\{-,-\}:E_{Q,V}\times E_{Q',V}\longrightarrow \bfk\]
	by noting that $E_{Q,V}=X_V\times Y_V$, and $E_{Q',V}=X_V\times Z_V$.  
	
	Fix a non-trivial additive character $\phi:\bfk\rightarrow \C^\times$. Following \cite{SV99}, for any $f\in \mathcal{F}(\Lambda^\imath)_V$, that is, a $G_V$-invariant function on $E_{Q,V}=X_V\times Y_V$, we define $\hat{f}:E_{Q',V}=X_V\times Z_V\rightarrow\C$ in $\cf('\Lambda^\imath)_V$ by setting
	\begin{align}
		\label{eq:FT}
		\hat{f}(x,w)=\sqq^{-\dim Y_V}\sum_{y\in Y_V}\frac{a_{(x,v)}}{a_{(x,w)}}f(x,y)\phi(\{ y,w\}),\quad \forall (x,w)\in X_V\times Z_V.
	\end{align}
	
	\begin{lemma}
		\label{lem:FT-char}
		For any $(x,v)\in X_{\bd}\times Y_\bd$, we have
		\begin{align}
			\widehat{\chi_{\mathfrak{O}_{(x,v)}}}=\sqq^{-\dim Y_V} \sum_{w\in Z_\bd}
			\phi(\{ v,w \}) \chi_{\mathfrak{O}_{(x,w)}}.
		\end{align}
	\end{lemma}
	
	\begin{proof}
		By definition, we have 
		\begin{align*}
			\widehat{\chi_{(x,v)}}=&\sum_{\tilde{x}\in X_\bd,w\in Z_\bd}\sqq^{-\dim Y_V} \sum_{y\in Y_V}\frac{a_{(\tilde{x},v)}}{a_{(\tilde{x},w)}} \phi(\{ y,w \})  \chi_{(x,v)}(\tilde{x},y)\chi_{(\tilde{x},w)}
			\\
			=&\sum_{w\in Z_\bd}\sqq^{-\dim Y_V} \frac{a_{(x,v)}}{a_{(x,w)}} \phi(\{ v,w \}) \chi_{(x,w)}.
		\end{align*}
		Since Fourier transform can restrict to $G_V$-invariant subset (see e.g. \cite[\S7]{SV99}), we have
		\begin{align*}
			\widehat{\chi_{\mathfrak{O}_{(x,v)}}}=&\sqq^{-\dim Y_V} \sum_{w\in Z_\bd}  \sum_{\tilde{x}\in X_V,\tilde{v}\in Y_V:V(\tilde{x},\tilde{v})\cong V(x,v)}\frac{a_{(\tilde{x},\tilde{v})}}{a_{(\tilde{x},w)}} \phi(\{ \tilde{v},w \}) \chi_{{(\tilde{x},w)}}
			\\
			=&\sqq^{-\dim Y_V} \sum_{w\in Z_\bd}  \sum_{\tilde{x}\in X_V,\tilde{v}\in Y_V:V(\tilde{x},\tilde{v})\cong V(x,v)}\frac{a_{(\tilde{x},\tilde{v})}}{a_{(\tilde{x},w)}|\mathfrak{O}_{(\tilde{x},w)}|} \phi(\{ \tilde{v},w \}) \chi_{\mathfrak{O}_{(\tilde{x},w)}}
			\\
			=&\sqq^{-\dim Y_V} \sum_{w\in Z_\bd}  \sum_{\tilde{x}\in X_V,\tilde{v}\in Y_V:V(\tilde{x},\tilde{v})\cong V(x,v)}\frac{a_{(\tilde{x},\tilde{v})}}{|G_V|} \phi(\{ \tilde{v},w \}) \chi_{\mathfrak{O}_{(\tilde{x},w)}}
			\\
			=&\sqq^{-\dim Y_V} \sum_{w\in Z_\bd}  \frac{a_{(x,v)}|\mathfrak{O}_{(x,v)}|}{|G_V|} \phi(\{ v,w \}) \chi_{\mathfrak{O}_{(x,w)}}
			\\
			=&\sqq^{-\dim Y_V} \sum_{w\in Z_\bd}   \phi(\{ v,w \}) \chi_{\mathfrak{O}_{(x,w)}}.\qedhere
		\end{align*}
	\end{proof}

	For any $\ov{\chi}_{\mathfrak{O}_y},\ov{\chi}_{\mathfrak{O}_x}$ with $V(x)\in\mod(\bfk Q)\subseteq\mod(\Lambda^\imath)$ and $V(y)\in\proj(\BH)\subseteq \mod(\Lambda^\imath)$, we have $\dim\cg\cf(\Lambda^\imath)_{V(y)}=1$, and $\ov{\chi}_{\mathfrak{O}_{y}}$ is its basis. Since $\mod(\BH)\subseteq \mod('\Lambda^\imath)$, we can view $V(y)\in \mod('\Lambda^\imath)$, and then $\ov{\chi}_{\mathfrak{O}_y}\in \cg\cf('\Lambda^\imath)$. 
	By Corollary \ref{cor:decomposition}, we define the $\C$-linear map 
	\begin{align}
		\Phi_{Q',Q}:\cg\cf(\Lambda^\imath)\longrightarrow\cg\cf('\Lambda^\imath), \quad \frac{\ov{\chi}_{\mathfrak{O}_y}}{|\mathfrak{S}_y|}\ast \ov{\chi}_{\mathfrak{O}_x}\mapsto \frac{\ov{\chi}_{\mathfrak{O}_y}}{|\mathfrak{S}_y'|}\ast \widehat{\ov{\chi}_{\mathfrak{O}_x}},
	\end{align}
	called the Fourier transform from $\cg\cf(\Lambda^\imath)$ to $\cg\cf('\Lambda^\imath)$. Here we use $\mathfrak{S}'_y$ to denote the set in \eqref{eq:Px} for $V(y)$ viewed as $'\Lambda^\imath$-module.

	Similarly, using the conjugate character $\ov{\phi}$, we obtain a $\C$-linear isomorphism $\ov{\Phi}_{Q',Q}:\cg\cf(\Lambda^\imath)\rightarrow \cg\cf('\Lambda^\imath)$.
	
	%We identify $\cf(\Lambda^\imath)$ and $\cg\cf('\Lambda^\imath)$ with $\widetilde{\ch}(\Lambda^\imath)$ and $\widetilde{\ch}('\Lambda^\imath)$ respectively in the following. Then we can view $\Phi_{Q',Q}$ and $\ov{\Phi}_{Q',Q}$ 
	%as maps $\widetilde{\ch}(\Lambda^\imath)\rightarrow \widetilde{\ch}('\Lambda^\imath)$. In general, $\Phi_{Q',Q}$ and $\ov{\Phi}_{Q',Q}$ 
	%are not homomorphisms of algebras:  $\widetilde{\ch}(\Lambda^\imath)\rightarrow \widetilde{\ch}('\Lambda^\imath)$. \red{By definition, 
	%the linear map   $\Phi_{Q',Q}: \cf(\Lambda^\imath)\rightarrow \cg\cf('\Lambda^\imath)$ induces a linear map $\cg\cf(\Lambda^\imath)\rightarrow \cg\cf('\Lambda^\imath)$, which is also denoted by $\Phi_{Q',Q}$.} Similarly for $\ov{\Phi}_{Q',Q}$. 
	Using Proposition \ref{prop:quotient}, we identify $\cg\cf(\Lambda^\imath)$ and $\cg\cf('\Lambda^\imath)$ with $\widehat{\ch}(\bfk Q,\varrho)$ and $\widehat{\ch}(\bfk Q',\varrho)$ respectively. The following lemma is obvious. %, where $\mathcal{I}'$ is the subspace defined in \eqref{def:I} for $\widetilde{\ch}('\Lambda^\imath)$.

	%It is obvious that the above construction holds for all quivers. Now, let us consider Dynkin quivers. 
	
	\begin{lemma}
		\label{lem:iso}
		Let $Q$ be a Dynkin quiver.
		For any $X\in E_V=E_{Q,V}$, define $\widetilde{X}\in E_V$ such that $\widetilde{X}_h=-X_h$ for $h\in \ce$ and $\widetilde{X}_h=X_h$ for $h\notin \ce$. %If $Q$ is a tree quiver, 
		Then we have $\widetilde{X}\cong X$.  
	\end{lemma}
	
	%\begin{proof}
	%   We assume that $\ce=\{h\}$. 
	%  Since $Q$ is a Dynkin quiver, we set $\ct(h)=\{i\in\I\mid  \text{there is a path from }\tt(h) \text{ to }i\}$, and $\cs(h)=Q_0\setminus\ct(h)$. Then there is an isomorphism $g=(g_i)_{i\in\I}:X\rightarrow \widetilde{X}$ where $g_i=\Id$ if $i\in\cs(h)$, and $g_i=-\Id$ if $i\in \ct(h)$.  
	%\end{proof}
	
	%It is remarkable that the above lemma holds for tree quivers.
	
	\begin{theorem}
		\label{thm:FThall}
		Let $(Q,\btau)$ be a Dynkin $\imath$quiver, and $(Q',\btau)$ be the Dynkin $\imath$quiver constructed from $(Q,\btau)$ by reversing the arrows in a subset $\ce$ of $Q_1$ satisfying $\btau(\ce)=\ce$. Then the Fourier transform $\Phi_{Q',Q}:\widehat{\ch}(\bfk Q,\btau)\rightarrow \widehat{\ch}(\bfk Q',\btau)$ is an isomorphism of $\C$-algebras with the inverse given by $\ov{\Phi}_{Q,Q'}:\widehat{\ch}(\bfk Q',\btau)\rightarrow \widehat{\ch}(\bfk Q,\btau)$.
	\end{theorem}

	\begin{proof}
		We denote $\Phi=\Phi_{Q',Q}$ and $\ov{\Phi}=\ov{\Phi}_{Q',Q}$ for simplicity. 
		By definition, we know $\Phi$ is an $\C$-linear isomorphism with $\ov{\Phi}$ as its inverse. 
		It remains to prove that $\Phi$ is an algebra homomorphism. 
		
		First, we show that $\Phi(f\ast g)=\Phi(f)\ast \Phi(g)$ for any $f=\ov{\chi}_{\mathfrak{O}_{(x_1,v_1)}}\in\cg\cf(\Lambda^\imath)_{V'}$, $g=\ov{\chi}_{\mathfrak{O}_{(x_2,v_2)}}\in \cg\cf(\Lambda^\imath)_{V''}$ such that $V',V''\in\mod(\BH)$ are semi-simple, and $x_1\in X_{V'}$, $v_1\in Y_{V'}$, $x_2\in X_{V''}$, $v_2\in Y_{V''}$.
		%Assume $V=V'\oplus V''$ as $\BH$-module.
		In this case, we have $\ov{\chi}_{\mathfrak{O}_{(x_1,v_1)}}= \chi_{\mathfrak{O}_{(x_1,v_1)}}$, $\ov{\chi}_{\mathfrak{O}_{(x_2,v_2)}}= \chi_{\mathfrak{O}_{(x_2,v_2)}}$. 
		We denote by $\bd',\bd''$ the dimension vectors of $V',V''$ respectively. Let $\bd=\bd'+\bd''$. Then we have $f\ast g\in\bigoplus_{[V]}\cg\cf(\Lambda^\imath)_{V}$ and $\Phi(f)\ast\Phi(g)\in\bigoplus_{ [V]}\cg\cf('\Lambda^\imath)_{V}$,
		where the sum is over $[V]\in\Iso(\mod(\BH))$ such that there exists a short exact sequence 
		$$0\longrightarrow V''\longrightarrow V\longrightarrow V'\longrightarrow0.$$
		By transitivity (\cite[Lemma 5.2]{SV99}) it is enough to consider the case that $\ce$ consists the orbit of one edge: $\ce=\{e,\btau (e)\}$. We assume $\btau (e)\neq e$ and $e:i\rightarrow j$ in the following, since the case $\varrho (e)=e$ is similar.
		
		By definition we have
		\begin{equation}
			\label{eq:RHS-def}
			\begin{aligned}
				\Phi(f)\ast\Phi(g)=&\sqq^{-\bd'_{i}\bd'_{j}-\bd'_{\btau (i)}\bd'_{\btau (j)}-\bd''_{i}\bd''_{j}-\bd''_{\btau (i)}\bd''_{\btau (j)}}\sum_{w_1\in Z_{\bd'},w_2\in Z_{\bd''}} 
				\\
				& \phi(\{ v_1,w_1\}+\{ v_2,w_2\}) \chi_{\mathfrak{O}_{(x_1,w_1)}}\ast \chi_{\mathfrak{O}_{(x_2,w_2)}}.
			\end{aligned}
		\end{equation}
		and then by Corollary \ref{cor:hallmult},  
		\begin{align*}
			\Phi(f)\ast\Phi(g)=&%\sum_{\substack{w_1,w_2,w,\delta\\ M',N',X'}}
			\sum_{[V(x,w)],\delta}\sum_{\substack{w_1\in Z_{\bd'},w_2\in Z_{\bd''}\\ [M'],[N'],[X']:\;\dimv{ X'}=\delta}} \sqq^{\widetilde{T}(\bd',\bd'',\delta)}\phi(\{ v_1,w_1\})\phi(\{ v_2,w_2\}) \frac{1}{a_{M'}a_{N'}a_{X'}}\\
			&P_{M',N'}^{V(x,w)}P_{X',M'}^{V(x_1,w_1)}P_{N',\btau (X')}^{V(x_2,w_2)}\cdot \frac{\ov{\chi}_{[\K_\delta]}}{|\mathfrak{P}'_{\K_\delta}|}\ast \frac{\chi_{\mathfrak{O}_{(x,w)}}}{a_{V(x,w)}}
		\end{align*}
		where $\widetilde{T}(\bd',\bd'',\delta)=-\bd'_{i}\bd'_{j}-\bd'_{\btau (i)}\bd'_{\btau (j)}-\bd''_{i}\bd''_{j}-\bd''_{\btau (i)}\bd''_{\btau (j)}+\langle \delta-\btau (\delta),\bd-\delta-\btau(\delta)\rangle_{Q'}-\langle \bd',\bd''\rangle_{Q'}+2\langle \bd'-\delta,\bd''-\btau(\delta)\rangle_{Q'}$. 
		On the other hand, 
		\begin{align*}
			\Phi(f\ast g)=&\sum_{[V(x,w)],\delta}\sum_{\substack{v\in Y_{\bd-\delta-\btau(\delta)}\\ [M],[N],[X]:\;\dimv{ X}=\delta}}%\sum_{\substack{v,w,\delta\\ M,N,X}}
			\sqq^{T(\bd',\bd'',\delta)}
			\phi(\{ v,w\})\frac{1}{a_{M}a_{N}a_{X}}
			\\
			&P_{M,N}^{V(x,v)}P_{X,M}^{V(x_1,v_1)}P_{N,\btau (X)}^{V(x_2,v_2)}\cdot \frac{\ov{\chi}_{[\K_\delta]}}{|\mathfrak{P}'_{\K_\delta}|}\ast \frac{\chi_{\mathfrak{O}_{(x,w)}}}{G_{V(x,w)}},
		\end{align*}
		where 
		$T(\bd',\bd'',\delta)=-(\bd-\delta-\btau(\delta))_{i}(\bd-\delta-\btau(\delta))_{j}-(\bd-\delta-\btau(\delta))_{\btau (i)}(\bd-\delta-\btau(\delta))_{\btau (j)}+\langle \delta-\btau(\delta),\bd-\delta-\btau(\delta)\rangle_Q-\langle \bd',\bd''\rangle_Q+2\langle \bd'-\delta,\bd''-\btau(\delta)\rangle_Q$.
		
		Note that 
		\begin{align}
			\label{eq:TT}
			\widetilde{T}(\bd',\bd'',\delta)=&T(\bd',\bd'',\delta) +2\bd_i'\bd_j''+2\bd_{\btau (i)}'\bd_{\btau (j)}''-2(\bd_i'+\bd_{\btau (i)}')(\delta_j+\delta_{\btau (j)})
			\\\notag
			&-2(\bd_j''+\bd_{\btau (j)}'')(\delta_i+\delta_{\btau (i)}) +2(\delta_i+\delta_{\btau (i)})(\delta_j+\delta_{\btau (j)}). 
		\end{align}
		For fixed $[V(x,w)]$ and $\delta$, we have (by identifying $K_0(\bfk Q)=\Z^\I=K_0(\bfk Q')$) 
		\[\dimv X=\dimv X'=\delta,\quad \dimv M=\dimv M'=\bd'-\delta,\quad \dimv N=\dimv N'=\bd''-\btau (\delta).\]
		Note that $|\Aut(A)|\cdot|\mathfrak{O}_{A}|=|G_{\dimv A}|$ for any $A\in\mod(\bfk Q)$ or $\mod(\bfk Q')$. 
		It then suffices to prove the incidence of the following two expressions for fixed $[V(x,w)]$ and $\delta$: 
		\begin{align}
			\label{eq:LHS}
			&\sum_{\substack{w_1\in Z_{\bd'},w_2\in Z_{\bd''}\\ M',N',X':\;\dimv X'=\delta}} \sqq^{\widetilde{T}(\bd',\bd'',\delta)}\phi(\{ v_1,w_1\})\phi(\{ v_2,w_2\}) P_{M',N'}^{V(x,w)}P_{X',M'}^{V(x_1,w_1)}P_{N',\btau (X')}^{V(x_2,w_2)},
			\\
			\label{eq:RHS}
			&\sum_{\substack{v\in Y_{\bd-\delta-\varrho\delta}\\ M,N,X:\;\dimv X=\delta}}\sqq^{T(\bd',\bd'',\delta)}\phi(\{ v,w\}) P_{M,N}^{V(x,v)}P_{X,M}^{V(x_1,v_1)}P_{N,\btau (X)}^{V(x_2,v_2)},
		\end{align}
		where the summations are over $M',N',X'\in\mod(\bfk Q')$ and $M,N,X\in\mod(\bfk Q)$.

		%\begin{align*}
		%LHS&=\sum_{\substack{M',N',X',w_1,w_2}}\sqq^{\widetilde{T}(\bd',\bd'',\delta)}\psi(\langle v_1,w_1\rangle)\psi(\langle v_2,w_2\rangle) P_{M',N'}^{V(x,w)}P_{X',M'}^{V(x_1,w_1)}P_{N',\tau X'}^{V(x_2,w_2)}\\
		%RHS&=\sum_{\substack{M,N,X,v}}\sqq^{T(\bd',\bd'',\delta)}\psi(\langle v,w\rangle) P_{M,N}^{V(x,v)}P_{X,M}^{V(x_1,v_1)}P_{N,\tau X}^{V(x_2,v_2)}
		%\end{align*}
		%where the summations are over %$\widehat{X}=\widehat{X}'=\delta$.
		
		Before going on, let us recall 
		the twisted Hall algebra $\widetilde{\ch}(\bfk Q')$ and its Fourier transforms given in \cite{Rin90,SV99}. 
		Recall that for any $A',B'\in\mod(\bfk Q')$, we have 
		\[[A']\cdot [B']=\sum_{[L']}\sqq^{\langle A',B'\rangle_{Q'}} G_{A',B'}^{L'}[L'].\]
		Following \cite{SV99}, let $\cf(Q')_\bd$ denote the set of all functions $E_{Q',V}\rightarrow \C$ for fixed $\I$-graded vector space $V=\oplus_{i\in\I} V_i$ of dimension vector $\bd$. Let $\cf(Q')=\oplus_{\bd}\cf(Q')_\bd$. Then $\cf(Q')$ is an associative algebra with its product given by
		\begin{align*}
			(f\cdot g)(z)=\frac{\sqq^{\langle \bd',\bd''\rangle_Q}}{|G_{V'}|\cdot |G_{V''}|}\sum_{x,y,\alpha,\beta} \frac{a_xa_y}{a_z}f(x)g(y), \quad \forall z\in E_{Q,\bd'+\bd''}
		\end{align*}
		for any $f\in\cf(Q')_{\bd'}$, $g\in\cf(Q')_{\bd''}$, 
		where the sum is over all exact sequence \eqref{eq:ses} in $\mod(\bfk Q')$.
		%\begin{align}
		%\label{eq:ses}
		%0\longrightarrow V(y)\stackrel{\alpha}{\longrightarrow} V(z)\stackrel{\beta}{\longrightarrow} V(x)\longrightarrow0.
		%\end{align}
		Then $\widetilde{\ch}(\bfk Q')$ is isomorphic to $(\cf(Q'),\cdot)$ by mapping $[A']$ to $\chi_{\mathfrak{O}_{A'}}$.
		Moreover, we have 
		\[\chi_{\mathfrak{O}_{A'}}\cdot\chi_{\mathfrak{O}_{B'}}=\sqq^{\langle A',B'\rangle_{Q'}}\sum_{[L]}G_{A',B'}^{L'}\;\chi_{\mathfrak{O}_{L'}}=\sqq^{\langle A',B'\rangle_{Q'}}\sum_{[L']}P_{A',B'}^{L'} \frac{\chi_{\mathfrak{O}_{L'}}}{a_{L'}}.\]
		%in $\widetilde{\ch}(\bfk Q')$. 

		We also have the convolution algebra $\cf(Q,\cdot)$ for the quiver $Q$. 
		Using the formula \eqref{eq:FT}, we can define the Fourier transform from $\cf(Q',\cdot)$ to $\cf(Q,\cdot)$; see also \cite{SV99}.
		
		Return to the proof. By viewing $\chi_{\mathfrak{O}_{M'}},\chi_{\mathfrak{O}_{N'}},\chi_{\mathfrak{O}_{X'}}$ in $\cf(Q')$,  using the definition of Fourier transforms, we have 
		\begin{align*}
			\widehat{\chi_{\mathfrak{O}_{X'}}\cdot\chi_{\mathfrak{O}_{M'}}}(x_1,v_1)&=\sqq^{-\dim(Y_{(x_1,v_1)})}\sum_{w_1\in Z_{\bd'}}\phi(\{ v_1,w_1\})\frac{a_{(x_1,w_1)}}{a_{(x_1,v_1)}}\cdot \sqq^{\{ X',M'\}_{Q'}}P_{X',M'}^{V(x_1,w_1)}\frac{1}{a_{(x_1,w_1)}}\\
			&=\frac{\sqq^{\langle X',M'\rangle_{Q'}-\dim(Y_{(x_1,v_1)})}}{a_{(x_1,v_1)}}\sum_{w_1\in Z_{\bd'}}\phi(\{ v_1,w_1\})P_{X',M'}^{V(x_1,w_1)}.
		\end{align*}
		On the other hand,
		\begin{align*}
			\widehat{\chi_{\mathfrak{O}_{X'}}}\cdot\widehat{\chi_{\mathfrak{O}_{M'}}}(x_1,v_1) 
			&= \sqq^{-\dim(Y_{X'})-\dim(Y_{M'})}\sum_{X_1,M}\phi(\{ v_{X_1},w_{X'}\}+\{v_M,w_{M'}\})(\chi_{\mathfrak{O}_X}\cdot\chi_{\mathfrak{O}_M})(x_1,v_1)\\
			&=\sqq^{\langle X_1,M\rangle_{Q}-\dim(Y_{X'})-\dim(Y_{M'})}\sum_{X_1,M}\frac{\phi(\{ v_{X_1},w_{X'}\})\phi(\{ v_M,w_{M'}\})}{a_{(x_1,v_1)}}P^{V(x_1,v_1)}_{X_1,M},
		\end{align*}
		where the summation is over $X_1=(x_{X_1},v_{X_1})\in E_{Q,\delta}=X_{\delta}\times Y_{\delta}$ and $M=(x_M,v_M)\in E_{Q,\bd'-\delta}=X_{\bd'-\delta}\times Y_{\bd'-\delta}$, and we write $X'=(x_{X'},v_{X'})\in E_{Q',\delta}=X_\delta\times Z_\delta$, $M'=(x_{M'},w_{M'})\in E_{Q',\delta}=X_{\bd'-\delta}\times Z_{\bd'-\delta}$. By \cite[Theorem 7.1]{SV99}, we know that the Fourier transform is an algebra homomorphism from $\cf(Q',\cdot)$ to $\cf(Q,\cdot)$. For $X',M'$, we deduce that
		\begin{align*}
			&\sum_{w_1\in Z_{\bd'}}\phi(\{ v_1,w_1\})P_{X',M'}^{V(x_1,w_1)}=\sqq^{-\langle X',M'\rangle_{Q'}+\dim(Y_{(x_1,v_1)})}a_{(x_1,v_1)}\widehat{\chi_{\mathfrak{O}_{X'}}\cdot\chi_{\mathfrak{O}_{M'}}}(x_1,v_1)\\
			&=\sqq^{-\langle X',M'\rangle_{Q'}+\dim(Y_{(x_1,v_1)})}a_{(x_1,v_1)}\widehat{\chi_{\mathfrak{O}_{X'}}}\cdot\widehat{\chi_{\mathfrak{O}_{M'}}}(x_1,v_1)\\
			&=\sqq^{\langle X_1,M\rangle_{Q}-\langle X',M'\rangle_{Q'}+\dim(Y_{(x_1,v_1)})-\dim(Y_{X'})-\dim(Y_{M'})}\sum_{X_1,M}\phi(\{ v_{X_1},w_{X'}\}+\{ v_M,w_{M'}\})P^{V(x_1,v_1)}_{X_1,M}\\
			&=\sqq^{2(\bd'_i\delta_j+\bd'_{\btau (i)}\delta_{\btau (j)})-2(\delta_i\delta_j+\delta_{\btau( i)}\delta_{\btau (j)})}\sum_{X_1,M}\phi(\{ v_{X_1},w_{X'}\}+\{ v_M,w_{M'}\})P^{V(x_1,v_1)}_{X_1,M},
		\end{align*}
		by noting that $\dimv X'=\delta$, and $\dimv M'=\bd'-\delta$. 
		
		Similarly, for $N',X'$, we find that 
		\begin{align*}
			\sum_{w_2\in Z_{\bd''}} \phi(\{ v_2,w_2\}) P^{V(x_2,w_2)}_{N',\btau (X')} = &\sqq^{2(\bd''_j\delta_{\btau (i)}+\bd''_{\btau (j)}\delta_i)-2(\delta_i\delta_j+\delta_{\btau (i)}\delta_{\btau (j)})} \\
			&\sum_{N,X_2}\phi(\{ v_N,w_{N'}\}+\{ v_{X_2},w_{X'}\}) P^{V(x_2,v_2)}_{N,\btau (X_2)},
		\end{align*}
		where we set $N=(x_N,v_N)\in E_{Q,\bd''-\varrho(\delta)}=X_{\bd''-\varrho(\delta)}\times Y_{\bd''-\varrho(\delta)}$ and $N'=(x_{N'},w_{N'})\in E_{Q',\bd''-\varrho(\delta)}=X_{\bd''-\varrho(\delta)}\times Z_{\bd''-\varrho(\delta)}$. 
		Plugging these two formulae into \eqref{eq:LHS}, we have \eqref{eq:LHS} equals to 
		\begin{align*}
			%\eqref{eq:LHS}&=
			&\sum_{\substack{M',N',X'\\ M,N,X_1,X_2}} \sqq^{\widetilde{T}_1(\bd',\bd'',\delta)} \phi(\{ v_M,w_{M'}\}+\{ v_N,w_{N'}\})\phi(\{ v_{X_1}+v_{X_2},w_{X'}\})  P_{M',N'}^{V(x,w)}P_{X_1,M}^{V(x_1,v_1)}P_{N,\btau (X_2)}^{V(x_2,v_2)}\\
			&=\sum_{\substack{M',N',M,N,X}}\sqq^{\widetilde{T}_1(\bd',\bd'',\delta)+2(\delta_i\delta_j+\delta_{\btau (i)}\delta_{\btau (j)})} \phi(\{ v_M,w_{M'}\}+\{ v_N,w_{N'}\}) P_{M',N'}^{V(x,w)}P_{X,M}^{V(x_1,v_1)}P_{N,\btau (\widetilde{X})}^{V(x_2,v_2)},
		\end{align*}
		where
		\begin{align*}
			\widetilde{T}_1(\bd',\bd'',\delta):=&\widetilde{T}(\bd',\bd'',\delta)+2(\bd'_i\delta_j+\bd'_{\btau (i)}\delta_{\btau (j)})-2(\delta_i\delta_j+\delta_{\btau (i)}\delta_{\btau (j)})
			\\
			&+2(\bd''_j\delta_{\btau (i)}+\bd''_{\btau (j)}\delta_i)-2(\delta_i\delta_j+\delta_{\btau (i)}\delta_{\btau (j)})
		\end{align*}
		and $\widetilde{X}=(x_{\widetilde{X}},v_{\widetilde{X}})$ satisfies $x_{\widetilde{X}}=x_X$, $v_{\widetilde{X}}=-v_X$ as in Lemma \ref{lem:iso}. 
		
		Similarly, %the formula \eqref{eq:RHS} can be equal to 
		\begin{align*}
			\eqref{eq:RHS}&=\sum_{M,N,X,M',N'}\sqq^{T_1(\bd',\bd'',\delta)}\phi(\{ v_M,w_{M'}\}+\{ v_N,w_{N'}\}) P_{M',N'}^{V(x,w)}P_{X,M}^{V(x_1,v_1)}P_{N,\btau (X)}^{V(x_2,v_2)}
		\end{align*}
		where
		\begin{align*}
			T_1(\bd',\bd'',\delta)
			:=&T(\bd',\bd'',\delta)+2\big(\bd'_i\bd''_j+\bd'_{\btau (i)}\bd''_{\btau (j)}-\bd'_i\delta_{\btau (j)}-\bd'_{\btau (i)}\delta_j\big)
			\\
			&-2\big(\bd''_j\delta_i+\bd''_{\btau (j)}\delta_{\btau (i)}-\delta_i\delta_{\btau (j)}-\delta_{\btau (i)}\delta_j\big).
		\end{align*}
		Then \eqref{eq:LHS} equals to \eqref{eq:RHS} provided the following two statements hold. First one should have
		\[\widetilde{T}_1(\bd',\bd'',\delta)+2(\delta_i\delta_j+\delta_{\btau (i)}\delta_{\btau (j)})=T_1(\bd',\bd'',\delta)\]
		which follows from \eqref{eq:TT}.
		Second, one should have $P^{V(x_2,v_2)}_{N,\btau(\widetilde{X})}=P^{V(x_2,v_2)}_{N,\btau (X)}$ if $\widetilde{X}\cong X$, which follows from Lemma \ref{lem:iso}.
		
		Finally, for any $\ov{\chi}_{\mathfrak{O}_{y_1}}$, $\ov{\chi}_{\mathfrak{O}_{y_2}}$ such that $V(y_1)=\K_{\beta_1}, V(y_2)=\K_{\beta_2}\in\proj(\BH)$ for $\beta_1,\beta_2\in\N^\I$, by Corollary \ref{cor:KX=XK}, we have 
		\begin{align*}
			&\Phi\big((\frac{\ov{\chi}_{\mathfrak{O}_{y_1}}}{|\mathfrak{S}_{y_1}|}\ast\ov{\chi}_{\mathfrak{O}_{(x_1,v_1)}})\ast (\frac{\ov{\chi}_{\mathfrak{O}_{y_2}}}{|\mathfrak{S}_{y_2}|}\ast\ov{\chi}_{\mathfrak{O}_{(x_2,v_2)}})\big)
			\\
			&= \sqq^{(\beta_1,\bd')_Q-(\btau (\beta_1),\bd')_Q}\Phi\big(\frac{\ov{\chi}_{\mathfrak{O}_{y_1}}}{|\mathfrak{S}_{y_1}|}\ast \frac{\ov{\chi}_{\mathfrak{O}_{y_2}}}{|\mathfrak{S}_{y_2}|}\ast \ov{\chi}_{\mathfrak{O}_{(x_1,v_1)}}\ast \ov{\chi}_{\mathfrak{O}_{(x_2,v_2)}}\big)
			\\
			&= \sqq^{(\beta_1,\bd')_Q-(\btau (\beta_1),\bd')_Q}\frac{\ov{\chi}_{\mathfrak{O}_{y_1}}}{|\mathfrak{S}'_{y_1}|}\ast \frac{\ov{\chi}_{\mathfrak{O}_{y_2}}}{|\mathfrak{S}'_{y_2}|}\ast \Phi\big(\ov{\chi}_{\mathfrak{O}_{(x_1,v_1)}}\ast \ov{\chi}_{\mathfrak{O}_{(x_2,v_2)}}\big)
			%\\
			%&=\sqq^{(\beta_1,\bd')_{Q'}-(\btau (\beta_1),\bd')_{Q'}}\frac{\ov{\chi}_{\mathfrak{O}_{y_1}}}{|\mathfrak{S}'_{y_1}|}\ast \frac{\ov{\chi}_{\mathfrak{O}_{y_2}}}{|\mathfrak{S}'_{y_2}|}\ast \Phi(\ov{\chi}_{\mathfrak{O}_{(x_1,v_1)}})\ast \Phi(\ov{\chi}_{\mathfrak{O}_{(x_2,v_2)}})
			\\
			&=\Phi(\frac{\ov{\chi}_{\mathfrak{O}_{y_1}}}{|\mathfrak{S}_{y_1}|}\ast \ov{\chi}_{\mathfrak{O}_{(x_1,v_1)}})\ast \Phi(\frac{\ov{\chi}_{\mathfrak{O}_{y_2}}}{|\mathfrak{S}_{y_2}|}\ast \ov{\chi}_{\mathfrak{O}_{(x_2,v_2)}})
		\end{align*}
		by noting that $(\alpha,\beta)_Q=(\alpha,\beta)_{Q'}$ of any $\alpha,\beta\in\Z^\I$.  The proof is completed.
	\end{proof}
	
	\begin{corollary} 
		\label{cor:FT-ihall}
		Let $(Q,\btau)$ be a Dynkin $\imath$quiver.  The Fourier transform $\Phi_{Q',Q}:\widehat{\ch}(\bfk Q,\btau)\rightarrow \widehat{\ch}(\bfk Q',\btau)$ induces an algebra isomorphism  
		$$\Phi_{Q',Q}:\widetilde{\ch}(\bfk Q,\btau)\longrightarrow \widetilde{\ch}(\bfk Q',\btau),$$
		which is called the Fourier transform of $\imath$Hall algebras.
	\end{corollary}
	
	Note that $\Phi_{Q',Q}([S_i])=[S'_i]$, and $\Phi_{Q',Q}([\K_i])=[\K'_i]$ for any $i\in\I$ by definition.
	Since $Q$ is a Dynkin quiver, by \cite[Proposition 5.5]{LW19}, we know $\widehat{\ch}(\bfk Q,\btau)$ (and also $\widetilde{\ch}(\bfk Q,\btau)$)
	is generated by $[S_i],[\K_i]$ ($i\in\I$). So we can see that $\Phi_{Q',Q}:\widetilde{\ch}(\bfk Q,\btau)\rightarrow \widetilde{\ch}(\bfk Q',\btau)$ and $\Phi_{Q',Q}:\widehat{\ch}(\bfk Q,\btau)\rightarrow \widehat{\ch}(\bfk Q',\btau)$
	are well defined when considering these Hall algebras over any field containing $\Z[\sqq,\sqq^{-1}]$, especially  the field $\Q(\sqq^{1/2})$. Moreover, for integral forms of $\imath$Hall algebras, we have the following corollary.
	
	\begin{corollary} 
		\label{cor:FT-ihall-integral}
		Let $(Q,\btau)$ be a Dynkin $\imath$quiver.  The Fourier transform $\Phi_{Q',Q}:\widehat{\ch}(\bfk Q,\btau)\rightarrow \widehat{\ch}(\bfk Q',\btau)$ induces algebra isomorphisms  
		$$\Phi_{Q',Q}:\widehat{\ch}(\bfk Q,\btau)_{\cz}\longrightarrow \widehat{\ch}(\bfk Q',\btau)_{\cz},\quad\Phi_{Q',Q}:\widetilde{\ch}(\bfk Q,\btau)_{\cz}\longrightarrow \widetilde{\ch}(\bfk Q',\btau)_{\cz}.$$
		%which is called the Fourier transform of $\imath$Hall algebras.
	\end{corollary}
	
	\begin{proof}
		By Lemma \ref{lem:FT-char}, we know that 
		$\Phi_{Q',Q}([M])=\sum_{[N]\in\Iso(\mod(\bfk Q'))}\sqq^{a_N}\xi_N[N]$ for any $M\in\mod(\bfk Q)$,
		where $\xi_N$ is a sum of some $q$-roots of unity (in $\C$). We have shown that $\xi_N\in\Q$. Then $\xi_N\in\Z$ by noting that $q$-roots of unity are algebraic integers. So $\Phi_{Q',Q}([M])\in\widetilde{\ch}(\bfk Q',\btau)_\cz$, and then $\Phi_{Q',Q}(\widetilde{\ch}(\bfk Q,\btau)_\cz)\subseteq \widetilde{\ch}(\bfk Q',\btau)_\cz$. Using the inverse $\ov{\Phi}_{Q,Q'}$
		of $\Phi_{Q',Q}$, the desired result follows.
	\end{proof}
	
	\begin{remark}
		For arbitrary quiver $Q$ (not necessarily of Dynkin type), if the set $\ce$ satisfies that any $h\in \ce$ does not lie on any cycle (not necessarily oriented), then the result of Lemma \ref{lem:iso} remains valid, so the Fourier transform $\Phi_{Q',Q}:\widehat{\ch}(\bfk Q,\btau)\rightarrow \widehat{\ch}(\bfk Q',\btau)$ is an isomorphism of $\C$-algebras. For example, this holds for any tree quivers.
	\end{remark}
	
	\begin{remark}
		By Example \ref{ex:diagquiver}, we can view the quiver algebra $\Lambda$ in \eqref{eq:La} to be the $\imath$quiver algebra of $(Q^{\rm dbl}=Q\sqcup Q^\diamond,\swa)$.
		Note that $\mod(\bfk Q)\times \mod(\bfk Q^\diamond)\subseteq \mod(\Lambda)$.  It is proved in \cite[Theorem 4.9]{LP21} that the $\imath$Hall algebra $\cs\cd\widetilde{\ch}(\Lambda)$ is isomorphic to the Drinfeld double of the extended Hall algebra defined in \cite{X97}. 
		
		Let $\ce$ be a subset of $Q$, and $Q'$ be the quiver constructed from $Q$ by reversing the arrows in $\ce$. Let $'\Lambda$ be the algebra defined in \eqref{eq:La} for $Q'$. 
		The Fourier transform gives an isomorphism of algebras from $\mathcal{SD}\widetilde{\mathcal{H}}(\Lambda)$ to $\mathcal{SD}\widetilde{\mathcal{H}}('\Lambda)$, provided that we use the conjugate character $\ov{\phi}$ to define the action on $\mod(\bfk Q^\diamond)\subset \mod(\Lambda)$ (cf. \cite[eq. (8.10)]{SV99}). This can also be observed from the proof of Theorem~\ref{thm:FThall}.
	\end{remark}

	%%%%%%%%%%
	\subsection{Fourier transform of dual canonical bases}
	
	In the sequel \cite{LP26} of this paper, we shall prove that the Fourier transform $\Phi_{Q',Q}$ preserves the dual canonical bases:
	\begin{proposition}
		[\cite{LP26}]\label{prop:dCB-Fourier}
		For a Dynkin $\imath$quiver $(Q,\btau)$, we have that the Fourier transform $$\Phi_{Q',Q}:\widehat{\ch}(\bfk Q,\btau)_{\cz}\longrightarrow \widehat{\ch}(\bfk Q',\btau)_{\cz},\quad\Phi_{Q',Q}:\widetilde{\ch}(\bfk Q,\btau)_{\cz}\longrightarrow \widetilde{\ch}(\bfk Q',\btau)_{\cz}.$$
		preserves the dual canonical bases of  of $\widehat{\ch}(\bfk Q,\btau)_{\cz}$ and $\widetilde{\ch}(\bfk Q,\btau)_{\cz}$.
	\end{proposition}
	
	From Proposition \ref{prop:dCB-Fourier}, we know that the dual canonical basis $\{\K_\alpha\diamond \mathfrak{L}_\lambda \mid \lambda\in\fp,\alpha\in\Z^\I\}$ of $\widetilde{\ch}(Q,\varrho)$ (and also $\widetilde{\ch}(Q,\varrho)_\cz$) does not depend on the orientation of $Q$. From them, we  can also obtain the similar statement for $\widehat{\ch}(Q,\btau)$ (and $\widehat{\ch}(Q,\btau)_\cz$), so we make the following definition. 
	
	\begin{definition}[Dual canonical basis of $\imath$quantum groups]
		\label{def:dual CB}
		The dual canonical basis for $\widetilde{\ch}(Q,\varrho)_\cz$ is transferred to a basis for $\tUi_\cz$ via the isomorphism in Lemma~\ref{lem:Hall-iQG}, which are called the {\em dual canonical basis} for $\tUi$.
	\end{definition}
	
	In particular, for quantum groups $\tU$, we have the following. 
	\begin{definition}[Dual canonical basis of quantum groups]
		\label{def:dual CB}
		The dual canonical basis for $\widetilde{\ch}(Q^{\dbl},\swa)_\cz$ is transferred to a basis for $\tU_\cz$ via the isomorphism in Lemma \ref{lem:bridgeland}, which are called the {\em dual canonical basis} for $\tU$.
	\end{definition}
	
	By combining with Corollary \ref{iHA reflection functor of dCB}, we have  the following result.
	
	\begin{theorem}
		\label{iQG dCB braid invariant}
		The dual canonical basis of $\tUi$ is preserved by the action of the  braid group operator $\tTT_i$ ($i\in\I$).
	\end{theorem}
	
	\begin{proof}
		For $i\in Q_0$, to prove that $\tTT_i$ preserves the dual canonical basis, we choose an orientation $(Q,\varrho)$ such that $i$ is a sink. Then there is a commutative diagram
		\[
		\begin{tikzcd}
			\tUi \ar[r,"\tTT_i"]\ar[d,swap,"\widetilde{\psi}_Q"]& \tUi \ar[d,"\widetilde{\psi}_{r_i Q}"]\\
			\tMHg \ar[r,"\Gamma_i"] &  \widetilde{\ch}(r_iQ,\btau)
		\end{tikzcd}
		\]
		By Corollary~\ref{iHA reflection functor of dCB}, the map $\Gamma_i$ preserves dual canonical bases of $\imath$Hall algebras. The assertion for $\tTT_i$ now follows from Proposition~\ref{prop:dCB-Fourier}.
	\end{proof}
	
	By Theorem \ref{iQG dCB braid invariant}, the dual canonical basis of $\tU$ is preserved by the action of the  braid group operator $T_i$ ($i\in\I$); compare with \cite{Lus96} of the invariance of (part of) the canonical basis of $\U^+$ under Lusztig's braid group actions.


\begin{thebibliography}{DDPW01}
		
		%\bibitem[As11]{As} H. Asashiba,
		%{\em A generalization of Gabriel's Galois covering functors and derived equivalences}, J. Algebra {\bf 334} (2011), 109--149.
		%\bibitem[Ach21]{Ach21} P. Achar, Perverse sheaves and applications to representation theory, vol. {\bf258}, American Mathematical Society, 2021.
		
		%\bibitem[AB89]{AB89} M. Auslander and R. Buchweitz,
		%{\em The homological theory of maximal Cohen-Macaulay approximation}, Colloque en l'honneur de Pierre Samuel (Orsay, 1987). M\'{e}m. Soc. Math. France (N.S.) {\bf 38} (1989), 5--37.
		
		%\bibitem[BKLW18]{BKLW} H. Bao, J. Kujawa, Y. Li and W. Wang,
		%{\it  Geometric Schur duality of classical type}, Transform. Groups {\bf 23} (2018), 329--389.
		%\href{https://arxiv.org/abs/1404.4000}{arXiv:1404.4000v3}
		
		
		
		%\bibitem[BSWW18]{BSWW} H. Bao, P. Shan, W. Wang and B. Webster,
		%{\em Categorification of quantum symmetric pairs I}, Quantum Topology {\bf 9} (2018), 643--714. %,  \href{https://arxiv.org/abs/1605.03780}{arXiv:1605.03780v2}
		
		\bibitem[BW18a]{BW18}
		H. Bao and W. Wang,
		{\em A new approach to Kazhdan-Lusztig theory of type $B$ via quantum symmetric pairs},  Ast\'erisque {\bf 402}, 2018, vii+134pp. %, \href{http://arxiv.org/abs/1310.0103}{arXiv:1310.0103v2}
		
		\bibitem[BW18b]{BW18b} H. Bao and W. Wang,
		{\em Canonical bases arising from quantum symmetric pairs}, Invent. Math. {\bf 213} (2018), 1099--1177. %\href{https://arxiv.org/abs/1610.09271}{arXiv:1610.09271v2}
		
		%\bibitem[BBD82]{BBD} 
		%A.A. Beilinson, J. Bernstein, P. Deligne, \emph{Faisceaux pervers}, Ast\'{e}risque {\bf100}, 1982.
		
		%\bibitem[BLM90]{BLM90}
		%A. Beilinson, G. Lusztig and R. MacPherson,
		%         {\em A geometric setting for the quantum deformation of $GL_n$},
		%Duke Math. J. {\bf 61} (1990), 655--677.
		
		\bibitem[BG17]{BG17} A. Berenstein and J. Greenstein, {\em Double canonical bases}, Adv. Math. {\bf 316} (2017), 381--468.
		
		\bibitem[BZ14]{BZ14} A. Berenstein and A. Zelevinsky, {\em Triangular Bases in Quantum Cluster Algebras}, IMRN {\bf2014} (2014), no. 6, 1651--1688.
		
		%\bibitem[Bra03]{Bra03} T. Braden, {\em Hyperbolic localization of intersection cohomology}, Transform. Groups {\bf8} (2003), no. 3, 209--216.
		
		\bibitem[Br13]{Br13} T. Bridgeland, {\em Quantum groups via Hall algebras of complexes}, Ann. Math. {\bf 177} (2013), 739--759.
		
		
		
		
		
		%\bibitem[CM06]{CM06} C. Cibils, E. Marcos, {\em Skew category, Galois covering and smash product of a $k$-category}, Proc. Amer. Math. Soc. {\bf134}
		%(2006), 39--50.
		
		\bibitem[DDPW08]{DDPW}
		B.~Deng, J.~Du, B.~Parshall and J.~Wang,
		Finite dimensional algebras and quantum groups,
		Mathematical Surveys and Monographs {\bf 150}.
		AMS, Providence, RI, 2008.
		
		%\bibitem[dCM02]{dM02} M. A. A. de Cataldo and L. Migliorini, {\em The Hard Lefschetz Theorem and the topology of semismall maps}, Ann. Sci. \'Ecole Norm. Sup. (4) {\bf35} (2002), no. 5, 759–772.
		
		\bibitem[EJ00]{EJ} E.E. Enochs and O.M.G. Jenda,
		Relative homological algebra, De Gruyter Exp. Math. {\bf 30}, Walter de Gruyter Co., 2000.
		
		%\bibitem[FLX23]{FLX23} J. Fang, Y. Lan and J. Xiao, {\em Sheaf realization of Bridgeland's Hall algebra of Dynkin type
		%}, \href{arXiv:2303.0499}{https://arxiv.org/abs/2303.04993}
		
		%\bibitem[FK88]{FK88} E. Freitag and R. Kiehl, Etale Cohomology and the Weil Conjecture, Springer, 1988.
		
		
		%\bibitem[Fu17]{Fu17} R. Fujita,
		%{\em Affine highest weight categories and quantum affine Schur-Weyl duality of Dynkin quiver types}, 	Represent. Theory {\bf26} (2022), 211--263.  %\href{https://arxiv.org/abs/1710.11288}{arXiv:1710.11288}
		
		%\bibitem[Ga79]{Ga79} P. Gabriel, {\em Auslander-Reiten sequences and representation-finite algebras}, Representation theory, I (Proc.
		%Workshop, Carleton Univ., Ottawa, Ont., 1979), Springer, Berlin, 1980,  1–-71.
		
		
		%\bibitem[Ga81]{Ga} P. Gabriel, {\em The universal cover of a representation finite algebra}, in: Representation of Algebras, in: Lecture Notes in Math., vol. {\bf903} (1981), 65--105.
		
		
		%\bibitem[GLS12]{GLS12}  C. Geiss, B. Leclerc, and J. Schr\"{o}er, {\em Generic bases for cluster algebras and the Chamber Ansatz}, J.
		%Amer. Math. Soc. {\bf25} (2012), no. 1, 21–76.
		
		\bibitem[GLS13]{GLS13}  C. Geiss, B. Leclerc, and J. Schr\"{o}er, {\em Cluster structures on quantum coordinate rings}, Selecta Math. (N.S.) {\bf19} (2013), 337--397.
		
		%\bibitem[GY21]{GY21} K.R. Goodearl and M.T. Yakimov, {\em  Integral quantum cluster structures}, Duke Math. J. {\bf170} (2021), no. 6, 1137--1200.
		%\bibitem[Gor13]{Gor13} M. Gorsky,
		%{\em Semi-derived Hall algebras and tilting invariance of Bridgeland-Hall algebras}, \href{https://arxiv.org/abs/1303.5879}{arXiv:1303.5879v2}
		
		\bibitem[Gor18]{Gor18} M. Gorsky,
		{\em Semi-derived and derived Hall algebras for stable categories}, IMRN, Vol. {\bf 2018}, 138--159. %\href{https://arxiv.org/abs/1409.6798}{arXiv:1409.6798}.
		
		\bibitem[Gr95]{Gr95} J.A. Green,
		{\em Hall algebras, hereditary algebras and quantum groups}, Invent. Math. {\bf 120} (1995), 361--377.
		
		%\bibitem[Ha87]{Ha1} D. Happel, {\em On the deived category of a finite-dimesional algebra}, Comment. Math. Helv. {\bf62} (3)(1987), 339--389.
		
		%\bibitem[Ha88]{Ha2} D. Happel, Triangulated Categories in the Representation Theory of Finite Dimensional Algebras. London Math. Soc. Lecture Notes Ser. {\bf119}, Cambridge Univ. Press, Cambridge, 1988.
		
		\bibitem[Ha91]{Ha3} D. Happel,
		{\em On Gorenstein algebras}, In: Progress in Math. {\bf 95}, Birkh\"{a}user Verlag, Basel, 1991, 389--404.
		
		\bibitem[HL15]{HL15} D. Hernandez and B. Leclerc, {\em Quantum Grothendieck rings and derived Hall algebras}, J. Reine Angew. Math. {\bf701} (2015), 77--126.
		
		%\bibitem[I05]{I05} O. Iyama, $\tau$-categories. I: Ladders, Algebr. Represent. Theory {\bf8} (2205), no.3, 297--321.
		
		%\bibitem[KKKO18]{KKKO18} S.-J. Kang, M. Kashiwara, M. Kim, and S.-j. Oh, {\em Monoidal categorification of cluster algebras}, J. Amer. Math. Soc. 31 (2018), 349--426.
		
		\bibitem[Ka91]{Ka91} M. Kashiwara,
		{\em On crystal bases of the $Q$-analogue of universal enveloping algebras}, Duke Math.~J.~{\bf 63} (1991), 456--516.
		
		%\bibitem[Ka14]{Ka14} S. Kato, {\em Poincar\'{e}-Birkhoff-Witt bases and Khovanov-Lauda-Rouquier algebras}, Duke Math. J. {\bf163} (2014), 619--663.
		
		%\bibitem[Ke90]{Ke1} B. Keller, {\em Chain complexes and stable categories}, Manus. Math. {\bf67} (1990), 379--417.
		
		\bibitem[Ke05]{Ke05} B. Keller,
		{\em On triangulated orbit categories}, Doc. Math. {\bf 10} (2005), 551--581.
		
		\bibitem[KS16]{KS16} B. Keller and S. Scherotzke, {\em Graded quiver varieties and derived categories}, J. Reine Angew. Math. {\bf713} (2016), 85--127.
		
		%\bibitem[KKK15]{KKK15} S.-J. Kang, M. Kashiwara, M. Kim,
		%{\em Symmetric quiver Hecke algebras and $R$-matrices of quantum affine algebras, II}, Duke Math.~J.~{\bf 164} (2015), 1549--1602.
		
		%\bibitem[KL09]{KL09}
		%M. Khovanov and A. Lauda,
		%\emph{A diagrammatic approach to  categorification of quantum groups. {I}}, Represent. Theory \textbf{13}  (2009), 309--347.
		
		%\bibitem[KL10]{KL10}
		%M. Khovanov and A. Lauda,
		%\emph{A categorification of quantum {${\rm sl}(n)$}}, Quantum Topology
		% \textbf{1} (2010),  1--92.  %\MR{2628852 (2011g:17028)}
		
		%\bibitem[KW01]{RW01} R. Kiehl and R. Weissauer, Weil conjectures, perverse sheaves and l'adic Fourier transform,  Ergebnisse 
		%der Mathematik und Ihrer Grenzgebiete. 3. Folge. A Series of Modern Surveys in Mathematics 
		%[Results in Mathematics and Related Areas. 3rd Series. A Series of Modern Surveys in Mathematics], 
		%vol. {\bf42}, Springer-Verlag, Berlin, 2001.
		
		
		\bibitem[Ko14]{Ko14} S. Kolb, {\em Quantum symmetric Kac-Moody pairs}, Adv. Math. {\bf267} (2014), 395--469.
		
		\bibitem[KP11]{KP11} S. Kolb and J. Pellegrini,
		{\em Braid group actions on coideal subalgebras of quantized enveloping algebras}, J. Algebra {\bf 336} (2011), 395--416.
		
		
		%Publications Math\'ematiques de L’Institut des Hautes Scientifiques {\bf 65}, 131--210 (1987).
		
		
		
		
		
		%\bibitem[Lau87]{Lau87} G. Laumon,  {\em Transformation de Fourier, constantes d’\'equations fonctionnelles et conjecture de Weil}, Inst. Hautes \'Etudes Sci. Publ. Math. {\bf65} (1987), 131--210. 
		
		\bibitem[LeP13]{LeP13} B. Leclerc and P.-G. Plamondon, {\em Nakajima varieties and repetitive algebras}, Publ.  RIMS {\bf49} (2013), 531--561.
		
		%\bibitem[LeBP90]{LBP} L. Le Bruyn and C. Procesi, {\em Semisimple representations of quivers}, Trans.
		%Amer. Math. Soc. {\bf 317} (1990), 585--598.
		
		\bibitem[Let99]{Let99} G. Letzter, {\em Symmetric pairs for quantized enveloping algebras}, J. Algebra {\bf220} (1999), 729--767.
		
		%\bibitem[LiW18]{LiW18} Y. Li and W. Wang,
		%{\em Positivity vs negativity of canonical bases}, %Proceedings for Lusztig's 70th birthday conference,
		%Bull. Inst. Math. Acad. Sin. (N.S.) {\bf13} (2) (2018), 143--198. %, %{\bf 13} (2018), 143--198.
		%\href{arXiv:1501.00688v4}{https://arxiv.org/abs/1501.00688v4}.
		
		\bibitem[LP26]{LP26} M. Lu and X. Pan, {\em Dual canonical bases of quantum groups and $\imath$quantum groups II: geometry}, preparing.
		
		
		\bibitem[LP24]{LP24} J. Lin and L. Peng, {\em Semi-derived Ringel-Hall algebras and Hall algebras of odd-periodic relative derived categories}, Sci. China Math. {\bf67} (2024), no.8, 1735--1760.
		
		
		
		\bibitem[LP21]{LP21} M. Lu and L. Peng,
		{\em Semi-derived Ringel-Hall algebras and Drinfeld double}, Adv. Math. {\bf 383} (2021), 107668.
		
		\bibitem[LR24]{LR24} M. Lu and S. Ruan,  {\em $\imath$Hall algebras of weighted projective lines and quantum symmetric pairs III: quasi-split type}, \href{https://arxiv.org/abs/2411.13078}{arXiv:2411.13078}
		
		
		%\bibitem[LW19b]{LW19b} M. Lu and W. Wang, {\em Hall algebras and quantum symmetric pairs II: Reflection functors}, \href{http://export.arxiv.org/abs/1904.01621}{arXiv:1901.01621}
		
		%\bibitem[LW19e]{LW19e} M. Lu and W. Wang, {\em Hall algebras and quantum symmetric pairs V: Bases}, in preparation.
		
		%\bibitem[Lu19]{Lu19} M. Lu,
		%{\em Modified Ringel-Hall algebras of 1-Gorenstein algebras},
		%Appendix A to \cite{LW19}, \href{http://arxiv.org/abs/1901.11446}{arXiv:1901.11446}
		
		\bibitem[LW21a]{LW21a} M. Lu and W. Wang, {\em Hall algebras and quantum symmetric pairs II: reflection functors},  Commun. Math. Phys.  {\bf 381} (2021), 799--855. %, \href{http://arxiv.org/abs/1904.01621}{arXiv:1904.01621}
		
		
		\bibitem[LW21b]{LW21b} M. Lu and W. Wang,
		{\em Hall algebras and quantum symmetric pairs III: Quiver varieties}, Adv. Math. {\bf393} (2021), 108071.
		
		
		\bibitem[LW22a]{LW19} M. Lu and W. Wang, {\em Hall algebras and quantum symmetric pairs I: Foundations}, Proc. London Math. Soc. {\bf 124} (1):  1--82.
		
		
		\bibitem[LW23]{LW20} M. Lu and W. Wang, {\em Hall algebras and quantum symmetric pairs of Kac-Moody type}, Adv. Math. {\bf430} (2023), 109215. 
		
		
		
		\bibitem[Lus90a]{Lus90a} G.~Lusztig,
		\textit{Finite-dimensional Hopf algebras arising from quantized universal enveloping algebra}, J. Amer. Math. Soc. {\bf 3} (1990), 257--296.
		
		\bibitem[Lus90b]{Lus90} G. Lusztig, {\em Canonical bases arising from quantized enveloping algebras}, J. Amer. Math.
		Soc. {\bf 3} (1990),  447--498.
		
		\bibitem[Lus91]{Lus91} G. Lusztig,  {\em Quivers, perverse sheaves, and quantized enveloping algebras}, J. Amer.
		Math. Soc. {\bf4} (2)(1991), 365--421, 1991.
		
		\bibitem[Lus93]{Lus93} G. Lusztig, Introduction to Quantum Groups, Birkh\"{a}user, Boston, 1993.
		
		\bibitem[Lus96]{Lus96} G. Lusztig, {\em Braid group action and canonical bases}, Adv. Math. {\bf122} (1996), 237--261.
		
		\bibitem[Lus98]{Lus98} G. Lusztig, {\em Canonical bases and Hall algebras}, in: Representation Theories and Algebraic Geometry, 
		Springer, 1998, 365--399.
		
		\bibitem[Na01]{Na01} H. Nakajima, {\em Quiver varieties and finite-dimensional representations of quantum affine algebras}, J. Amer. Math. Soc. {\bf14} (2001), 145--238 (electronic).
		
		\bibitem[Na04]{Na04} H. Nakajima, {\em Quiver varieties and $t$-analogs of $q$-characters of quantum affine algebras}, Ann. Math. {\bf160} (2004), 1057--1097.
		
		%\bibitem[NV04]{NV04} C. N\v{a}st\v{a}sescu and F. Van Oystaeyen, Methods of Graded Rings, Lecture Notes in
		%Mathematics, {\bf1836}, Springer-Verlag, Berlin, 2004.
		
		%\bibitem[Q14]{Qin14} F. Qin, 
		%{\em $t$-Analog of $q$-Characters, bases of quantum cluster algebras, and a correction technique}, IMRN, Vol. {\bf 2014}, No. 22, 6175--6232.
		
		\bibitem[Q16]{Qin} F. Qin, {\em Quantum groups via cyclic quiver varieties I}, Compos. Math. {\bf152} (2016), 299--326.
		
		%\bibitem[Q21]{Qin21} F. Qin, {\em Dual canonical bases and quantum cluster algebras}, \href{https://arxiv.org/abs/2003.13674v3}{arxiv:2003.13674}
		
		
		%\bibitem[Qui73]{Q} D. Quillen, {\em Higher algebraic K-theory, I}, Springer Lecture Notes in Mathematics {\bf 341} (1973), 85--147.
		
		%\bibitem[Rie86]{Rie86} C. Riedtmann, {\em Degenerations for representations of quivers with relations}, Ann. scient. \'{E}c. Norm. Sup. {\bf 19} (1986), 275--301.
		
		\bibitem[Rin90]{Rin90} C. Ringel, {\em Hall algebras and quantum groups}, Invent. Math. {\bf101} (1990), 583--591.
		
		\bibitem[Rin96]{Rin3} C.M. Ringel,
		{\em PBW-bases of quantum groups}, J. reine angrew. Math. {\bf 470} (1996), 51--88.
		
		%\bibitem[Rin96b]{Rin4}
		%C.M. Ringel, {\em Green's theorem on Hall algebras} in Representations of Algebras and Related Topics, CMS Conference Proceedings {\bf 19}, Providennce, 1996,185--245.
		
		%\bibitem[R08]{R08}
		%R.~Rouquier,  {\em 2-Kac-Moody Algebras},  \href{https://arxiv.org/abs/0812.5023}{arXiv:0812.5023}.
		
		
		
		%\bibitem[R12]{R12}
		%\bysame,
		%{\em  Quiver Hecke algebras and 2-Lie algebras},  Algebra Colloq. \textbf{19} (2012),  359--410.
		
		%\bibitem[Sch17]{Sch17} S. Scherotzke, {\em Desingularization of Quiver Grassmannians via Nakajima Categories}, Algebr. Represent. Theor. {\bf 20} (2017), 231--243.
		
		\bibitem[Sch19]{Sch19} S. Scherotzke, {\em Generalized quiver varieties and triangulated categories}, Math. Z. {\bf 292} (2019), 1453--1478. %,  \href{https://arxiv.org/abs/1405.4729}{arXiv:1405.4729v6}
		
		%\bibitem[SS16]{SS16}  S. Scherotzke and N. Sibilla, {\em Quiver varieties and Hall algebras}, Proc. London Math. Soc. {\bf 112} (2016), 1002--1018.
		
		
		\bibitem[SV99]{SV99} B. Sevenhant, M. Van den Bergh, {\em On the double of the Hall algebra of a quiver}, J. Algebra {\bf 221} (1999), 135--160.
		
		
		%\bibitem[Sh22]{Sh22} L. Shen, {\em Duals of Semisimple Poisson–Lie Groups and Cluster Theory of Moduli Spaces of G-local
		% Systems},  
		% IMRN {\bf 18} (2022), 14295--14318.
		
		%\bibitem[Sh22]{Sh22} L. Shen, {\em Cluster nature of quantum groups}, \href{https://arxiv.org/abs/2209.06258}{arXiv:2209.06258}
		
		%\bibitem[So24]{S24} J. Song, {\em Cluster realisations of $\imath$quantum groups of type AI}, Proc. London Math. Soc. {\bf 129} (4) (2024).
		
		%\bibitem[VV03]{VV} M. Varagnolo and E. Vasserot, {\em Perverse sheaves and quantum Grothendieck rings}, Studies in memory
		%of Issai Schur (Chevaleret/Rehovot, 2000), Progress in Mathematics {\bf210} (Birkhauser Boston, Boston, MA,
		%2003), 345--365.
		
		\bibitem[WZ23]{WZ23} W. Wang and W. Zhang,
		{\em An intrinsic approach to relative braid group symmetries on $\imath$quantum groups}, Proc. London Math. Soc. \textbf{127} (2023), 1338--1423.
		
		\bibitem[X97]{X97} J. Xiao, {\em Drinfeld double and Ringel-Green theory of Hall algebras}, J. Algebra {\bf 190} (1997), 100--144.
		
		%\bibitem[XY01]{XY01} J. Xiao, S. Yang, {\em BGP-reflection functors and Lusztig’s symmetries: a Ringel–Hall algebra approach 
		%to quantum groups}, J. Algebra {\bf241} (1) (2001), 204--246.
		
		%\bibitem[XZ05]{XZ05} J. Xiao and B. Zhu, {\em Locally finite triangulated categories}, J. Algebra {\bf 290} (2005), 473--490.
	\end{thebibliography}
\end{document}